\documentclass[10pt]{amsart}
\usepackage{amssymb}
\usepackage{amssymb,amsthm,amsmath}
\usepackage[numbers,sort&compress]{natbib}
\usepackage{color}
\usepackage{graphicx}
\usepackage{tikz}
\usepackage{bbm}
\usepackage{mathrsfs}
\usepackage{amsfonts}
\usepackage{bm}

\title[long time stability]{Long time stability of KAM tori for the nonlinear Schr\"odinger equation}

\author[X. He]{Xiaolong He}
\address[Xiaolong He]
{Department of mathematics, Hangzhou Normal University, Hangzhou, 311121, China}
 \email{xlhe@hznu.edu.cn}

\author[J. Shi]{Jia Shi}
\address[Jia Shi]
{School of Mathematical Sciences,
Fudan University,
Shanghai 200433, China} \email{15110180007@fudan.edu.cn}

\author[X. Yuan]{Xiaoping Yuan}
\address[Xiaoping Yuan]
{School of Mathematical Sciences,
Fudan University,
Shanghai 200433, China} \email{xpyuan@fudan.edu.cn}

\keywords{KAM tori, long time stability, nonlinear Schr\"odinger equation.}

\theoremstyle{plain}
\newtheorem{thm}{Theorem}[section]
 
 \newtheorem{lem}[thm]{Lemma}
 \newtheorem{prop}[thm]{Proposition}
 \newtheorem{defn}[thm]{Definition}
 \newtheorem{rem}[thm]{Remark}
 
 \numberwithin{equation}{section}

\begin{document}


\begin{abstract}
In this paper, we prove the long time stability of KAM tori for the nonlinear Schr\"odinger equation on
the  torus with arbitrary
dimensions.
\end{abstract}
\maketitle

\section{Introduction}

In this paper, we consider the long time stability of the invariant tori for the $d$-dimensional nonlinear Schr\"odinger (NLS) equation
\begin{equation}\label{1}
-\mathrm{i}\dot{u}=-\Delta u+V(x)\ast u+\varepsilon\frac{\partial F}{\partial \bar{u}}(|u|^{2}), \quad u=u(t,x)
\end{equation}
under the periodic boundary condition $x\in\mathbb{T}^{d}, d\geq 1$.
The convolution function
$V: \mathbb{T}^{d}\rightarrow \mathbb{C}$ is analytic and the Fourier coefficient $\hat{V}(a)$
takes real value,
when expanding $V$ into Fourier series
$V(x)=\sum_{a\in\mathbb{Z}^{d}}\hat{V}(a)e^{\mathrm{i}\langle a,x\rangle}$.
The nonlinearity $F$ is real analytic.

The NLS equation \eqref{1} is a Hamiltonian PDE.
The  KAM theory is a well-known approach   to
establish the existence of  the invariant tori for Hamiltonian PDEs. The invariant tori so constructed
are often referred to as  the \emph{KAM tori}.
 The ``KAM for  PDE'' theory started
in late 1980's and originally applied  to the  one  spatial dimensional PDEs,
which is now  well understood. See for example \cite{K, Kuk00, KP,P1,P2,W, LY11, BBP13, BBHM18}
and the references therein.

However, the  KAM theory for space-multidimensional Hamiltonian PDEs is
at its early stage.
The first breakthrough  was made by Bourgain
\cite{B1} on the two dimensional NLS equation, in which he developed
Craig and Wayne's scheme on periodic problems. Using
 new techniques of Green's function estimates
in the spectral theory, Bourgain proved the persistence of invariant tori for space-multidimensional NLS and nonlinear wave (NLW) equations \cite{B2}.
The above mentioned method is now known as the Craig-Wayne-Bourgain (CWB) method.
See also \cite{BB13, BB20, Wang16, Wang19, Wang20} and the references therein.
The classical KAM approach for space-multidimensional Hamiltonian PDEs was developed by
Eliasson and Kuksin  \cite{EK} on NLS equation.
They take  a sequence of symplectic transformations such that the transformed Hamiltonian guarantees the existence of the invariant tori.
Moreover, the KAM approach in \cite{EK} also provides  the reducibility and linear stability of the obtained invariant
tori. See also \cite{EGK,GY,GXY11, PP1,PP2,Y} for the KAM approach on the space-multidimensional  PDEs.

To ensure that the obtained KAM tori can be observed in physics and other real applications,
one has to prove that those KAM tori are stable in some sense, among which
the simplest one is the \emph{linear stability}. Let us recall the definition of the linear stability of the invariant tori.
Consider a nonlinear differential equation
\begin{equation}\label{01}
  \dot{x}=X(x),
\end{equation}
which has an invariant torus $\mathcal{T}$ carrying the quasi-periodic flow $x(t)=x_{0}(t)$.
We say that the invariant torus $\mathcal{T}$ is linearly stable if the equilibrium of the linearized equation
\begin{equation*}
  \dot{y}=D X(x_{0}(t)) y
\end{equation*}
of \eqref{01} along $\mathcal{T}$
is Lyapunov stable. A more general definition is that
$\mathcal{T}$ is linearly stable if all the Lyapunov exponents of $\mathcal{T}$
equal to zero.

The classical KAM tori are linearly stable as a matter of the KAM approach.
To see this, we consider a Hamiltonian perturbation
\begin{equation*}
  H= N+P= N+ P^{\textrm{low}}+ P^{\textrm{high}}
\end{equation*}
of the integrable
part
\begin{equation*}
  N= \langle \omega, y\rangle+ \sum_{j\in \mathbb{Z}^{d}} \Omega_{j} z_{j} \bar{z}_{j},
\end{equation*}
where
\begin{equation*}
\begin{aligned}
P^{\textrm{low}}=& R^{x}+ \langle R^{y}, y\rangle+\langle R^{z}, z\rangle
+\langle R^{\bar{z}}, \bar{z}\rangle
+ \langle R^{zz} z, z\rangle+ \langle R^{z\bar{z}} z, \bar{z}\rangle
+\langle R^{\bar{z}\bar{z}}\bar{z}, \bar{z}\rangle
\end{aligned}
\end{equation*}
and
\begin{equation*}
  P^{\textrm{high}}= O(|y|^{2}+ \|z\| \cdot |y|+ \|z\|^{3}).
\end{equation*}
The classical KAM approach (for $d=1$) aims at taking a sequence of symplectic
transformations to eliminate all  terms in $P^{\textrm{low}}$,
except for the averages $\langle \widehat{R^{y}}(0), y\rangle$
and $\sum_{i=j} \widehat{R^{z\bar{z}}_{ij}}(0) z_{i} \bar{z}_{j}$.
In particular, the quadratic  terms in
$P^{\textrm{low}}$ are reduced to  $\sum_{i=j} \widehat{R^{z\bar{z}}_{ij}}(0) z_{i} \bar{z}_{j}$
of constant coefficients, which can be put into the integrable part $N$ for the next iteration.
In this way, the linearized equation of the obtained KAM tori can be reduced to
\begin{equation}\label{001}
  i \dot{z}= (\Omega+ \varepsilon [R]) z,
\end{equation}
where $\Omega=\textrm{diag}(\Omega_{j}: j\in\mathbb{Z})$
and $[R]= \textrm{diag}((R^{\infty;z\bar{z}}_{ij})^{\wedge}(0): i=j\in \mathbb{Z})$
are diagonal and constant.\
Obviously, the equilibrium $z=0$ of \eqref{001} is Lyapunov stable, and
thus the KAM tori are linearly stable.

Unfortunately, there is a difficulty in extending the classical KAM approach for
$d=1$ to the case of $d>1$. Taking NLS equation for example, the normal frequency satisfies
$\Omega_{j}\sim |j|^{2}, j\in\mathbb{Z}^{d}$ after writing NLS equation into an infinitely
dimensional Hamiltonian system as above. It follows that the normal frequencies may have
unbounded multiplicities since  $\# \{j'\in\mathbb{Z}^{d}: |j'|=|j|\}\sim |j|^{d-1}\rightarrow \infty$ as $|j|\rightarrow \infty$.
This feature leads to serious resonances in solving the homological equations,
which might impede the convergence of the symplectic transformations.
Eliasson-Kuksin \cite{EK} analyzed carefully the separation property of the normal frequencies
and provided  insight on the T\"{o}plitz-Lipschitz property of the Hamiltonian.
Using the ``super Newton iteration" (rather than the usual Newton iteration  in KAM theory for $d=1$), they succeeded in eliminating $P^{\textrm{low}}$, but leaving an infinitely dimensional
block-diagonal and constant matrix in the quadratic term of $z,\bar{z}$ behind.
As a result, they proved that there are plenty of KAM tori for NLS equation with $d>1$,
whose all Lyapunov exponents equal to zero and hence are  linearly stable.
As for NLW equation, the normal frequency $\Omega_{j}=|j|=\sqrt{j_{1}^{2}+\cdots+j_{d}^{2}}$
does not have a good separation property like NLS equation.
Although Bourgain \cite{B2} had applied the CWB method
to prove the existence of KAM tori for
NLW equation with $d>1$, the linear stability of those KAM tori remains  open.
Recently, by modifying the CWB method, the authors \cite{HSSY}  obtained not only
the existence but also the linear stability of the KAM tori for the Hamiltonian system
with finite degrees of freedom.

In many cases, we cannot determine the stability of the nonlinear system  from its linearized equation directly.
Typically,
in linear plane dynamical systems, a center equilibrium can become
a focus after certain perturbation.
There are also examples that a linearly stable model can be triggered
by an initial perturbation to exhibit chaotic dynamics \cite{GG}.
This prompts us to study the nonlinear stability, among which the \emph{long time stability} is
of particular interest in PDEs.
In finitely dimensional Hamiltonian system,  the best result concerning the long
time stability  is the Nekhoroshev estimate \cite{Neh}.
Consider a $n$ degree of freedom Hamiltonian $H= N(y)+ \varepsilon R(y, x)$,
where $(y,x)\in \mathbb{R}^{n}\times \mathbb{T}^{n}$ is the action-angle variable.
Assume the functions $N$ and $R$ are analytic in $(y,x)$ in some open domain. The Nekhoroshev estimate tells us that
the variation of the actions of all orbits remains small over a finite, but exponentially
long  time interval.
More precisely, for sufficiently small $\varepsilon$, one has
\begin{equation}\label{Nek1}
| y(t)- y(0)|\lesssim \varepsilon^{a} \quad \textrm{for}~|t|\lesssim \exp (\varepsilon^{-b}),
\end{equation}
where the constants $a, b$ depend on the degree of the freedom.
In particular, if $N$ is convex, one can get $a=b=\frac{1}{2n}$.
See P\"{o}schel \cite{Pos93}.
Noticing also that the instabilities such as Arnold diffusion \cite{Arn64} may occur with the degree of freedom $n\geq 3$ and transfer of energy may appear in NLS equation \cite{CK}, one should not expect some orbits are stable forever.
Consequently, it is reasonable to apply the Nekhoroshev estimate on the long time stability of orbits
to describe the stability of the Hamiltonian system.

For NLS equation (or generally the Hamiltonian PDEs), the degree of freedom of the Hamiltonian is infinite.
One immediately gets that $a=b=0$ and the Nekhoroshev estimate in \eqref{Nek1} no longer
works. Instead, Bourgain \cite{Bou96-GAFA} suggested  investigating the long time behavior of orbits in the neighborhood
of the equilibrium and relaxed the stable time interval from
$|t|< \exp (-\varepsilon^{b})$ to $|t|< \varepsilon^{-M}$ for large $M$.
From then on, there are lots of literature devoted to the long time stability of the equilibrium
for  Hamiltonian PDEs.
See \cite{Bam03, DS04, BG, FG13, YZ14, BMP20,CLW20,CMW20, BG21}.
We emphasize that in \cite{BG} Bambusi and Gr\'{e}bert
introduced the tame property of the vector field, which simplifies the proof considerably.

In contrast to the equilibrium, the KAM tori are much more complicated
solutions  of NLS equation. It is known that the KAM tori are sup-exponentially long time stable
for the finitely dimensional Hamiltonian system \cite{BFG88}.
For Hamiltonian PDEs, the study of the long time stability of KAM tori is
limited to the case of $d=1$.
For instance, \cite{CLY} and \cite{CGL15} studied the long time stability
of the KAM tori for NLS and NLW equations, respectively.
For cubic defocusing NLS equation on $\mathbb{T}^{2}$, Maspero-Procesi \cite{MP18} studied  the large-time
stability (with the stable time interval  $|t|<\delta^{-2}$)  of small  finite gap solutions, which depend only on one spatial variable.
For $d>2$, as far as we know, there seems no results in this respect.
The main result of the present paper is that the majority of the KAM tori obtained by Eliasson-Kuksin
\cite{EK} are stable in long time $|t|< \delta^{-M}$ for large
$M$. More precisely, we have the following theorem.

\begin{thm}\label{t}
Under the assumptions for equation \eqref{1},
if $\varepsilon>0$ is sufficiently small, then for typical $V$ $($in the sense of measure$)$, the nonlinear Schr\"odinger equation \eqref{1}
possesses  a linearly stable KAM torus $\mathcal{T}=\mathcal{T}_{V}$ in the Sobolev space $H^{p}(\mathbb{T}^{d})$. Moreover,
letting $ M \approx\varepsilon^{-\frac{2}{3}}$ and $p\geq80(4d)^{4d}(M+7)^{4}+1$, there exists
a small $\delta_{0}$ depending on $p, M$ and $\emph{\textrm{dim}}~\mathcal{T}$ such that
for any $0<\delta<\delta_{0}$ and any solution $u(t,x)$ of \eqref{1} with the initial datum
$u(0,\cdot)$ satisfying
\begin{equation*}
  d_{H^{p}(\mathbb{T}^{d})}(u(0,\cdot), \mathcal{T}):=\inf_{w\in \mathcal{T}}
  \|u(0,\cdot)- w\|_{H^{p}(\mathbb{T}^{d})}<\delta,
\end{equation*}
we have
\begin{equation*}
   d_{H^{p}(\mathbb{T}^{d})}(u(t,\cdot), \mathcal{T})<2 \delta,
   \quad \forall~ |t|< \delta^{-M}.
\end{equation*}
In other words, the KAM tori for the nonlinear Schr\"odinger equation \eqref{1} are stable
in long time.
\end{thm}

Theorem \ref{t} consists of two results (Theorem \ref{t1} and Theorem
\ref{t2}) after writing  NLS equation as the infinitely dimensional Hamiltonian system.
By excluding some parameters, we establish the KAM theorem (see Theorem \ref{t1}) to
guarantee the existence and linear stability of the KAM tori,
which have already been obtained in \cite{EK}. However, to
study the long time stability of the KAM tori, we modify the proof
in \cite{EK} by taking  Kolomogorov's iterative scheme such that
the transformed Hamiltonian after the KAM iteration is still defined
on an open domain. By the further parameter exclusion, we show
the long time stability of the majority of the obtained KAM tori (see Theorem \ref{t2}), by establishing
the partial normal form of the transformed Hamiltonian.

The proof of Theorem \ref{t} is given at the end of section \ref{sect 2}. We clarify the main ideas in the proof of the theorem.
\begin{enumerate}
  \item[i)] Write \eqref{1} as an infinitely dimensional Hamiltonian system
  \begin{equation*}
    H=\sum_{a\in \mathcal{A}}\omega_{a}r_{a}+\frac{1}{2}\sum_{a\in\mathcal{L}}\Omega_{a}(\xi_{a}^{2}+\eta_{a}^{2})+f
  \end{equation*}
  and split
  \begin{equation*}
    f=f^{\textrm{low}}+f^{\textrm{high}}.
  \end{equation*}
   The KAM approach developed in \cite{EK} aims at  eliminating $f^{\textrm{low}}$,
  but leaving some resonant terms. In this way, the transformed Hamiltonian $H_{\infty}=H\circ \psi_{\infty}
  $ takes the form of
  \begin{equation*}
    H_{\infty}=\langle \omega',r\rangle+\frac{1}{2}\langle \zeta,(\Omega+Q)\zeta\rangle+f^{\textrm{high}}_{\infty}.
  \end{equation*}
  One sees that $r=0, \zeta=0$ is the KAM torus for $H_{\infty}$.

   Recall that the domain  $D(\mu_{j},\sigma_{j})$ of the  symplectic transformation $\psi_{j}$
  in \cite{EK}
  degenerates into a singleton since $\mu_{j}\rightarrow 0, \sigma_{j}\rightarrow 0$
  as $j\rightarrow \infty$. Since $\psi_{j}$ is  quadratic in $r$ and $\zeta$,
  one surely can extend the domain of $\psi_{j}$ to the initial domain $D(\mu_{0},\sigma_{0})$.
  In addition, one can  take $D(\mu_{0},\sigma_{0})$ as the domain of the vector
  field for the  \emph{linearized equation} of the Hamiltonian $H_{\infty}$, based on which
  one is able to study the Lyapunov stability of the equilibrium on $D(\mu_{0}, \sigma_{0})$.
  Along this line, Eliasson-Kuksin \cite{EK} developed a powerful KAM approach for NLS equation
  with $d>1$ to show not only the existence of the KAM tori, but also their linear stability.

  However, when studying the long time stability of those KAM tori, we have to take the domain of the high order term $f^{\textrm{high}}_{\infty}$ into consideration. Since
  the domain $\cap_{j=1}^{\infty} D(\mu_{j},\sigma_{j})$ of $\psi_{\infty}$ is a singleton, we can only define $f^{\textrm{high}}_{\infty}$ on $D(0,0)$, on which
  the KAM tori are indeed constructed for the original NLS equation.
  We emphasize that the domain of $f^{\textrm{high}}_{\infty}$
  usually cannot be extended to $D(\mu_{0},\sigma_{0})$
  since $f^{\textrm{high}}_{\infty}$ is not a polynomial function.
  For that reason, we introduce Kolmogorov's iterative scheme in the framework
  of \cite{EK} by modifying the homological equations such that
  $D(\mu_{j}, \sigma_{j})\supset D(\frac{\mu_{0}}{2}, \frac{\sigma_{0}}{2})$
  for all $j$. Then we can define $f^{\textrm{high}}_{\infty}$
  on an open set $D(\frac{\mu_{0}}{2}, \frac{\sigma_{0}}{2})$ to take
  normal form computations.

  \item[ii)] To establish the long time stability of the KAM tori so constructed,
  we shall take symplectic transformation of $f^{\textrm{high}}_{\infty}$ to obtain a suitable Birkhoff normal form. We will not put
      the frequency shift produced in the symplectic transformation into the homological equations, and we will finally get
      \begin{equation*}
        f^{\textrm{high; new}}_{\infty}=O(\|z\|^{M+1}).
      \end{equation*}
        It then follows that the KAM torus ($r=0, z=0$)
        is stable in a long time interval of length $\delta^{-M}$.

        In this process,  the tame property for space-multidimensional NLS equation can be preserved during the KAM iteration.
        Moreover,
        we will take the advantage of the momentum conservation.
The corresponding Hamiltonian $f$ consists of monomials
\begin{equation*}
e^{\mathrm{i}\langle k,\varphi\rangle}\prod_{a\in \mathcal{A}}r_{a}^{n_{a}}\prod_{a\in\mathcal{L}}u_{a}^{l_{a}}v_{a}^{m_{a}}
\end{equation*}
satisfying
\begin{equation*}
-\sum_{a\in \mathcal{A}}k_{a}a+\sum_{a\in\mathcal{L}}(l_{a}-m_{a})a=0.
\end{equation*}
We need to verify momentum conservation in the KAM iteration.
The persistence under Poisson bracket can be  checked directly.
Since the homological equations are of constant coefficients,
the persistence under solving homological equations can also be directly checked.
By the momentum conservation, we can deal with the frequency shift to establish
the long time stability.

\end{enumerate}

We end up this section with several remarks.

\begin{rem}
  In this paper, we benefit a lot from the  momentum conservation, which comes  from the
  $x$-independent nonlinearity $F$. See also \emph{\cite{GY,PP1,PP2}}. For the general case, there are extra difficulties in dealing with a block-diagonal shift of frequency.
\end{rem}

\begin{rem}
  As mentioned before, the existence of KAM tori $($quasi-periodic solution$)$ for space-multidimensional NLW equation
  can be obtained by the CWB method \emph{\cite{B2}}. However, on the one hand, a counterpart of KAM approach
  for NLW equation like Eliasson-Kuksin \emph{\cite{EK}} $($on NLS equation$)$ is still not available.
  See \emph{\cite{EGK}}. On the other hand, the CWB method
  does not provide a normal form of the Hamiltonian in the neighborhood of the KAM torus.
  As a result, the linear stability of KAM tori (quasiperiodic solutions) for NLW equation with $d>1$ is not clear, let alone the long time stability.
\end{rem}

The paper is organized as follows. In section \ref{sect 2}, we introduce
some  notations as the preliminary and present our main  results.
In section \ref{sect 3}, we formulate and solve the homological equation.
In section \ref{sect 4}, we prove the KAM theorem to show the existence of the KAM tori
for NLS equation.
In section \ref{sect 5}, we construct a partial normal form  to show the long time stability of the obtained KAM tori.


\section{Main results}\label{sect 2}

In this section, we present the main results of the infinitely dimensional
Hamiltonian system.  To begin with, we introduce some notations as the preliminary.

\subsection{Preliminary}
In this part, we collect some notations and definitions, which are frequently
used throughout the paper. In subsection \ref{sec 2.1}, we write NLS equation as
an infinitely dimensional Hamiltonian system. In subsection \ref{sec 2.2}, we introduce
the tame property of the Hamiltonian vector field. In subsection \ref{sec 2.3}, we
introduce the T\"{o}plitz-Lipschitz property. Finally, in subsection \ref{sec 2.4},
we introduce the normal form matrix.

\subsubsection{Hamiltonian formulation of NLS equation}\label{sec 2.1}

In order to prove Theorem \ref{t}, we write the nonlinear Schr\"odinger equation (\ref{1}) as an infinitely dimensional Hamiltonian system.
We keep the notations consistent with those in \cite{EK}.

Write
\begin{equation*}
u(x)=\sum_{a\in\mathbb{Z}^{d}}u_{a}e^{\mathrm{i}\langle a,x\rangle}, \ \overline{u(x)}=\sum_{a\in\mathbb{Z}^{d}}v_{a}e^{\mathrm{i}\langle -a,x\rangle},
\end{equation*}
and let
\begin{equation*}
\zeta_{a}=\left(
        \begin{array}{c}
          \xi_{a} \\
          \eta_{a} \\
        \end{array}
      \right)=\frac{1}{\sqrt{2}}\left(
                                  \begin{array}{c}
                                    u_{a}+v_{a} \\
                                    -\mathrm{i}(u_{a}-v_{a}) \\
                                  \end{array}
                                \right).
\end{equation*}
Then the nonlinear Schr\"odinger equation (\ref{1}) becomes a real Hamiltonian system
with the symplectic structure $d\xi\wedge d\eta$ and the Hamiltonian
\begin{equation*}
  \frac{1}{2}\sum_{a\in\mathbb{Z}^{d}}(|a|^{2}+\hat{V}(a))(\xi_{a}^{2}+\eta_{a}^{2})+\varepsilon \int_{\mathbb{T}^{d}}F(|u(x)|^{2})dx.
\end{equation*}
Let $\mathcal{A}$ be a finite subset of $\mathbb{Z}^{d}$ and $\mathcal{L}=\mathbb{Z}^{d}\setminus\mathcal{A}$.
Introduce action-angle variables $(\varphi_{a},r_{a})$, $a\in \mathcal{A}$,
\begin{equation*}
\xi_{a}=\sqrt{2(r_{a}+q_{a})}\cos\varphi_{a}, \eta_{a}=\sqrt{2(r_{a}+q_{a})}\sin\varphi_{a}, \ q_{a}>0.
\end{equation*}
Let
\begin{equation*}
\omega_{a}=|a|^{2}+\hat{V}(a), a\in \mathcal{A}  \quad\textrm{and}\quad \Omega_{a}=|a|^{2}+\hat{V}(a), a\in \mathcal{L}.
\end{equation*}
We have the Hamiltonian
\begin{equation*}
h+f=\sum_{a\in \mathcal{A}}\omega_{a}r_{a}+\frac{1}{2}\sum_{a\in\mathcal{L}}\Omega_{a}(\xi_{a}^{2}+\eta_{a}^{2})+\varepsilon \int_{\mathbb{T}^{d}}F(|u(x)|^{2})dx.
\end{equation*}
Assume $f$ is real analytic on
\begin{equation*}
D(\rho,\mu,\sigma)=\{(\varphi,r,\zeta)\in(\mathbb{C}/2\pi\mathbb{Z})^{\mathcal{A}}\times\mathbb{C}^{\mathcal{A}}\times l^{2}_{p}:|\Im\varphi|\leq\rho,|r|\leq\mu,\|\zeta\|_p\leq\sigma\},
\end{equation*}
where
\begin{equation*}
\|\zeta\|^{2}_{p}=\sum_{a\in\mathcal{L}}(|\xi^{2}_{a}|+|\eta^{2}_{a}|)\langle a\rangle^{2p}, \ \langle a\rangle=\max(|a|,1).
\end{equation*}

\subsubsection{The $p$-tame norm of the Hamiltonian vector field}\label{sec 2.2}
In this paper, $\|\cdot\|$ is an operator norm or $l^{2}$ norm. $|\cdot|$ will in general denote a sup norm.
For  $a \in \mathbb{Z}^{d}$, we use $|a|$ for the $l^{2}$ norm.
Let $\mathcal{A}$ be a finite subset of $\mathbb{Z}^{d}$ and $\mathcal{L}=\mathbb{Z}^{d}\setminus\mathcal{A}$.
Denote $\langle \zeta,\zeta'\rangle=\sum(\xi_{a}\xi'_{a}+\eta_{a}\eta'_{a})$ and $J=\left(
     \begin{array}{cc}
       0 & 1 \\
       -1 & 0 \\
     \end{array}
   \right)
$.

For $\gamma\geq 0$, we denote
\begin{equation*}
l^{2}_{p,\gamma}=\{\zeta=(\xi,\eta)\in\mathbb{C}^{\mathcal{L}}\times\mathbb{C}^{\mathcal{L}}:\|\zeta\|_{p,\gamma}<\infty\},
\end{equation*}
where
\begin{equation*}
\|\zeta\|^{2}_{p,\gamma}=\sum_{a\in\mathcal{L}}(|\xi^{2}_{a}|+|\eta^{2}_{a}|)e^{2\gamma|a|}\langle a\rangle^{2p}, \ \langle a\rangle=\max(|a|,1).
\end{equation*}
When $\gamma= 0$, we simply write $l^{2}_{p}$ and $\|\zeta\|_p$.
The phase space of the Hamiltonian dynamical system is defined by
\begin{equation*}
\mathcal{P}^{p}=(\mathbb{C}/2\pi\mathbb{Z})^{\mathcal{A}}\times\mathbb{C}^{\mathcal{A}}\times l^{2}_{p}.
\end{equation*}

Let $U\subset\mathbb{R}^{\mathbb{Z}^{d}}$ be a parameter set with positive measure (in the sense of Gauss or Kolmogorov).
We define $p$-tame norm as in \cite{CLY}.

\begin{defn}
Let
\begin{equation*}
D(\rho)=\{\varphi\in(\mathbb{C}/2\pi\mathbb{Z})^{\mathcal{A}}:|\Im\varphi|\leq\rho\},
\end{equation*}
and
$f:D(\rho)\times U\rightarrow \mathbb{C}$ be analytic in $\varphi\in D(\rho)$ and $C^1$ $($in the sense of Whitney$)$ in $w\in U$ with
\begin{equation*}
f(\varphi;w)=\sum_{k\in\mathbb{Z}^{\mathcal{A}}}\hat{f}(k;w)e^{\mathrm{i}\langle k,\varphi\rangle}.
\end{equation*}
Define the norm
\begin{equation*}
\|f\|_{D(\rho)\times U}=\sup_{w\in U,a \in \mathbb{Z}^{d}}\sum_{k\in\mathbb{Z}^{\mathcal{A}}}\left(|\hat{f}(k;w)|+|\partial_{w_{a}}\hat{f}(k;w)|\right)e^{|k|\rho}.
\end{equation*}
\end{defn}

\begin{defn}
Let
\begin{equation*}
D(\rho,\mu)=\{(\varphi,r)\in(\mathbb{C}/2\pi\mathbb{Z})^{\mathcal{A}}\times\mathbb{C}^{\mathcal{A}}:|\Im\varphi|\leq\rho,|r|\leq\mu\},
\end{equation*}
and
$f:D(\rho,\mu)\times U\rightarrow \mathbb{C}$ be analytic in $(\varphi,r)\in D(\rho,\mu)$ and $C^1$ in $w\in U$ with
\begin{equation*}
f(\varphi,r;w)=\sum_{\alpha\in\mathbb{N}^{\mathcal{A}}}f^{\alpha}(\varphi;w)r^{\alpha}.
\end{equation*}
Define the norm
\begin{equation*}
\|f\|_{D(\rho,\mu)\times U}=\sum_{\alpha\in\mathbb{N}^{\mathcal{A}}}\|f^{\alpha}(\varphi;w)\|_{D(\rho)\times U}\mu^{|\alpha|}.
\end{equation*}
\end{defn}

\begin{defn}
Let
\begin{equation*}
D(\rho,\mu,\sigma)=\{(\varphi,r,\zeta)\in(\mathbb{C}/2\pi\mathbb{Z})^{\mathcal{A}}\times\mathbb{C}^{\mathcal{A}}\times l^{2}_{p}:|\Im\varphi|\leq\rho,|r|\leq\mu,\|\zeta\|_p\leq\sigma\},
\end{equation*}
and
$f:D(\rho,\mu,\sigma)\times U\rightarrow \mathbb{C}$ be analytic in $(\varphi,r,\zeta)\in D(\rho,\mu,\sigma)$ and $C^1$ in $w\in U$ with
\begin{equation*}
f(\varphi,r,\zeta;w)=\sum_{\alpha\in\mathbb{N}^{\mathcal{A}},\beta\in\mathbb{N}^{\bar{\mathcal{L}}}}f^{\alpha\beta}(\varphi;w)r^{\alpha}\zeta^{\beta},
\end{equation*}
where $\bar{\mathcal{L}}=\mathcal{L}_{-1}\sqcup\mathcal{L},\mathcal{L}_{-1}=\mathcal{L}.$ For $a\in\mathcal{L}_{-1}, \zeta_{a}=\xi_{a}$, and for $a\in\mathcal{L}, \zeta_{a}=\eta_{a}$.
Define the modulus
\begin{equation*}
\lfloor f \rceil_{D(\rho,\mu)\times U}(\zeta)=\sum_{\beta\in\mathbb{N}^{\bar{\mathcal{L}}}}\|f^{\beta}(\varphi,r;w)\|_{D(\rho,\mu)\times U}\zeta^{\beta},
\end{equation*}
where
\begin{equation*}
f^{\beta}(\varphi,r;w)=\sum_{\alpha\in\mathbb{N}^{\mathcal{A}}}f^{\alpha\beta}(\varphi;w)r^{\alpha}.
\end{equation*}
\end{defn}

For a homogeneous polynomial $f(\zeta)$ of degree $h>0$, it is associated with a symmetric $h$-linear form $\tilde{f}(\zeta^{(1)},\ldots,\zeta^{(h)})$ such that
$\tilde{f}(\zeta,\ldots,\zeta)=f(\zeta)$. For a monomial
\begin{equation*}
f(\zeta)=f^{\beta}\zeta^{\beta}=f^{\beta}\zeta_{j_{1}}\cdots\zeta_{j_{h}},
\end{equation*}
define
\begin{equation*}
\tilde{f}(\zeta^{(1)},\ldots,\zeta^{(h)})=\widetilde{f^{\beta}\zeta^{\beta}}=\frac{1}{h!}\sum_{\tau_{h}}f^{\beta}\zeta_{j_{1}}^{(\tau_{h}(1))}\cdots\zeta_{j_{h}}^{(\tau_{h}(h))},
\end{equation*}
where $\tau_{h}$ is an $h$-permutation. For a homogeneous polynomial
\begin{equation*}
f(\zeta)=\sum_{|\beta|=h}f^{\beta}\zeta^{\beta},
\end{equation*}
define
\begin{equation*}
\tilde{f}(\zeta^{(1)},\ldots,\zeta^{(h)})=\sum_{|\beta|=h}\widetilde{f^{\beta}\zeta^{\beta}}.
\end{equation*}

Now we can define $p$-tame norm of a Hamiltonian vector field. We first consider a Hamiltonian
\begin{equation*}
f(\varphi,r,\zeta;w)=f_{h}=\sum_{\alpha\in\mathbb{N}^{\mathcal{A}},\beta\in\mathbb{N}^{\bar{\mathcal{L}}},|\beta|=h}f_{h}^{\alpha\beta}(\varphi;w)r^{\alpha}\zeta^{\beta}.
\end{equation*}
Let $f_{\zeta}=(f_{\eta},-f_{\xi})$, and the Hamiltonian vector field $X_f$ is $(f_{r},-f_{\varphi},f_{\zeta})$.
For $h\geq 1$, denote
\begin{equation*}
\|(\zeta^h)\|_{p,1}=\frac{1}{h}\sum_{j=1}^{h}\|\zeta^{(1)}\|_{1}\cdots\|\zeta^{(j-1)}\|_{1}\|\zeta^{(j)}\|_{p}\|\zeta^{(j+1)}\|_{1}\cdots\|\zeta^{(h)}\|_{1}.
\end{equation*}

\begin{defn}
Let
\begin{equation*}
|||f_{\zeta}|||^{T}_{p,D(\rho,\mu)\times U}=\left\{
\begin{array}{cl}
  \sup_{0\neq\zeta^{(j)}\in l^2_p, 1\leq j\leq h-1}\frac{\|\widetilde{\lfloor f_{\zeta} \rceil}_{D(\rho,\mu)\times U}
(\zeta^{(1)},\ldots,\zeta^{(h-1)})\|_{p}}{\|(\zeta^{h-1})\|_{p,1}}, & h\geq 2\\
\|\widetilde{\lfloor f_{\zeta} \rceil}_{D(\rho,\mu)\times U}\|_{p}, & h=0,1.
\end{array}
\right.
\end{equation*}
Define the $p$-tame norm of $f_{\zeta}$ by
\begin{equation*}
|||f_{\zeta}|||^{T}_{p,D(\rho,\mu,\sigma)\times U}=\max(|||f_{\zeta}|||^{T}_{p,D(\rho,\mu)\times U},|||f_{\zeta}|||^{T}_{1,D(\rho,\mu)\times U})\sigma^{h-1}.
\end{equation*}
\end{defn}

\begin{defn}
Let
\begin{equation*}
|||f_{r}|||_{D(\rho,\mu)\times U}=
\left\{
\begin{array}{cl}
\sup_{0\neq\zeta^{(j)}\in l^2_1, 1\leq j\leq h}\frac{|\widetilde{\lfloor f_{r} \rceil}_{D(\rho,\mu)\times U}
(\zeta^{(1)},\ldots,\zeta^{(h)})|}{\|(\zeta^{h})\|_{1,1}}, &  h\geq 1,\\
|\widetilde{\lfloor f_{r} \rceil}_{D(\rho,\mu)\times U}|, & h=0.
\end{array}
\right.
\end{equation*}
Define the norm of $f_{r}$ by
\begin{equation*}
|||f_{r}|||_{D(\rho,\mu,\sigma)\times U}=|||f_{r}|||_{D(\rho,\mu)\times U}\sigma^{h}.
\end{equation*}
The norm of $f_{\varphi}$ is defined as that of $f_{r}$.
\end{defn}

\begin{defn}
Define the $p$-tame norm of the Hamiltonian vector field $X_f$ by
\begin{equation*}
|||X_f|||^{T}_{p,D(\rho,\mu,\sigma)\times U}=|||f_{r}|||_{D(\rho,\mu,\sigma)\times U}+\frac{1}{\mu}|||f_{\varphi}|||_{D(\rho,\mu,\sigma)\times U}+\frac{1}{\sigma}|||f_{\zeta}|||^{T}_{p,D(\rho,\mu,\sigma)\times U}.
\end{equation*}
\end{defn}

\begin{defn}
For a Hamiltonian
\begin{equation*}
f(\varphi,r,\zeta;w)=\sum_{h\geq 0}f_{h}, \quad f_{h}=\sum_{\alpha\in\mathbb{N}^{\mathcal{A}},\beta\in\mathbb{N}^{\bar{\mathcal{L}}},|\beta|=h}f_{h}^{\alpha\beta}(\varphi;w)r^{\alpha}\zeta^{\beta},
\end{equation*}
define the $p$-tame norm of the Hamiltonian vector field $X_f$ by
\begin{equation*}
|||X_f|||^{T}_{p,D(\rho,\mu,\sigma)\times U}=\sum_{h\geq 0}|||X_{f_{h}}|||^{T}_{p,D(\rho,\mu,\sigma)\times U}.
\end{equation*}
\end{defn}

\begin{rem}
The $p$-tame norm can also be defined in complex coordinates
\begin{equation*}
z=\left(
    \begin{array}{c}
      u \\
      v \\
    \end{array}
  \right)
  =C^{-1}\left(
           \begin{array}{c}
             \xi \\
             \eta \\
           \end{array}
         \right),
         C=\frac{1}{\sqrt{2}}\left(
         \begin{array}{cc}
           1 & 1 \\
           -\mathrm{i} & \mathrm{i} \\
         \end{array}
       \right).
\end{equation*}
\end{rem}

Following the proof of Theorem 3.1 in \cite{CLY}, we have the following proposition.

\begin{prop}\label{p}
If $0<\tau<\rho, 0<\tau'<\frac{\sigma}{2}$, then
\begin{equation*}
|||X_{\{f,g\}}|||^{T}_{p,D(\rho-\tau,(\sigma-\tau')^{2},\sigma-\tau')\times U}\leq
C\max\left(\frac{1}{\tau},\frac{\sigma}{\tau'}\right)|||X_{f}|||^{T}_{p,D(\rho,\sigma^{2},\sigma)\times U}|||X_{g}|||^{T}_{p,D(\rho,\sigma^{2},\sigma)\times U},
\end{equation*}
where $C>0$ is a constant depending on $\# \mathcal{A}$.
\end{prop}

We define the weighted norm of the Hamiltonian vector field $X_f$ by
\begin{equation*}
|||X_f|||_{\mathcal{P}^{p},D(\rho,\mu,\sigma)\times U}=\sup_{(\varphi,r,\zeta;w)\in D(\rho,\mu,\sigma)\times U}\|X_{f}\|_{\mathcal{P}^{p},D(\rho,\mu,\sigma)},
\end{equation*}
where
\begin{equation*}
\|X_{f}\|_{\mathcal{P}^{p},D(\rho,\mu,\sigma)}=|f_{r}|+\frac{1}{\mu}|f_{\varphi}|+\frac{1}{\sigma}\|f_{\zeta}\|_{p}.
\end{equation*}
Following the proof of Theorem 3.5 in \cite{CLY}, we have
\begin{equation*}
|||X_f|||_{\mathcal{P}^{p},D(\rho,\mu,\sigma)\times U}\leq|||X_f|||^{T}_{p,D(\rho,\mu,\sigma)\times U}.
\end{equation*}

\subsubsection{T\"{o}plitz-Lipschitz property}\label{sec 2.3}
Recall the definition of T\"{o}plitz-Lipschitz matrices in \cite{EK}.
Let $gl(2,\mathbb{C})$ be the space of all complex $2\times2$-matrices.
For $A=\left(
         \begin{array}{cc}
           a & b \\
           c & d \\
         \end{array}
       \right)
\in gl(2,\mathbb{C})$, denote $\pi A=\frac{1}{2}\left(
         \begin{array}{cc}
           a+d & b-c \\
           c-b & a+d \\
         \end{array}
       \right)
$ and $[A]=\left(
         \begin{array}{cc}
           |a| & |b| \\
           |c| & |d| \\
         \end{array}
       \right)
$. Now consider an infinite-dimensional $gl(2,\mathbb{C})$-valued matrix
\begin{equation*}
A:\mathcal{L}\times\mathcal{L}\rightarrow gl(2,\mathbb{C}), \ (a,b)\mapsto A_{a}^{b}.
\end{equation*}
For any $\mathcal{D}\subset\mathcal{L}\times\mathcal{L}$, define
\begin{equation*}
|A|_\mathcal{D}=\sup_{(a,b)\in\mathcal{D}}\|A_{a}^{b}\|,
\end{equation*}
where $\|\cdot\|$ is the operator norm. Define $(\pi A)_{a}^{b}=\pi A_{a}^{b}$ and $(\mathcal{E}_{\gamma}^{\pm}A)_{a}^{b}=[A_{a}^{b}]e^{\gamma|a\mp b|}$.
Define the norm
\begin{equation*}
|A|_\gamma=\max(|\mathcal{E}_{\gamma}^{+}\pi A|_{\mathcal{L}\times\mathcal{L}},|\mathcal{E}_{\gamma}^{-}(1-\pi) A|_{\mathcal{L}\times\mathcal{L}}).
\end{equation*}
Define
\begin{equation*}
\mathcal{T}_{\Delta}^{\pm}A=A\mid_{\{(a,b)\in\mathcal{L}\times\mathcal{L}:|a\mp b|\leq\Delta\}}, \ \mathcal{T}_{\Delta}A=\mathcal{T}_{\Delta}^{+}\pi A+\mathcal{T}_{\Delta}^{-}(1-\pi) A.
\end{equation*}

A matrix $A:\mathcal{L}\times\mathcal{L}\rightarrow gl(2,\mathbb{C})$ is T\"{o}plitz at $\infty$, if for all $a,b,c$, the two limits
\begin{equation*}
\lim_{t\rightarrow +\infty}A_{a+tc}^{b\pm tc}\exists =A_{a}^{b}(\pm,c).
\end{equation*}

For $c\neq 0$, define $(\mathcal{M}_{c}A)_{a}^{b}=\left(\max(\frac{|a|}{|c|},\frac{|b|}{|c|})+1\right)[A_{a}^{b}]$.
For $\Lambda\geq 0$, define the Lipschitz domain $D_{\Lambda}^{+}(c)\subset\mathcal{L}\times\mathcal{L}$ be the set of all $(a,b)$ such that
there exist $a',b'\in \mathbb{Z}^{d}, t\geq 0$ such that
\begin{equation*}
|a=a'+tc|\geq\Lambda(|a'|+|c|)|c|, \ |b=b'+tc|\geq\Lambda(|b'|+|c|)|c|, \ \frac{|a|}{|c|},\frac{|b|}{|c|}\geq 2\Lambda^{2}.
\end{equation*}
Define  $(a,b)\in D_{\Lambda}^{-}(c)$ if and only if $(a,-b)\in D_{\Lambda}^{+}(c)$.
Define the Lipschitz constants
\begin{equation*}
\textrm{Lip}_{\Lambda,\gamma}^{\pm}A=\sup_{c}
|\mathcal{E}_{\gamma}^{\pm}\mathcal{M}_{c}(A-A(\pm,c))|_{D_{\Lambda}^{\pm}(c)},
\end{equation*}
and the Lipschitz norm
\begin{equation*}
\langle A\rangle_{\Lambda,\gamma}=\max(\textrm{Lip}_{\Lambda,\gamma}^{+}\pi A,\textrm{Lip}_{\Lambda,\gamma}^{-}(1-\pi) A)+|A|_\gamma.
\end{equation*}
For $d=2$, the matrix $A$ is T\"{o}plitz-Lipschitz if it is T\"{o}plitz at $\infty$ and $\langle A\rangle_{\Lambda,\gamma}<\infty$ for some $\Lambda,\gamma$.
For $d>2$, we can define T\"{o}plitz-Lipschitz matrices inductively (see Section 2.4 in \cite{EK}).

\begin{defn}
Let
\begin{equation*}
D^{\gamma}(\sigma)=\{\zeta\in l^{2}_{p,\gamma}:\|\zeta\|_{p,\gamma}\leq\sigma\},
\end{equation*}
and
$f:D^{0}(\sigma)\rightarrow \mathbb{C}$ be real analytic. We say that
$f$ is T\"{o}plitz at $\infty$ if $\partial^{2}_{\zeta}f(\zeta)$ is T\"{o}plitz at $\infty$ for all $\zeta\in D^{0}(\sigma)$.
Define the norm $[f]_{\Lambda,\gamma,\sigma}$ to be the smallest $C$ such that
\begin{equation*}
\begin{aligned}
& |f(\zeta)|\leq C,~ \forall~ \zeta\in D^{0}(\sigma), \\ & \|\partial_{\zeta}f(\zeta)\|_{p,\gamma'}\leq\frac{C}{\sigma},
~\langle\partial^{2}_{\zeta}f(\zeta)\rangle_{\Lambda,\gamma'}
\leq\frac{C}{\sigma^{2}},~ \forall~\zeta\in D^{\gamma'}(\sigma),~\forall~\gamma'\leq\gamma.
\end{aligned}
\end{equation*}
\end{defn}

\begin{defn}
Let $A(w):\mathcal{L}\times\mathcal{L}\rightarrow gl(2,\mathbb{C})$ be $C^1$ in $w\in U$. Define
\begin{equation*}
|A|_{\gamma;U}=\sup_{w\in U}(|A(w)|_{\gamma},|\partial_{w}A(w)|_{\gamma}).
\end{equation*}
If $A(w),\partial_{w}A(w)$ are T\"{o}plitz at $\infty$ for all $w\in U$, define
\begin{equation*}
\langle A\rangle_{\Lambda,\gamma;U}=\sup_{w\in U}(\langle A(w)\rangle_{\Lambda,\gamma},\langle \partial_{w}A(w)\rangle_{\Lambda,\gamma}).
\end{equation*}
When $\gamma=0$, we simply write $|A|_{U},\langle A\rangle_{\Lambda;U}$.
\end{defn}

\subsubsection{The Normal form matrix}\label{sec 2.4}
For $\Delta\geq 0$, define an equivalence relation on  $\mathcal{L}$ generated by the pre-equivalence relation
\begin{equation*}
a\sim b \Leftrightarrow |a|=|b|,|a-b|\leq\Delta.
\end{equation*}
Let $[a]_\Delta$ be the equivalence class (block) of $a$ and $\mathcal{E}_\Delta$ be the set of equivalence classes.
Let $d_\Delta$ be the supremum of all block diameters, then by Proposition 4.1 in \cite{EK}, $d_\Delta\preceq \Delta^{\frac{(d+1)!}{2}}$.

A matrix $A:\mathcal{L}\times\mathcal{L}\rightarrow gl(2,\mathbb{C})$ is on normal form, denoted $\mathcal{NF}_{\Delta}$, if
$A$ is real valued, symmetric, $\pi A=A$ and block-diagonal over $\mathcal{E}_\Delta$, i.e., $ A_{a}^{b}=0, \forall~[a]_\Delta\neq [b]_\Delta$.
A matrix $Q:\mathcal{L}\times\mathcal{L}\rightarrow \mathbb{C}$ is on normal form, denoted $\mathcal{NF}_{\Delta}$, if
$Q$ is Hermitian and block-diagonal over $\mathcal{E}_\Delta$.

For a normal form matrix $A$,
\begin{equation*}
\frac{1}{2}\langle\zeta,A\zeta\rangle=\frac{1}{2}\langle\xi,A_{1}\xi\rangle+\langle\xi,A_{2}\eta\rangle+\frac{1}{2}\langle\eta,A_{1}\eta\rangle,
\end{equation*}
where $A_{1}+\mathrm{i}A_{2}$ is a Hermitian matrix.
Let
\begin{equation*}
z=\left(
    \begin{array}{c}
      u \\
      v \\
    \end{array}
  \right)
  =C^{-1}\left(
           \begin{array}{c}
             \xi \\
             \eta \\
           \end{array}
         \right),
         C=\frac{1}{\sqrt{2}}\left(
         \begin{array}{cc}
           1 & 1 \\
           -\mathrm{i} & \mathrm{i} \\
         \end{array}
       \right),
\end{equation*}
and define $C^{T}AC:\mathcal{L}\times\mathcal{L}\rightarrow gl(2,\mathbb{C})$ by $(C^{T}AC)_{a}^{b}=C^{T}A_{a}^{b}C$.
If $A$ is on normal form, then
\begin{equation*}
\frac{1}{2}\langle z,C^{T}ACz\rangle=\langle u,Qv\rangle,
\end{equation*}
where $Q$ is the normal form matrix associated to $A$.

\subsection{Main results}

Let
\begin{equation*}
h(r,\zeta;w)=\langle\omega(w),r\rangle+\frac{1}{2}\langle\zeta,(\Omega(w)+H(w))\zeta\rangle,
\end{equation*}
where $\Omega(w)$ is a real diagonal matrix with diagonal elements $\Omega_{a}(w)I$,
$H(w),\partial_{w}H(w)$ are T\"{o}plitz at $\infty$ and $\mathcal{NF}_{\Delta}$ for all $w\in U$.

Assume
\begin{equation}\label{as1}
\partial_{w_{a}}\omega_{b}(w)=\delta_{ab}, \ a\in \mathbb{Z}^{d}, b\in \mathcal{A}, w\in U,
\end{equation}
\begin{equation}\label{as2}
\partial_{w_{a}}\Omega_{b}(w)=\delta_{ab}, \ a\in \mathbb{Z}^{d}, b\in \mathcal{L}, w\in U.
\end{equation}

Assume there exist constants $c_1,c_2,c_3,c_4,c_5>0$ such that
\begin{equation}\label{as3}
|\Omega_{a}(w)-|a|^{2}|\leq c_1e^{-c_2|a|}, \ a\in \mathcal{L}, w\in U,
\end{equation}
\begin{equation}\label{as4}
|\Omega_{a}(w)|\geq c_3, \ a\in \mathcal{L}, w\in U,
\end{equation}
\begin{equation}\label{as5}
|\Omega_{a}(w)+\Omega_{b}(w)|\geq c_3, \ a,b\in \mathcal{L}, w\in U,
\end{equation}
\begin{equation}\label{as6}
|\Omega_{a}(w)-\Omega_{b}(w)|\geq c_3, \   |a|\neq |b|,      \ a,b\in \mathcal{L}, w\in U,
\end{equation}
\begin{equation}\label{as7}
\|H(w)\|\leq \frac{c_3}{4}, \  w\in U,
\end{equation}
\begin{equation}\label{as8}
\|\partial_{w}H(w)\|\leq c_4, w\in U,
\end{equation}
\begin{equation}\label{as9}
\langle H \rangle_{\Lambda;U}\leq c_5.
\end{equation}

Let
\begin{equation*}
D^{\gamma}(\rho,\mu,\sigma)=\{(\varphi,r,\zeta)\in(\mathbb{C}/2\pi\mathbb{Z})^{\mathcal{A}}\times\mathbb{C}^{\mathcal{A}}\times l^{2}_{p,\gamma}:|\Im\varphi|\leq\rho,|r|\leq\mu,\|\zeta\|_{p,\gamma}\leq\sigma\},
\end{equation*}
and
$f:D^{\gamma}(\rho,\mu,\sigma)\times U\rightarrow \mathbb{C}$ be real analytic in $(\varphi,r,\zeta)\in D^{\gamma}(\rho,\mu,\sigma)$ and $C^1$ in $w\in U$.
Define
\begin{equation*}
[f]_{\Lambda,\gamma,\sigma;U,\rho,\mu} =\sup_{(\varphi,r)\in D(\rho,\mu)}[f(\varphi,r,\cdot;\cdot)]_{\Lambda,\gamma,\sigma;U}.
\end{equation*}

\begin{thm}\label{t1}

Consider the Hamiltonian $h+f$, where
\begin{equation*}
h(r,\zeta;w)=\langle\omega(w),r\rangle+\frac{1}{2}\langle\zeta,(\Omega(w)+H(w))\zeta\rangle
\end{equation*}
satisfy \eqref{as1}-\eqref{as9}, $H(w),\partial_{w}H(w)$ are T\"{o}plitz at $\infty$ and $\mathcal{NF}_{\Delta}$ for all $w\in U$,
\begin{equation}\label{k1}
|||X_{f}|||^{T}_{p,D(\rho,\mu,\sigma)\times U}\leq\varepsilon,
\end{equation}
\begin{equation}\label{k2}
[f]_{\Lambda,\gamma,\sigma;U,\rho,\mu}\leq\varepsilon.
\end{equation}
Assume $\gamma,\sigma,\rho,\mu<1$, $\Lambda,\Delta\geq 3$, $\rho=\sigma$, $\mu=\sigma^{2}$, $d_{\Delta}\gamma\leq1$.
Then there is a subset $U_{\infty}\subset U$ such that if
\begin{equation*}
\varepsilon\preceq\min\left(\gamma,\rho,\frac{1}{\Delta},\frac{1}{\Lambda}\right)^{\exp},
\end{equation*}
then for all $w\in U_{\infty}$, there is a real analytic symplectic map
\begin{equation*}
\Phi:D(\frac{\rho}{2},\frac{\mu}{4},\frac{\sigma}{2})\rightarrow D(\rho,\mu,\sigma)
\end{equation*}
such that
\begin{equation*}
(h+f)\circ\Phi=h_{\infty}+f_{\infty},
\end{equation*}
where
\begin{equation*}
h_{\infty}=\langle \omega_{\infty}(w),r\rangle+\frac{1}{2}\langle\zeta,(\Omega(w)+H_{\infty}(w))\zeta\rangle,
\end{equation*}
\begin{equation*}
f_{\infty}=O(|r|^{2}+|r|\|\zeta\|_{p}+\|\zeta\|^{3}_{p})
\end{equation*}
with the estimates
\begin{equation}\label{k3}
|||X_{f_{\infty}}|||^{T}_{p,D(\frac{\rho}{2},\frac{\mu}{2},\frac{\sigma}{2})\times U_{\infty}}
\leq c\varepsilon^{\frac{2}{3}},
\end{equation}
\begin{equation}\label{k4}
|\omega_{\infty}(w)-\omega(w)|+|\partial_{w}(\omega_{\infty}(w)-\omega(w))|\leq c\varepsilon^{\frac{2}{3}},
\end{equation}
\begin{equation}\label{k5}
\|H_{\infty}(w)-H(w)\|+\|\partial_{w}( H_{\infty}(w)-H(w))\|\leq c\varepsilon^{\frac{2}{3}},
\end{equation}
\begin{equation}\label{k6}
\mathrm{meas}(U\setminus U_{\infty})\preceq \varepsilon^{\exp'}.
\end{equation}
The exponents $\exp$, $\exp'$ depend on $d, \# \mathcal{A}, p$, and the constant $c$ depends on
$d,\# \mathcal{A}, p, c_1,\cdots,c_5$.
\end{thm}

The proof of Theorem \ref{t1} is delayed to section \ref{sect 4}.

\begin{thm}\label{t2}

Given any $1\leq M \leq (4c\varepsilon^{\frac{2}{3}})^{-1}$ and $p\geq80(4d)^{4d}(M+7)^{4}+1$, there exists a set $\tilde{U}\subset U_{\infty}$ such that
 for any
$\delta>0$ and $w\in \tilde{U}$, the KAM tori obtained in Theorem \ref{t1} is stable in long time $\delta^{-M}$.
Moreover, $\mathrm{meas}(U_{\infty}\setminus\tilde{U})\preceq \delta^{\exp}$,
where the positive exponent $\exp$ depends on $d, \# \mathcal{A}, p, M$.
\end{thm}

The proof of Theorem \ref{t2} is given at the end of section \ref{sect 5}.
Based on Theorem \ref{t1} and Theorem \ref{t2}, we are able to prove
Theorem \ref{t} on the long time stability of the KAM tori for NLS equation.

\bigskip

\noindent\textbf{Proof of Theorem \ref{t}.}
Recall the Hamiltonian formulation of NLS equation \eqref{1} in subsection \ref{sec 2.1}.
Let
$\omega_{a}=|a|^{2}+\hat{V}(a), a\in \mathcal{A},$
$\Omega_{a}=|a|^{2}+\hat{V}(a), a\in \mathcal{L}$,
and take $w_{a}=\hat{V}(a)$. Then
we have
\begin{equation*}
h=\sum_{a\in \mathcal{A}}\omega_{a}r_{a}+\frac{1}{2}\sum_{a\in\mathcal{L}}\Omega_{a}(\xi_{a}^{2}+\eta_{a}^{2}),\quad
f=\varepsilon \int_{\mathbb{T}^{d}}F(|u(x)|^{2})dx.
\end{equation*}
The T\"{o}plitz-Lipschitz property of $f$ follows from Theorem 7.2 in \cite{EK} and the tame property follows from Section 3.5 in \cite{BG}.
By Theorem \ref{t1}, if $\varepsilon>0$ is sufficiently small, then for typical $V$ (in the sense of measure), the $d$-dimensional nonlinear Schr\"odinger equation (\ref{1})
has a quasi-periodic solution.
By Theorem \ref{t2}, assume $u_{0}(t,x)$ with initial value $u_{0}(0,x)$ is a quasi-periodic solution for the equation (\ref{1}),
then for any solution $u(t,x)$ with initial value $u(0,x)$ satisfying
\begin{equation*}
\|u(0,\cdot)-u_{0}(0,\cdot)\|_{H^{p}(\mathbb{T}^{d})}<\delta, \ \forall~ 0<\delta\ll1,
\end{equation*}
we have
\begin{equation*}
\|u(t,\cdot)-u_{0}(t,\cdot)\|_{H^{p}(\mathbb{T}^{d})}<C\delta, \ \forall~ 0<|t|<\delta^{-M}.
\end{equation*}
In other words, the obtained KAM tori for the nonlinear Schr\"odinger equation (\ref{1}) are of long time stability.
\qed

\section{The homological equations}\label{sect 3}

In this section, we formulate and solve the homological equation in the KAM iteration.
To obtain an open and uniform domain for the transformed Hamiltonian, we apply Kolmogorov's iterative
scheme. As a result, the homological equation is complicated than that in \cite{EK}, but it can be solved
by the method developed in \cite{EK}.

Write
\begin{equation*}
f(\varphi,r,\zeta;w)=f^{low}+f^{high},
\end{equation*}
where
\begin{equation*}
f^{low}=f^{\varphi}+f^{0}+f^{1}+f^{2}
=F^{\varphi}(\varphi;w)+\langle F_{0}(\varphi;w),r\rangle+\langle F_{1}(\varphi;w),\zeta\rangle+\frac{1}{2}\langle F_{2}(\varphi;w)\zeta,\zeta\rangle.
\end{equation*}
Define
\begin{equation*}
\mathcal{T}_{\Delta}f^{low}=
\sum_{|k|\leq\Delta}\left(\hat{F^{\varphi}}(k;w)+\langle\hat{F_{0}}(k;w),r\rangle+\langle\hat{F_{1}}(k;w),\zeta\rangle
+\frac{1}{2}\langle\mathcal{T}_{\Delta}\hat{F_{2}}(k;w)\zeta,\zeta\rangle\right)e^{\mathrm{i}\langle k,\varphi\rangle}.
\end{equation*}

Let $\Delta'>1$ and $0<\kappa<1$. Assume there exists $U' \subset U$ such that
for all $w\in U'$, $0<|k|\leq\Delta'$, the following properties hold:
\begin{itemize}
\item Diophantine condition:
\begin{equation}\label{sd1}
|\langle k,\omega(w)\rangle|\geq \kappa;
\end{equation}

\item The first Melnikov condition:
\begin{equation}\label{sd2}
|\langle k,\omega(w)\rangle+\alpha(w)|\geq \kappa, \quad \forall~\alpha(w)\in\mathrm{spec}(((\Omega+H)(w))_{[a]_\Delta}), \ \forall~ [a]_\Delta;
\end{equation}

\item The second Melnikov condition with the same sign:
\begin{equation}\label{sd3}
|\langle k,\omega(w)\rangle+\alpha(w)+\beta(w)|\geq \kappa, \quad \forall
\left\{
\begin{aligned}
\alpha(w)\in\mathrm{spec}(((\Omega+H)(w))_{[a]_\Delta}),\\
\beta(w)\in\mathrm{spec}(((\Omega+H)(w))_{[b]_\Delta}),
\end{aligned}
\right.
~\forall~ [a]_\Delta, [b]_\Delta;
\end{equation}

 \item The  second Melnikov condition with opposite signs:
\begin{equation}\label{sd4}
|\langle k,\omega(w)\rangle+\alpha(w)-\beta(w)|\geq \kappa, \quad
 \forall
 \left\{\begin{aligned}
 \alpha(w)\in\mathrm{spec}(((\Omega+H)(w))_{[a]_\Delta}),\\
 \beta(w)\in\mathrm{spec}(((\Omega+H)(w))_{[b]_\Delta}),
 \end{aligned}\right.
\end{equation}
and for any $\mathrm{dist}([a]_\Delta, [b]_\Delta)\leq \Delta'+2d_\Delta.$
\end{itemize}

We have the following result on the solution of the homological equation.

\begin{prop}\label{p1}

Consider the Hamiltonian $h+f$, where
\begin{equation*}
h(r,\zeta;w)=\langle\omega(w),r\rangle+\frac{1}{2}\langle\zeta,(\Omega(w)+H(w))\zeta\rangle
\end{equation*}
satisfy \eqref{as1}-\eqref{as9}, $H(w),\partial_{w}H(w)$ are T\"{o}plitz at $\infty$ and $\mathcal{NF}_{\Delta}$ for all $w\in U$,
\begin{equation*}
f(\varphi,r,\zeta;w)=f^{low}+f^{high}
\end{equation*}
satisfy
\begin{equation}\label{ho1}
|||X_{f^{low}}|||^{T}_{p,D(\rho,\mu,\sigma)\times U}\leq\varepsilon, \  |||X_{f^{high}}|||^{T}_{p,D(\rho,\mu,\sigma)\times U}\leq1,
\end{equation}
\begin{equation}\label{ho2}
[f^{low}]_{\Lambda,\gamma,\sigma;U,\rho,\mu}\leq\varepsilon, \ [f^{high}]_{\Lambda,\gamma,\sigma;U,\rho,\mu}\leq 1.
\end{equation}
Assume $\gamma,\sigma,\rho,\mu<1$, $\Lambda,\Delta\geq 3$, $\rho=\sigma$, $\mu=\sigma^{2}$, $d_{\Delta}\gamma\leq1$.
Let $U' \subset U$ satisfy \eqref{sd1}-\eqref{sd4}.

Then for all $w\in U'$, the homological equation
\begin{equation}\label{ho3}
\{h,s\}=-\mathcal{T}_{\Delta'}f^{low}-\mathcal{T}_{\Delta'}\{f^{high},s\}^{low}+h_1
\end{equation}
has solutions
\begin{equation}\label{ho4}
s(\varphi,r,\zeta;w)=s^{low}=s^{\varphi}+s^{0}+s^{1}+s^{2},
\end{equation}
\begin{equation}\label{ho5}
h_{1}(r,\zeta;w)=a_{1}(w)+\langle \chi_{1}(w),r\rangle+\frac{1}{2}\langle\zeta,H_{1}(w)\zeta\rangle
\end{equation}
with the estimates
\begin{equation}\label{ho6}
|||X_{s^{\varphi}}|||^{T}_{p,D(\rho-\tau,\sigma^{2},\sigma)\times U'}\preceq \frac{\varepsilon}{\tau\kappa^{2}},
\end{equation}
\begin{equation}\label{ho7}
|||X_{s^{1}}|||^{T}_{p,D(\rho-3\tau,(\sigma-3\tau)^{2},\sigma-3\tau)\times U'}
\preceq\frac{d_{\Delta}^{d}\varepsilon}{\tau^{3}\kappa^{4}},
\end{equation}
\begin{equation}\label{ho8}
|||X_{s^{0}}|||^{T}_{p,D(\rho-5\tau,(\sigma-5\tau)^{2},\sigma-5\tau)\times U'}\preceq \frac{d_{\Delta}^{d}\varepsilon}{\tau^{5}\kappa^{6}},
\end{equation}
\begin{equation}\label{ho9}
|||X_{s^{2}}|||^{T}_{p,D(\rho-5\tau,(\sigma-5\tau)^{2},\sigma-5\tau)\times U'}\preceq
\frac{d_{\Delta}^{3d}\varepsilon}{\tau^{5}\kappa^{6}},
\end{equation}
\begin{equation}\label{ho10}
|||X_{s}|||^{T}_{p,D(\rho-5\tau,(\sigma-5\tau)^{2},\sigma-5\tau)\times U'}\preceq\frac{d_{\Delta}^{3d}\varepsilon}{\tau^{5}\kappa^{6}},
\end{equation}
\begin{equation}\label{ho11}
|||X_{h_{1}}|||^{T}_{p,D(\rho-5\tau,(\sigma-5\tau)^{2},\sigma-5\tau)\times U'}
\preceq\frac{d_{\Delta}^{d}\varepsilon}{\tau^{4}\kappa^{4}},
\end{equation}
where $0<\tau<\frac{\rho}{100}$, $a\preceq b$ means there exists a constant $c>0$ depending on
$d,\# \mathcal{A}, p, c_1,\cdots,c_5$ such that $a\leq cb$.

The new Hamiltonian
\begin{equation}\label{ho12}
(h+f)\circ X^{t}_{s}\mid_{t=1}=h+h_1+f_1
\end{equation}
with
\begin{equation}\label{ho13}
\begin{aligned}
f_1=&(1-\mathcal{T}_{\Delta'})f^{low}+f^{high}+(1-\mathcal{T}_{\Delta'})
\{f^{high},s\}^{low}+\{f^{high},s\}^{high}\\
&+\int_{0}^{1}(1-t)\{\{h,s\},s\}\circ X^{t}_{s}dt+\int_{0}^{1}\{f^{low},s\}\circ X^{t}_{s}dt\\
&+\int_{0}^{1}(1-t)\{\{f^{high},s\},s\}\circ X^{t}_{s}dt
\end{aligned}
\end{equation}
satisfies
\begin{equation}\label{ho14}
|||X_{f_{1}^{low}}|||^{T}_{p,D(\rho-8\tau,(\sigma-8\tau)^{2},\sigma-8\tau)\times U'}
\preceq\frac{d_{\Delta}^{d}e^{-\frac{1}{2}\tau\Delta'}}{\kappa^{4}\tau^{\# \mathcal{A}+4}}\varepsilon
+\frac{(\Delta\Delta')^{\exp}}{\sigma^{6}\kappa^{4}}\frac{e^{-\frac{1}{2}\gamma\Delta'}}{\gamma^{d+p}\tau^{\# \mathcal{A}+3}}\varepsilon+\frac{d_{\Delta}^{6d}\varepsilon^{2}}{\tau^{12}\kappa^{12}},
\end{equation}
\begin{equation}\label{ho15}
|||X_{f_{1}^{high}}|||^{T}_{p,D(\rho-8\tau,(\sigma-8\tau)^{2},\sigma-8\tau)\times U'}\preceq 1+\frac{d_{\Delta}^{3d}\varepsilon}{\tau^{6}\kappa^{6}}+\frac{d_{\Delta}^{6d}\varepsilon^{2}}{\tau^{12}\kappa^{12}},
\end{equation}
where the exponent $\exp$ depends on $d,\# \mathcal{A}, p$.

Moreover, the following estimates hold.

\begin{itemize}
\item[i)] The solution $s$ and the remainder $h_{1}$ satisfy
\begin{equation}\label{ho16}
[s]_{\Lambda'+d_\Delta+2,\gamma,\sigma';U',\rho',\mu'}
\preceq\frac{1}{\kappa^{7}}(\Delta\Delta')^{\exp}\frac{1}{\rho-\rho'}\left(\frac{1}{\sigma-\sigma'}\frac{1}{\sigma}+\frac{1}{\rho-\rho'}\frac{1}{\mu-\mu'}\right)\frac{1}{\mu}\varepsilon,
\end{equation}
\begin{equation}\label{ho17}
[h_{1}]_{\Lambda'+d_\Delta+2,\gamma,\sigma';U',\rho',\mu'}
\preceq\frac{1}{\kappa^{4}}(\Delta\Delta')^{\exp}\frac{1}{\rho-\rho'}\left(\frac{1}{\sigma-\sigma'}\frac{1}{\sigma}+\frac{1}{\rho-\rho'}\frac{1}{\mu-\mu'}\right)\frac{1}{\mu}\varepsilon,
\end{equation}
where $\rho'<\rho$, $\mu'<\mu$, $\sigma'<\sigma$, $\Lambda'\geq \mathrm{cte}.\max(\Lambda,d^{2}_{\Delta},d^{2}_{\Delta'})$, and
the constant $\mathrm{cte}$. is the one  in \emph{\cite[Proposition 6.7]{EK}}.

\item[ii)] There is the measure estimate
\begin{equation}\label{ho18}
\mathrm{meas}(U\backslash U')\preceq\max(\Lambda,\Delta,\Delta')^{\exp}\kappa^{(\frac{1}{d+1})^{d}}.
\end{equation}
\end{itemize}

\end{prop}

The measure estimate (\ref{ho18}) follows directly from  \cite[Proposition 6.6 and 6.7 ]{EK},
which is not repeated here.
The rest of this section is devoted to the  proof of  Proposition \ref{p1}.
For the sake of notations, we shall not indicate
the dependence on the parameter $w$ of functions when it is known from the text.

\smallskip

\subsection{Formulation of the homological equations.}

Write
\begin{equation*}
f^{low}=f^{\varphi}+f^{0}+f^{1}+f^{2}=F^{\varphi}(\varphi;w)+\langle F_{0}(\varphi;w),r\rangle+\langle F_{1}(\varphi;w),\zeta\rangle+\frac{1}{2}\langle F_{2}(\varphi;w)\zeta,\zeta\rangle,
\end{equation*}
\begin{equation*}
s=s^{low}=s^{\varphi}+s^{0}+s^{1}+s^{2}=S^{\varphi}(\varphi;w)+\langle S_{0}(\varphi;w),r\rangle+\langle S_{1}(\varphi;w),\zeta\rangle+\frac{1}{2}\langle S_{2}(\varphi;w)\zeta,\zeta\rangle.
\end{equation*}
By the calculations  in \cite[Section 4.1.2]{CLY}, we obtain
\begin{equation}\label{es1}
\{f^{high},s\}^{low}=\{f^{high},s\}^{0}+\{f^{high},s\}^{1}+\{f^{high},s\}^{2},
\end{equation}
\begin{equation}\label{es2}
\{f^{high},s\}^{1}=\{f^{high},s^{\varphi}\}^{1},
\end{equation}
\begin{equation}\label{es3}
\{f^{high},s\}^{0}=\{f^{high},s^{\varphi}+s^{1}\}^{0},   \ \{f^{high},s\}^{2}=\{f^{high},s^{\varphi}+s^{1}\}^{2}.
\end{equation}
Let $g=\{f^{high},s\}$. Write
\begin{equation*}
g^{low}=g^{0}+g^{1}+g^{2}=\langle G_{0}(\varphi;w),r\rangle+\langle G_{1}(\varphi;w),\zeta\rangle+\frac{1}{2}\langle G_{2}(\varphi;w)\zeta,\zeta\rangle.
\end{equation*}

In Fourier modes, the homological equation (\ref{ho3}) decomposes into
\begin{equation}\label{es4}
-\mathrm{i}\langle k,\omega(w)\rangle\hat{S}^{\varphi}(k;w)=-\hat{F}^{\varphi}(k;w)+\delta_{0}^{k}a_{1}(w),
\end{equation}
\begin{equation}\label{es5}
-\mathrm{i}\langle k,\omega(w)\rangle\hat{S}_{1}(k;w)+J(\Omega(w)+H(w))\hat{S}_{1}(k;w)=-\hat{F}_{1}(k;w)-\hat{G}_{1}(k;w),
\end{equation}
\begin{equation}\label{es6}
-\mathrm{i}\langle k,\omega(w)\rangle\hat{S}_{0}(k;w)=-\hat{F}_{0}(k;w)-\hat{G}_{0}(k;w)+\delta_{0}^{k}\chi_{1}(w),
\end{equation}
\begin{equation}\label{es7}
\begin{array}{c}
-\mathrm{i}\langle k,\omega(w)\rangle\hat{S}_{2}(k;w)+(\Omega(w)+H(w))J\hat{S}_{2}(k;w)
-\hat{S}_{2}(k;w)J(\Omega(w)+H(w))\\
=-\hat{F}_{2}(k;w)-\hat{G}_{2}(k;w)+\delta_{0}^{k}H_{1}(w).
\end{array}
\end{equation}
We solve the equations (\ref{es4})-(\ref{es7}) in the order (\ref{es4}) $\rightarrow$ (\ref{es5}) $\rightarrow$ (\ref{es6}) $\rightarrow$ (\ref{es7}).

\subsection{Solution of the homological equation \eqref{es4}.}
The homological equation \eqref{es4} is very standard in the KAM theory.
From (\ref{es4}), we obtain
\begin{equation*}\label{es8}
a_{1}(w)=\hat{F}^{\varphi}(0;w)~\quad \textrm{and}\quad \hat{S}^{\varphi}(k;w)=\frac{\hat{F}^{\varphi}(k;w)}{\mathrm{i}\langle k,\omega(w)\rangle},~ k\neq0.
\end{equation*}
By the Diophantine condition (\ref{sd1}), we have
\begin{equation}\label{es9}
|\hat{S}^{\varphi}(k;w)|\leq\frac{1}{\kappa}|\hat{F}^{\varphi}(k;w)|.
\end{equation}
Differentiating (\ref{es4}), we obtain a similar homological equation
\begin{equation}\label{es10}
-\mathrm{i}\partial_{w}(\langle k,\omega(w)\rangle)\hat{S}^{\varphi}(k;w)-\mathrm{i}\langle k,\omega(w)\rangle\partial_{w}\hat{S}^{\varphi}(k;w)=-\partial_{w}\hat{F}^{\varphi}(k;w)
\end{equation}
for $\partial_{w}\hat{S}^{\varphi}(k;w)$
and there is also
\begin{equation*}\label{es11}
|\partial_{w}\hat{S}^{\varphi}(k;w)|\leq\frac{1}{\kappa}(|k|\cdot|\hat{S}^{\varphi}(k;w)|
+|\partial_{w}\hat{F}^{\varphi}(k;w)|),
\end{equation*}
which together with \eqref{es9} implies
\begin{equation*}\label{es12}
|\hat{S}^{\varphi}(k;w)|+|\partial_{w}\hat{S}^{\varphi}(k;w)|\leq\frac{|k|+1}{\kappa^{2}}(|\hat{F}^{\varphi}(k;w)|+|\partial_{w}\hat{F}^{\varphi}(k;w)|).
\end{equation*}
It follows that
\begin{equation*}\label{es13}
\begin{aligned}
\|S^{\varphi}_{\varphi}\|_{D(\rho-\tau)\times U'}&=\sum_{k\in\mathbb{Z}^{\mathcal{A}}}(|\hat{S}^{\varphi}_{\varphi}(k;w)|+|\partial_{w}\hat{S}^{\varphi}_{\varphi}(k;w)|)e^{|k|(\rho-\tau)}\\
 &\leq\sum_{k\in\mathbb{Z}^{\mathcal{A}}}\frac{|k|+1}{\kappa^{2}}(|\hat{F}^{\varphi}_{\varphi}(k;w)|+|\partial_{w}\hat{F}^{\varphi}_{\varphi}(k;w)|)e^{|k|(\rho-\tau)}\\
 &\preceq \frac{1}{\tau\kappa^{2}}\|F^{\varphi}_{\varphi}\|_{D(\rho)\times U'}.
\end{aligned}
\end{equation*}
As a result, we have
\begin{equation}\label{es14}
|||X_{s^{\varphi}}|||^{T}_{p,D(\rho-\tau,\sigma^{2},\sigma)\times U'}\preceq \frac{1}{\tau\kappa^{2}}|||X_{f^{\varphi}}|||^{T}_{p,D(\rho,\sigma^{2},\sigma)\times U'}.
\end{equation}

\bigskip

\subsection{Solution of the homological equation \eqref{es5}.}
 For simplicity, we write (\ref{es5}) as
\begin{equation*}\label{es15}
\mathrm{i}\langle k,\omega(w)\rangle S+J(\Omega+H)S=F+G.
\end{equation*}
Transforming into complex coordinates
\begin{equation*}
z=\left(
    \begin{array}{c}
      u \\
      v \\
    \end{array}
  \right)
  =C^{-1}\left(
           \begin{array}{c}
             \xi \\
             \eta \\
           \end{array}
         \right),
         C=\frac{1}{\sqrt{2}}\left(
         \begin{array}{cc}
           1 & 1 \\
           -\mathrm{i} & \mathrm{i} \\
         \end{array}
       \right),
\end{equation*}
and letting $S'=C^{-1}S=\left(
           \begin{array}{c}
             S'_{1} \\
             S'_{2} \\
           \end{array}
         \right),F'=C^{-1}F=\left(
           \begin{array}{c}
             F'_{1} \\
             F'_{2} \\
           \end{array}
         \right),G'=C^{-1}G=\left(
           \begin{array}{c}
             G'_{1} \\
             G'_{2} \\
           \end{array}
         \right)$,
we obtain the equivalent equations
\begin{equation}\label{es16}
\left\{ ~
\begin{aligned}
& \mathrm{i}\langle k,\omega(w)\rangle S'_{1}-\mathrm{i}(\Omega+H^{T})S'_{1}=F'_{1}+G'_{1},\\
& \mathrm{i}\langle k,\omega(w)\rangle S'_{2}+\mathrm{i}(\Omega+H)S'_{2}=F'_{2}+G'_{2}.
\end{aligned}
\right.
\end{equation}

\textbf{Solution of \eqref{es16}.}
We only solve $S'_{1}$ in \eqref{es16} since $S'_{2}$ can be solved accordingly.
By the first Melnikov condition (\ref{sd2}), we have
\begin{equation*}\label{es18}
\|S'_{1[a]_\Delta}\|\leq\frac{1}{\kappa}\|F'_{1[a]_\Delta}+G'_{1[a]_\Delta}\|.
\end{equation*}
Using similar arguments to that of $S^{\varphi}$, we get
\begin{equation}\label{es21}
\|S'_{1[a]_\Delta}\|+\|\partial_{w}S'_{1[a]_\Delta}\|\preceq\frac{|k|+1}{\kappa^{2}}(\|F'_{1[a]_\Delta}+G'_{1[a]_\Delta}\|+\|\partial_{w}F'_{1[a]_\Delta}+\partial_{w}G'_{1[a]_\Delta}\|).
\end{equation}

\smallskip

\textbf{Estimate of $s^{1}_{\varphi}$.}
Recall
\begin{equation*}\label{es22}
s^{1}=\sum_{a\in\mathcal{L}}(S'_{1a}(\varphi;w)u_{a}+S'_{2a}(\varphi;w)v_{a}),
\end{equation*}
and consider
\begin{equation*}\label{es23}
\lfloor s^{1}_{\varphi}\rceil_{D(\rho-3\tau,\sigma^{2})\times U'}(z)=\sum_{a\in\mathcal{L}}(\|S'_{1a\varphi}(\varphi;w)\|_{D(\rho-3\tau)\times U'}u_{a}+\|S'_{2a\varphi}(\varphi;w)\|_{D(\rho-3\tau)\times U'}v_{a}).
\end{equation*}
For $z=(u,v)$, define $\tilde{z}=(\tilde{u},\tilde{v})$ such that for all $a\in [a]_\Delta$, $\tilde{u}_{a}=\|u_{[a]_\Delta}\|,\tilde{v}_{a}=\|v_{[a]_\Delta}\|$.
By (\ref{es21}), we see the first sum in $\lfloor s^{1}_{\varphi}\rceil$ satisfies
\begin{equation*}\label{es24}
\begin{aligned}
&~ |\sum_{a\in\mathcal{L}}\|S'_{1a\varphi}(\varphi)\|_{D(\rho-3\tau)\times U'}u_{a}|\\
\leq &
\sum_{a\in\mathcal{L}}\sum_{k\in\mathbb{Z}^{\mathcal{A}}}
(|\hat{S}'_{1a\varphi}(k)|+|\partial_{w}\hat{S}'_{1a\varphi}(k)|) \cdot e^{(\rho-3\tau)|k|}\cdot |u_{a}|\\
= & \sum_{k\in\mathbb{Z}^{\mathcal{A}}}\sum_{[a]_\Delta\in\mathcal{E}_{\Delta}}
\sum_{a\in[a]_\Delta}(|\hat{S}'_{1a\varphi}(k)|+|\partial_{w}\hat{S}'_{1a\varphi}(k)|)
\cdot |u_{a}| \cdot e^{(\rho-3\tau)|k|}\\
\preceq & \sum_{k\in\mathbb{Z}^{\mathcal{A}}}\sum_{a\in\mathcal{L}}\frac{|k|+1}{\kappa^{2}}
\left(|\hat{F}'_{1a\varphi}(k)+\hat{G}'_{1a\varphi}(k)|
+|\partial_{w}\hat{F}'_{1a\varphi}(k)+\partial_{w}\hat{G}'_{1a\varphi}(k)|\right)
\tilde{u}_{a} e^{(\rho-3\tau)|k|}\\
\preceq &\frac{1}{\tau\kappa^{2}}\sum_{a\in\mathcal{L}}\|F'_{1a\varphi}(\varphi)
+G'_{1a\varphi}(\varphi)\|_{D(\rho-2\tau)\times U'} \cdot \tilde{u}_{a}.
\end{aligned}
\end{equation*}
Similarly, we have
\begin{equation*}\label{es25}
|\sum_{a\in\mathcal{L}}\|S'_{2a\varphi}(\varphi;w)\|_{D(\rho-3\tau)\times U'}v_{a}|
\preceq\frac{1}{\tau\kappa^{2}}\sum_{a\in\mathcal{L}}\|F'_{2a\varphi}(\varphi;w)+G'_{2a\varphi}(\varphi;w)\|_{D(\rho-2\tau)\times U'}\tilde{v}_{a}.
\end{equation*}
It follows that
\begin{equation*}\label{es26}
|\lfloor s^{1}_{\varphi}\rceil_{D(\rho-3\tau,\sigma^{2})\times U'}(z)|\preceq\frac{1}{\tau\kappa^{2}}\lfloor f^{1}_{\varphi}+g^{1}_{\varphi}\rceil_{D(\rho-2\tau,\sigma^{2})\times U'}(\tilde{z}).
\end{equation*}
Moreover, we see from
\begin{equation*}\label{es27}
|||s^{1}_{\varphi}|||_{D(\rho-3\tau,\sigma^{2})\times U'}=\sup_{0\neq z\in l_{1}^{2}}\frac{|\lfloor s^{1}_{\varphi}\rceil_{D(\rho-3\tau,\sigma^{2})\times U'}(z)|}{\|z\|_{1}},
\end{equation*}
and  $\|\tilde{z}\|_{1}\preceq d_{\Delta}^{d}\|z\|_{1}$
that
\begin{equation}\label{es28}
|||s^{1}_{\varphi}|||_{D(\rho-3\tau,\sigma^{2})\times U'}\preceq\frac{d_{\Delta}^{d}}{\tau\kappa^{2}}|||f^{1}_{\varphi}+g^{1}_{\varphi}|||_{D(\rho-2\tau,\sigma^{2})\times U'}.
\end{equation}

\smallskip

\textbf{Estimate of $s^{1}_{z}$.}
Consider
\begin{equation*}\label{es29}
|||s^{1}_{z}|||^{T}_{p,D(\rho-3\tau,\sigma^{2})\times U'}=\|\lfloor s^{1}_{z}\rceil_{D(\rho-3\tau,\sigma^{2})\times U'}\|_{p},
\end{equation*}
\begin{equation*}\label{es30}
\|\lfloor s^{1}_{z}\rceil_{D(\rho-3\tau,\sigma^{2})\times U'}\|_{p}^{2}=\sum_{a\in\mathcal{L}}(\|S'_{1a}(\varphi)\|^{2}_{D(\rho-3\tau)\times U'}+\|S'_{2a}(\varphi)\|^{2}_{D(\rho-3\tau)\times U'})\langle a\rangle^{2p}.
\end{equation*}
By (\ref{es21}), we have
\begin{equation*}\label{es31}
\left(\sum_{a\in[a]_\Delta}\|S'_{1a}(\varphi)\|^{2}_{D(\rho-3\tau)\times U'}\right)^{\frac{1}{2}}=
\left[\sum_{a\in[a]_\Delta}(\sum_{k\in\mathbb{Z}^{\mathcal{A}}}(|\hat{S}'_{1a}(k)|
+|\partial_{w}\hat{S}'_{1a}(k)|)e^{(\rho-3\tau)|k|})^{2}\right]^{\frac{1}{2}}
\end{equation*}
\begin{align*}
&\leq\sum_{k\in\mathbb{Z}^{\mathcal{A}}}[\sum_{a\in[a]_\Delta}(|\hat{S}'_{1a}(k)|
+|\partial_{w}\hat{S}'_{1a}(k)|)^{2}]^{\frac{1}{2}}e^{(\rho-3\tau)|k|}\\
&\preceq\sum_{k\in\mathbb{Z}^{\mathcal{A}}}\frac{|k|+1}{\kappa^{2}}
\left(\|\hat{F}'_{1[a]}(k)+\hat{G}'_{1[a]}(k)\|+\|\partial_{w}\hat{F}'_{1[a]}(k)
+\partial_{w}\hat{G}'_{1[a]}(k)\|\right)e^{(\rho-3\tau)|k|}\\
&\preceq\sum_{k\in\mathbb{Z}^{\mathcal{A}}}\sum_{a\in[a]_\Delta}\frac{|k|+1}{\kappa^{2}}
\left(|\hat{F}'_{1a}(k)+\hat{G}'_{1a}(k)|+|\partial_{w}\hat{F}'_{1a}(k)
+\partial_{w}\hat{G}'_{1a}(k)|\right)e^{(\rho-3\tau)|k|}\\
&\preceq\frac{1}{\tau\kappa^{2}}\sum_{a\in[a]_\Delta}\|F'_{1a}(\varphi)
+G'_{1a}(\varphi)\|_{D(\rho-2\tau)\times U'},
\end{align*}
which implies the first sum in $\lfloor s^{1}_{z}\rceil$ satisfies
\begin{equation*}\label{es32}
\sum_{a\in[a]_\Delta}\|S'_{1a}(\varphi;w)\|^{2}_{D(\rho-3\tau)\times U'}\langle a\rangle^{2p}
\preceq(\frac{d_{\Delta}^{d}}{\tau\kappa^{2}})^{2}\sum_{a\in[a]_\Delta}\|F'_{1a}(\varphi;w)+G'_{1a}(\varphi;w)\|^{2}_{D(\rho-2\tau)\times U'}\langle a\rangle^{2p}.
\end{equation*}
Similarly, the other sum in $\lfloor s^{1}_{z}\rceil$ satisfies
\begin{equation*}\label{es33}
\sum_{a\in[a]_\Delta}\|S'_{2a}(\varphi;w)\|^{2}_{D(\rho-3\tau)\times U'}\langle a\rangle^{2p}
\preceq(\frac{d_{\Delta}^{d}}{\tau\kappa^{2}})^{2}\sum_{a\in[a]_\Delta}\|F'_{2a}(\varphi;w)+G'_{2a}(\varphi;w)\|^{2}_{D(\rho-2\tau)\times U'}\langle a\rangle^{2p}.
\end{equation*}
Then we immediately get
\begin{equation}\label{es34}
|||s^{1}_{z}|||^{T}_{p,D(\rho-3\tau,\sigma^{2})\times U'}\preceq\frac{d_{\Delta}^{d}}{\tau\kappa^{2}}|||f^{1}_{z}+g^{1}_{z}|||^{T}_{p,D(\rho-2\tau,\sigma^{2})\times U'}.
\end{equation}

Combining (\ref{es28}) and (\ref{es34}), we have
\begin{equation}\label{es35}
|||X_{s^{1}}|||^{T}_{p,D(\rho-3\tau,\sigma^{2},\sigma)\times U'}\preceq \frac{d_{\Delta}^{d}}{\tau\kappa^{2}}|||X_{f^{1}+g^{1}}|||^{T}_{p,D(\rho-2\tau,\sigma^{2},\sigma)\times U'}.
\end{equation}

\bigskip

\subsection{Solution of the homological equations \eqref{es6}-\eqref{es7}.}
Solving  equation (\ref{es6}) as  equation (\ref{es4}), we obtain
\begin{equation}\label{es36}
|||X_{s^{0}}|||^{T}_{p,D(\rho-5\tau,\sigma^{2},\sigma)\times U'}\preceq \frac{1}{\tau\kappa^{2}}|||X_{f^{0}+g^{0}}|||^{T}_{p,D(\rho-4\tau,\sigma^{2},\sigma)\times U'}.
\end{equation}

Now we consider the equation (\ref{es7}). For simplicity, we write (\ref{es7}) as
\begin{equation*}\label{es37}
\mathrm{i}\langle k,\omega(w)\rangle S+(\Omega+H)JS-SJ(\Omega+H)=F+G-H_{1}.
\end{equation*}
Changing into complex coordinates $z=C^{-1}\zeta$ and
letting $S'=C^{T}SC=\left(
                  \begin{array}{cc}
                    S_{1}' &  S_{2}' \\
                    S'^{T}_{2} & S_{3}' \\
                  \end{array}
                \right)
,F'=C^{T}FC,G'=C^{T}GC$, $H_{1}'=C^{T}H_{1}C$, we obtain the equivalent equations as follows
\begin{equation}\label{es38}
\mathrm{i}\langle k,\omega(w)\rangle S_{1}'+\mathrm{i}(\Omega+H)S_{1}'+\mathrm{i}S_{1}'(\Omega+H^{T})=F'_{1}+G'_{1},
\end{equation}
\begin{equation}\label{es39}
\mathrm{i}\langle k,\omega(w)\rangle S_{2}'+\mathrm{i}(\Omega+H)S_{2}'-\mathrm{i}S_{2}'(\Omega+H)=F'_{2}+G'_{2}-H_{12}',
\end{equation}
\begin{equation}\label{es40}
\mathrm{i}\langle k,\omega(w)\rangle S'^{T}_{2}-\mathrm{i}(\Omega+H^{T})S'^{T}_{2}+\mathrm{i}S'^{T}_{2}(\Omega+H^{T})=F'^{T}_{2}+G'^{T}_{2}-H'^{T}_{12},
\end{equation}
\begin{equation}\label{es41}
\mathrm{i}\langle k,\omega(w)\rangle S_{3}'-\mathrm{i}(\Omega+H^{T})S_{3}'-\mathrm{i}S_{3}'(\Omega+H)=F'_{3}+G'_{3}.
\end{equation}

\textbf{Solutions of \eqref{es38}-\eqref{es41}.}
Consider first $S'_{2}$ in \eqref{es39}-\eqref{es40}.
When $k\neq0$, we have $H_{12}'=0$. In a similar way as before, we obtain from
the second Melnikov condition \eqref{sd4} that
\begin{equation*}\label{es45}
\|S'^{[b]_\Delta}_{2[a]_\Delta}\|+\|\partial_{w}S'^{[b]_\Delta}_{2[a]_\Delta}\|
\preceq\frac{|k|+1}{\kappa^{2}}(\|F'^{[b]_\Delta}_{2[a]_\Delta}+G'^{[b]_\Delta}_{2[a]_\Delta}\|
+\|\partial_{w}F'^{[b]_\Delta}_{2[a]_\Delta}+\partial_{w}G'^{[b]_\Delta}_{2[a]_\Delta}\|).
\end{equation*}
 When $k=0$, the above estimate  also holds since $H_{12}'=(F'_{2}+G'_{2})\mid_{\{(a,b)\in\mathcal{L}\times\mathcal{L}:|a|=|b|\}}$.
Using the second Melnikov condition (\ref{sd3}) with the same sign, we can
also solve $S_{1}'$ and $S_{3}'$. In conclusion, we have
\begin{equation}\label{es46}
\|S'^{[b]_\Delta}_{\nu [a]_\Delta}\|+\|\partial_{w}S'^{[b]_\Delta}_{\nu [a]_\Delta}\|
\preceq\frac{|k|+1}{\kappa^{2}}(\|F'^{[b]_\Delta}_{\nu [a]_\Delta}+G'^{[b]_\Delta}_{\nu [a]_\Delta}\|
+\|\partial_{w}F'^{[b]_\Delta}_{\nu [a]_\Delta}+\partial_{w}G'^{[b]_\Delta}_{\nu [a]_\Delta}\|),
\end{equation}
with $\nu\in\{1, 2, 3\}$.

\textbf{Estimate of $s^{2}_{\varphi}$.}
Recalling
\begin{equation*}\label{es48}
s^{2}=\frac{1}{2}\sum_{a,b\in\mathcal{L}}(S'^{b}_{1a}(\varphi;w)u_{a}u_{b}+2S'^{b}_{2a}(\varphi;w)u_{a}v_{b}+S'^{b}_{3a}(\varphi;w)v_{a}v_{b}),
\end{equation*}
we consider
\begin{equation*}\label{es49}
\begin{aligned}
\lfloor s^{2}_{\varphi}\rceil_{D(\rho-5\tau,\sigma^{2})\times U'}(z)=&\frac{1}{2}
\sum_{a,b\in\mathcal{L}}(\|S'^{b}_{1a\varphi}(\varphi;w)\|_{D(\rho-5\tau)\times U'}u_{a}u_{b}\\
+&2\|S'^{b}_{2a\varphi}(\varphi;w)\|_{D(\rho-5\tau)\times U'}u_{a}v_{b}+\|S'^{b}_{3a\varphi}(\varphi;w)\|_{D(\rho-5\tau)\times U'}v_{a}v_{b}),
\end{aligned}
\end{equation*}
and and the associated multilinear form $\widetilde{s^{2}_{\varphi}}$
\begin{equation}\label{es50}
\begin{aligned}
\widetilde{\lfloor s^{2}_{\varphi}\rceil}&_{D(\rho-5\tau,\sigma^{2})\times U'}(z^{(1)},z^{(2)})=\frac{1}{2}
\sum_{a,b\in\mathcal{L}}(\|S'^{b}_{1a\varphi}(\varphi;w)\|_{D(\rho-5\tau)\times U'}u^{(1)}_{a}u^{(2)}_{b}\\
+&\|S'^{b}_{2a\varphi}(\varphi;w)\|_{D(\rho-5\tau)\times U'}(u^{(1)}_{a}v^{(2)}_{b}+u^{(2)}_{a}v^{(1)}_{b})+\|S'^{b}_{3a\varphi}(\varphi;w)\|_{D(\rho-5\tau)\times U'}v^{(1)}_{a}v^{(2)}_{b}).
\end{aligned}
\end{equation}
By (\ref{es46}), we know that
$|\sum_{a\in[a]_\Delta,b\in[b]_\Delta}\|S'^{b}_{1a\varphi}(\varphi;w)\|_{D(\rho-5\tau)\times U'}u^{(1)}_{a}u^{(2)}_{b}|$
is less than
\begin{align*}
&\sum_{a\in[a],b\in[b]}\sum_{k\in\mathbb{Z}^{\mathcal{A}}}(|\hat{S}'^{b}_{1a\varphi}(k)|+|\partial_{w}\hat{S}'^{b}_{1a\varphi}(k)|)e^{(\rho-5\tau)|k|}|u^{(1)}_{a}u^{(2)}_{b}| \\
&\preceq\sum_{k\in\mathbb{Z}^{\mathcal{A}}}\sum_{a\in[a],b\in[b]}\frac{|k|+1}{\kappa^{2}}
\left(|\hat{F}'^{b}_{1a\varphi}(k)+\hat{G}'^{b}_{1a\varphi}(k)|
+|\partial_{w}\hat{F}'^{b}_{1a\varphi}(k)+\partial_{w}\hat{G}'^{b}_{1a\varphi}(k)|\right)\tilde{u}^{(1)}_{a}\tilde{u}^{(2)}_{b}e^{(\rho-5\tau)|k|}\\
&\preceq\frac{1}{\tau\kappa^{2}}\sum_{a\in[a],b\in[b]}\|F'^{b}_{1a\varphi}(\varphi)+G'^{b}_{1a\varphi}(\varphi)\|_{D(\rho-4\tau)\times U'}\tilde{u}^{(1)}_{a}\tilde{u}^{(2)}_{b},
\end{align*}
which implies
\begin{equation*}\label{es52}
|\sum_{a,b\in\mathcal{L}}\|S'^{b}_{1a\varphi}(\varphi)\|_{D(\rho-5\tau)\times U'}u^{(1)}_{a}u^{(2)}_{b}|
\preceq\frac{1}{\tau\kappa^{2}}\sum_{a,b\in\mathcal{L}}\|F'^{b}_{1a\varphi}(\varphi)+G'^{b}_{1a\varphi}(\varphi)\|_{D(\rho-4\tau)\times U'}\tilde{u}^{(1)}_{a}\tilde{u}^{(2)}_{b}.
\end{equation*}
There are similar estimate for the other three summations in the R.H.S. of \eqref{es50}.

Then we have
\begin{equation*}\label{es56}
|\widetilde{\lfloor s^{2}_{\varphi}\rceil}_{D(\rho-5\tau,\sigma^{2})\times U'}(z^{(1)},z^{(2)})|\preceq\frac{1}{\tau\kappa^{2}}
\widetilde{\lfloor f^{2}_{\varphi}+g^{2}_{\varphi}\rceil}_{D(\rho-4\tau,\sigma^{2})\times U'}(\tilde{z}^{(1)},\tilde{z}^{(2)}).
\end{equation*}
Since
\begin{equation*}\label{es57}
|||s^{2}_{\varphi}|||_{D(\rho-5\tau,\sigma^{2})\times U'}
=\sup_{0\neq z^{(1)},z^{(2)}\in l_{1}^{2}}\frac{|\widetilde{\lfloor s^{2}_{\varphi}\rceil}_{D(\rho-5\tau,\sigma^{2})\times U'}(z^{(1)},z^{(2)})|}{\|z^{(1)}\|_{1}\|z^{(2)}\|_{1}},
\end{equation*}
we see from  $\|\tilde{z}\|_{1}\preceq d_{\Delta}^{d}\|z\|_{1}$ that
\begin{equation}\label{es58}
|||s^{2}_{\varphi}|||_{D(\rho-5\tau,\sigma^{2})\times U'}\preceq\frac{d_{\Delta}^{2d}}{\tau\kappa^{2}}|||f^{2}_{\varphi}+g^{2}_{\varphi}|||_{D(\rho-4\tau,\sigma^{2})\times U'}.
\end{equation}

\textbf{Estimate of $s^{2}_{z}$.}
Now we estimate
\begin{equation*}\label{es59}
|||s^{2}_{z}|||^{T}_{p,D(\rho-5\tau,\sigma^{2})\times U'}=\sup_{0\neq z\in l_{p}^{2}}\frac{\|\lfloor s^{2}_{z}\rceil_{D(\rho-5\tau,\sigma^{2})\times U'}(z)\|_{p}}{\|z\|_{p}},
\end{equation*}
in which $\|\lfloor s^{2}_{z}\rceil_{D(\rho-5\tau,\sigma^{2})\times U'}(z)\|_{p}^{2}$ equals to the following sum
\begin{equation}\label{es60}
\begin{aligned}
& \sum_{a\in\mathcal{L}}
\left|\sum_{b\in\mathcal{L}}(\|S'^{b}_{1a}(\varphi;w)\|_{D(\rho-5\tau)\times U'}u_{b}+\|S'^{b}_{2a}(\varphi;w)\|_{D(\rho-5\tau)\times U'}v_{b})\right|^{2}\langle a\rangle^{2p}\\
+& \sum_{a\in\mathcal{L}}\left|\sum_{b\in\mathcal{L}}(\|S'^{a}_{2b}(\varphi;w)\|_{D(\rho-5\tau)\times U'}u_{b}+\|S'^{b}_{3a}(\varphi;w)\|_{D(\rho-5\tau)\times U'}v_{b})\right|^{2}\langle a\rangle^{2p}.
\end{aligned}
\end{equation}
By (\ref{es46}), we see that
\begin{equation*}\label{es61}
\left[\sum_{a\in[a]_\Delta}\left(\sum_{b\in\mathcal{L}}\|S'^{b}_{1a}(\varphi;w)\|_{D(\rho-5\tau)\times U'}|u_{b}|\right)^{2}\right]^{\frac{1}{2}}
\end{equation*}
\begin{align*}
&=\left[\sum_{a\in[a]}\left(\sum_{k\in\mathbb{Z}^{\mathcal{A}}}\sum_{b\in\mathcal{L}}(|\hat{S}'^{b}_{1a}(k)|+|\partial_{w}\hat{S}'^{b}_{1a}(k)|)e^{(\rho-5\tau)|k|}|u_{b}|\right)^{2}\right]^{\frac{1}{2}} \\
&\leq\sum_{k\in\mathbb{Z}^{\mathcal{A}}}\left[\sum_{a\in[a]}\left(\sum_{b\in\mathcal{L}}(|\hat{S}'^{b}_{1a}(k)|+|\partial_{w}\hat{S}'^{b}_{1a}(k)|)e^{(\rho-5\tau)|k|}|u_{b}|\right)^{2}\right]^{\frac{1}{2}}\\
&\leq\sum_{k\in\mathbb{Z}^{\mathcal{A}}}\left[\sum_{b\in\mathcal{L}}\left(\sum_{a\in[a]}(|\hat{S}'^{b}_{1a}(k)|+|\partial_{w}\hat{S}'^{b}_{1a}(k)|)^{2}\right)^{\frac{1}{2}}|u_{b}|\right]e^{(\rho-5\tau)|k|}\\
&\preceq\sum_{k\in\mathbb{Z}^{\mathcal{A}}}\frac{|k|+1}{\kappa^{2}}e^{(\rho-5\tau)|k|}\sum_{[b]\in\mathcal{E}_{\Delta}}(\|\hat{F}'^{[b]}_{1[a]}(k)+\hat{G}'^{[b]}_{1[a]}(k)\|
+\|\partial_{w}\hat{F}'^{[b]}_{1[a]}(k)+\partial_{w}\hat{G}'^{[b]}_{1[a]}(k)\|)\|u_{[b]}\| \\
&\preceq\sum_{k\in\mathbb{Z}^{\mathcal{A}}}\frac{|k|+1}{\kappa^{2}}e^{(\rho-5\tau)|k|}\sum_{a\in[a]}\sum_{b\in\mathcal{L}}(|\hat{F}'^{b}_{1a}(k)+\hat{G}'^{b}_{1a}(k)|
+|\partial_{w}\hat{F}'^{b}_{1a}(k)+\partial_{w}\hat{G}'^{b}_{1a}(k)|)\tilde{u}_{b}\\
&\preceq\frac{1}{\tau\kappa^{2}}\sum_{a\in[a]}\sum_{b\in\mathcal{L}}\|F'^{b}_{1a}(\varphi)+G'^{b}_{1a}(\varphi)\|_{D(\rho-4\tau)\times U'}\tilde{u}_{b}.
\end{align*}
As a result, we have
\begin{equation*}\label{es62}
\begin{aligned}
\sum_{a\in[a]_\Delta}&\left(\sum_{b\in\mathcal{L}}\|S'^{b}_{1a}(\varphi)\|_{D(\rho-5\tau)\times U'}|u_{b}|\right)^{2}\\
&\preceq(\frac{d^{d}_{\Delta}}{\tau\kappa^{2}})^{2}\sum_{a\in[a]_{\Delta}}\left(\sum_{b\in\mathcal{L}}\|F'^{b}_{1a}(\varphi)+G'^{b}_{1a}(\varphi)\|_{D(\rho-4\tau)\times U'}\tilde{u}_{b}\right)^{2}.
\end{aligned}
\end{equation*}
Similar estimates hold for the other three summations in \eqref{es60}.
Then we have
\begin{equation*}\label{es66}
\|\lfloor s^{2}_{z}\rceil_{D(\rho-5\tau,\sigma^{2})\times U'}(z)\|_{p}\preceq\frac{d^{d}_{\Delta}}{\tau\kappa^{2}}\|\lfloor f^{2}_{z}+g^{2}_{z}\rceil_{D(\rho-4\tau,\sigma^{2})\times U'}(\tilde{z})\|_{p},
\end{equation*}
which together with $\|\tilde{z}\|_{p}\preceq d_{\Delta}^{d}\|z\|_{p}$ implies
\begin{equation}\label{es67}
|||s^{2}_{z}|||^{T}_{D(\rho-5\tau,\sigma^{2})\times U'}\preceq\frac{d_{\Delta}^{2d}}{\tau\kappa^{2}}|||f^{2}_{z}+g^{2}_{z}|||^{T}_{D(\rho-4\tau,\sigma^{2})\times U'}.
\end{equation}

Finally, combining  (\ref{es58}) and  (\ref{es67}), we have
\begin{equation*}\label{es68}
|||X_{s^{2}}|||^{T}_{p,D(\rho-5\tau,\sigma^{2},\sigma)\times U'}\preceq \frac{d_{\Delta}^{2d}}{\tau\kappa^{2}}|||X_{f^{2}+g^{2}}|||^{T}_{p,D(\rho-4\tau,\sigma^{2},\sigma)\times U'}.
\end{equation*}

\smallskip

\subsection{Verification of the estimates \eqref{ho6}-\eqref{ho11}.}
In this part, we shall verify the estimates of the vector fields $X_{s}$ and $X_{h_{1}}$.

Recall that $s=s^{\varphi}+s^{0}+s^{1}+s^{2}$.
By (\ref{ho1}), (\ref{es14}), we have
\begin{equation*}\label{es69}
|||X_{s^{\varphi}}|||^{T}_{p,D(\rho-\tau,\sigma^{2},\sigma)\times U'}\preceq \frac{1}{\tau\kappa^{2}}|||X_{f^{\varphi}}|||^{T}_{p,D(\rho,\sigma^{2},\sigma)\times U'}\preceq \frac{\varepsilon}{\tau\kappa^{2}},
\end{equation*}
which together with
Proposition \ref{p} and  (\ref{ho1}) implies
\begin{equation}\label{es70}
\begin{aligned}
|||X_{\{f^{high},s^{\varphi}\}}|||^{T}&_{p,D(\rho-2\tau,(\sigma-\tau)^{2},\sigma-\tau)\times U'}\\
&\preceq
\frac{1}{\tau}|||X_{f^{high}}|||^{T}_{p,D(\rho-\tau,\sigma^{2},\sigma)\times U'}|||X_{s^{\varphi}}|||^{T}_{p,D(\rho-\tau,\sigma^{2},\sigma)\times U'}\preceq\frac{\varepsilon}{\tau^{2}\kappa^{2}}.
\end{aligned}
\end{equation}
By (\ref{ho1}), (\ref{es2}), (\ref{es35}) and (\ref{es70}), we have
\begin{equation*}\label{es71}
\begin{aligned}
&|||X_{s^{1}}|||^{T}_{p,D(\rho-3\tau,(\sigma-3\tau)^{2},\sigma-3\tau)\times U'}\preceq \frac{d_{\Delta}^{d}}{\tau\kappa^{2}}|||X_{f^{1}+g^{1}}|||^{T}_{p,D(\rho-2\tau,(\sigma-3\tau)^{2},\sigma-3\tau)\times U'}\\
\preceq &\frac{d_{\Delta}^{d}}{\tau\kappa^{2}}(|||X_{f^{1}}|||^{T}_{p,D(\rho-2\tau,(\sigma-3\tau)^{2},\sigma-3\tau)\times U'}+|||X_{g^{1}}|||^{T}_{p,D(\rho-2\tau,(\sigma-3\tau)^{2},\sigma-3\tau)\times U'})
\preceq\frac{d_{\Delta}^{d}\varepsilon}{\tau^{3}\kappa^{4}}.
\end{aligned}
\end{equation*}

Similar to $X_{s^{1}}$, we can estimate $X_{s^{0}}$ and $X_{s^{2}}$ in sequence and finally get
\begin{equation}\label{es75}
|||X_{s}|||^{T}_{p,D(\rho-5\tau,(\sigma-5\tau)^{2},\sigma-5\tau)\times U'}\preceq\frac{d_{\Delta}^{3d}\varepsilon}{\tau^{5}\kappa^{6}}.
\end{equation}
Moreover, the vector field $X_{h_{1}}$ satisfies
\begin{equation}\label{es76}
\begin{aligned}
|||X_{h_{1}}|||^{T}&_{p,D(\rho-5\tau,(\sigma-5\tau)^{2},\sigma-5\tau)\times U'}
\leq|||X_{f^{0}+g^{0}}|||^{T}_{p,D(\rho-5\tau,(\sigma-5\tau)^{2},\sigma-5\tau)\times U'}\\
&+|||X_{f^{2}+g^{2}}|||^{T}_{p,D(\rho-5\tau,(\sigma-5\tau)^{2},\sigma-5\tau)\times U'}
\preceq\frac{d_{\Delta}^{d}\varepsilon}{\tau^{4}\kappa^{4}}.
\end{aligned}
\end{equation}

\smallskip

\subsection{Estimate of the new Hamiltonian.} In this part, we shall verify the
properties \eqref{ho12}-\eqref{ho15}.
Using Taylor's formula,
we obtain from the homological equation (\ref{ho3}) that
\begin{equation*}\label{es77}
(h+f)\circ X^{t}_{s}\mid_{t=1}=(h+f^{low}+f^{high})\circ X^{t}_{s}\mid_{t=1}
\end{equation*}
\begin{align*}
=&h+\{h,s\}+\int_{0}^{1}(1-t)\{\{h,s\},s\}\circ X^{t}_{s}dt+f^{low}+\int_{0}^{1}\{f^{low},s\}\circ X^{t}_{s}dt \\
&+f^{high}+\{f^{high},s\}+\int_{0}^{1}(1-t)\{\{f^{high},s\},s\}\circ X^{t}_{s}dt \\
=&h+h_{1}+(1-\mathcal{T}_{\Delta'})f^{low}+f^{high}+(1-\mathcal{T}_{\Delta'})\{f^{high},s\}^{low}+\{f^{high},s\}^{high} \\
&+\int_{0}^{1}(1-t)\{\{h,s\},s\}\circ X^{t}_{s}dt+\int_{0}^{1}\{f^{low},s\}\circ X^{t}_{s}dt+\int_{0}^{1}(1-t)\{\{f^{high},s\},s\}\circ X^{t}_{s}dt.
\end{align*}
Then we have $f_{1}=f_{1}^{low}+f_{1}^{high}$, where
\begin{equation}\label{es113}
\begin{aligned}
f_{1}^{low}=&(1-\mathcal{T}_{\Delta'})f^{low}+(1-\mathcal{T}_{\Delta'})\{f^{high},s\}^{low}+\left(\int_{0}^{1}(1-t)\{\{h,s\},s\}\circ X^{t}_{s}dt\right)^{low}\\
&+\left(\int_{0}^{1}\{f^{low},s\}\circ X^{t}_{s}dt\right)^{low}+\left(\int_{0}^{1}(1-t)\{\{f^{high},s\},s\}\circ X^{t}_{s}dt\right)^{low},\\
 f_{1}^{high}=&f^{high}+\{f^{high},s\}^{high}+\left(\int_{0}^{1}(1-t)\{\{h,s\},s\}\circ X^{t}_{s}dt\right)^{high}\\
&+\left(\int_{0}^{1}\{f^{low},s\}\circ X^{t}_{s}dt\right)^{high}+\left(\int_{0}^{1}(1-t)\{\{f^{high},s\},s\}\circ X^{t}_{s}dt\right)^{high}.
\end{aligned}
\end{equation}

\bigskip

\textbf{Estimate of the vector field generated by $(1-\mathcal{T}_{\Delta'})f^{low}$.}
Observing that
\begin{equation*}\label{es78}
\begin{aligned}
\|(1-\mathcal{T}_{\Delta'})&F^{\varphi}_{\varphi}\|_{D(\rho-\tau)\times U'}=\sum_{|k|>\Delta'}(|\hat{F}^{\varphi}_{\varphi}(k;w)|+|\partial_{w}\hat{F}^{\varphi}_{\varphi}(k;w)|)e^{|k|(\rho-\tau)}\\
&\leq\sum_{|k|>\Delta'}e^{-\tau|k|}\|F^{\varphi}_{\varphi}\|_{D(\rho)\times U'}
\preceq \frac{1}{\tau^{\# \mathcal{A}}}e^{-\frac{1}{2}\tau\Delta'}\|F^{\varphi}_{\varphi}\|_{D(\rho)\times U'},
\end{aligned}
\end{equation*}
we obtain
\begin{equation*}\label{es79}
|||X_{(1-\mathcal{T}_{\Delta'})f^{\varphi}}|||^{T}_{p,D(\rho-\tau,\sigma^{2},\sigma)\times U'}
\preceq \frac{1}{\tau^{\# \mathcal{A}}}e^{-\frac{1}{2}\tau\Delta'}|||X_{f^{\varphi}}|||^{T}_{p,D(\rho,\sigma^{2},\sigma)\times U'}\preceq\frac{1}{\tau^{\# \mathcal{A}}}e^{-\frac{1}{2}\tau\Delta'}\varepsilon.
\end{equation*}
The same estimates hold for $|||X_{(1-\mathcal{T}_{\Delta'})f^{0}}|||^{T}_{p,D(\rho-\tau,\sigma^{2},\sigma)\times U'}$
and $|||X_{(1-\mathcal{T}_{\Delta'})f^{1}}|||^{T}_{p,D(\rho-\tau,\sigma^{2},\sigma)\times U'}$.
Then we turn to
\begin{equation*}\label{es83}
(1-\mathcal{T}_{\Delta'})f^{2}=(1-\mathcal{T}^{1}_{\Delta'})f^{2}+(1-\mathcal{T}^{2}_{\Delta'})f^{2},
\end{equation*}
where
\begin{equation*}\label{es82}
\begin{aligned}
(1-\mathcal{T}^{1}_{\Delta'})f^{2}=\frac{1}{2}\langle \zeta,(1-\mathcal{T}_{\Delta'})F_{2}(\varphi;w)\zeta\rangle,\quad
(1-\mathcal{T}^{2}_{\Delta'})f^{2}=\frac{1}{2}\sum_{|k|>\Delta'}\langle \zeta,\mathcal{T}_{\Delta'}\hat{F}_{2}(k;w)\zeta\rangle e^{\mathrm{i}\langle k,\varphi\rangle}.
\end{aligned}
\end{equation*}
It is easy to see
\begin{equation}\label{es84}
|||X_{(1-\mathcal{T}^{2}_{\Delta'})f^{2}}|||^{T}_{p,D(\rho-\tau,\sigma^{2},\sigma)\times U'}\preceq\frac{1}{\tau^{\# \mathcal{A}}}e^{-\frac{1}{2}\tau\Delta'}\varepsilon.
\end{equation}

By (\ref{ho2}), we have
\begin{equation*}\label{es85}
\sup_{(\varphi,w)\in D(\rho)\times U}(|F_{2}(\varphi;w)|_{\gamma},|\partial_{w}F_{2}(\varphi;w)|_{\gamma})\leq\frac{\varepsilon}{\sigma^{2}}.
\end{equation*}
Hence
\begin{equation*}\label{es86}
|\hat{F}'^{b}_{2a}(k;w)|+|\partial_{w}\hat{F}'^{b}_{2a}(k;w)|\preceq \frac{\varepsilon}{\sigma^{2}}e^{-\gamma|a-b|-\rho|k|},
\end{equation*}
which implies
\begin{equation*}\label{es87}
\begin{aligned}
\|F'^{b}_{2a\varphi}(\varphi;w)\|_{D(\rho-\tau)\times U}=&\sum_{k\in\mathbb{Z}^{\mathcal{A}}}(|\hat{F}'^{b}_{2a\varphi}(k;w)|
+|\partial_{w}\hat{F}'^{b}_{2a\varphi}(k;w)|)e^{(\rho-\tau)|k|}\\
\preceq& \sum_{k\in\mathbb{Z}^{\mathcal{A}}}|k|e^{-\tau|k|}\frac{\varepsilon}{\sigma^{2}}e^{-\gamma|a-b|}
\preceq \frac{1}{\tau^{\# \mathcal{A}+1}}\frac{\varepsilon}{\sigma^{2}}e^{-\gamma|a-b|}.
\end{aligned}
\end{equation*}
Using Young's inequality (2) in \cite{EK}, we obtain
\begin{equation*}\label{es88}
|||X_{(1-\mathcal{T}^{1}_{\Delta'})f^{2}}|||^{T}_{p,D(\rho-\tau,\sigma^{2},\sigma)\times U'}\preceq\frac{1}{\gamma^{d+p}}\frac{1}{\tau^{\# \mathcal{A}+1}}\frac{\varepsilon}{\sigma^{2}}e^{-\frac{1}{2}\gamma\Delta'},
\end{equation*}
which together with \eqref{es84} leads to
\begin{equation*}\label{es89}
|||X_{(1-\mathcal{T}_{\Delta'})f^{2}}|||^{T}_{p,D(\rho-\tau,\sigma^{2},\sigma)\times U'}
\preceq\frac{1}{\tau^{\# \mathcal{A}}}e^{-\frac{1}{2}\tau\Delta'}\varepsilon+\frac{1}{\gamma^{d+p}}\frac{1}{\tau^{\# \mathcal{A}+1}}\frac{\varepsilon}{\sigma^{2}}e^{-\frac{1}{2}\gamma\Delta'}.
\end{equation*}
In conclusion, we have
\begin{equation}\label{es90}
|||X_{(1-\mathcal{T}_{\Delta'})f^{low}}|||^{T}_{p,D(\rho-\tau,\sigma^{2},\sigma)\times U'}
\preceq\frac{1}{\tau^{\# \mathcal{A}}}e^{-\frac{1}{2}\tau\Delta'}\varepsilon+\frac{1}{\gamma^{d+p}}\frac{1}{\tau^{\# \mathcal{A}+1}}\frac{\varepsilon}{\sigma^{2}}e^{-\frac{1}{2}\gamma\Delta'}.
\end{equation}

\textbf{Estimate of  the vector field generated by $(1-\mathcal{T}_{\Delta'})g^{low}=(1-\mathcal{T}_{\Delta'}) \{f^{high}, s\}^{low}$.}
From  equation (\ref{es4}), we obtain
\begin{equation}\label{es91}
[s^{\varphi}]_{\Lambda,\gamma,\sigma;U',\rho,\mu}\preceq\frac{1}{\kappa^{2}}(\Delta')^{\exp}\varepsilon.
\end{equation}

Applying  \cite[Proposition 6.6]{EK} to  equation (\ref{es5}), we obtain
\begin{equation*}\label{es92}
\|\hat{S}_1(k;\cdot)\|_{p,\gamma;U'}\preceq\frac{1}{\kappa^{2}}(\Delta\Delta')^{\exp}\|\hat{F}_1(k;\cdot)+\hat{G}_1(k;\cdot)\|_{p,\gamma;U'}.
\end{equation*}
Noticing that
\begin{equation*}\label{es93}
\{f^{high},s^{\varphi}\}=-\langle \partial_{r}f^{high},\partial_{\varphi}s^{\varphi}\rangle,
\end{equation*}
it follows from
 (\ref{ho2}), (\ref{es91}) and  \cite[Equation (42)]{EK} that
\begin{equation}\label{es94}
[\{f^{high},s^{\varphi}\}]_{\Lambda,\gamma,\sigma;U',\rho^{(1)},\mu^{(1)}}\preceq\frac{1}{\rho-\rho^{(1)}}\frac{1}{\mu-\mu^{(1)}}\frac{1}{\kappa^{2}}(\Delta')^{\exp}\varepsilon,
\end{equation}
which together with
(\ref{ho2}) and  \eqref{es1} implies
\begin{equation*}\label{es95}
[s^{1}]_{\Lambda,\gamma,\sigma;U',\rho^{(1)},\mu^{(1)}}\preceq\frac{1}{\kappa^{4}}(\Delta\Delta')^{\exp}\frac{1}{\rho-\rho^{(1)}}\frac{1}{\mu-\mu^{(1)}}\varepsilon.
\end{equation*}
Since $s^{1}$ is independent of $r$, there is
\begin{equation}\label{es96}
[s^{1}]_{\Lambda,\gamma,\sigma;U',\rho^{(1)},\mu}\preceq\frac{1}{\kappa^{4}}(\Delta\Delta')^{\exp}\frac{1}{\rho-\rho^{(1)}}\frac{1}{\mu}\varepsilon.
\end{equation}

Next we estimate
\begin{equation*}\label{es97}
\{f^{high},s^{1}\}=-\langle \partial_{r}f^{high},\partial_{\varphi}s^{1}\rangle+\langle \partial_{\zeta}f^{high},J\partial_{\zeta}s^{1}\rangle.
\end{equation*}
By (\ref{ho2}), (\ref{es96}) and the Cauchy estimate (42) in \cite{EK}, we have
\begin{equation*}\label{es98}
[\langle \partial_{r}f^{high},\partial_{\varphi}s^{1}\rangle]_{\Lambda,\gamma,\sigma;U',\rho^{(2)},\mu^{(1)}}
\preceq\frac{1}{\kappa^{4}}(\Delta\Delta')^{\exp}\frac{1}{\rho-\rho^{(1)}}\frac{1}{\rho^{(1)}-\rho^{(2)}}\frac{1}{\mu-\mu^{(1)}}\frac{1}{\mu}\varepsilon.
\end{equation*}
Applying further Proposition 3.1 (ii) in \cite{EK}, we have
\begin{equation*}\label{es99}
[\langle \partial_{\zeta}f^{high},J\partial_{\zeta}s^{1}\rangle]_{\Lambda,\gamma,\sigma^{(1)};U',\rho^{(1)},\mu}
\preceq\frac{1}{\kappa^{4}}(\Delta\Delta')^{\exp}\frac{1}{\rho-\rho^{(1)}}\frac{1}{\sigma-\sigma^{(1)}}\frac{1}{\sigma}\frac{1}{\mu}\varepsilon,
\end{equation*}
which implies
\begin{equation*}\label{es100}
[\{f^{high},s^{1}\}]_{\left\{\substack{\Lambda,\gamma,\sigma^{(1)};\\ U',\rho^{(2)},\mu^{(1)}}\right\}}
\preceq\frac{1}{\kappa^{4}}(\Delta\Delta')^{\exp}\bm{\delta}\varepsilon,
\end{equation*}
where
\begin{equation}\label{bm delta}
  \bm{\delta}=\frac{1}{\rho-\rho^{(1)}}
\left(\frac{1}{\sigma-\sigma^{(1)}}\frac{1}{\sigma}+\frac{1}{\rho^{(1)}-\rho^{(2)}}
\frac{1}{\mu-\mu^{(1)}}\right)\frac{1}{\mu}.
\end{equation}
By (\ref{es1}) and
(\ref{es94}),
we have
\begin{equation}\label{es101}
[g^{low}]_{\left\{\substack{\Lambda,\gamma,\sigma^{(1)};\\ U',\rho^{(2)},\mu^{(1)}}\right\}}
\preceq\frac{1}{\kappa^{4}}(\Delta\Delta')^{\exp}\bm{\delta} \varepsilon,
\end{equation}
which implies
\begin{equation}\label{es102}
\sup_{(\varphi,w)\in D(\rho^{(2)})\times U'}\left\{
\begin{aligned} |G_{2}(\varphi;w)|_{\gamma},\\
|\partial_{w}G_{2}(\varphi;w)|_{\gamma}\end{aligned}\right\}
\preceq\frac{(\Delta\Delta')^{\exp}}{(\sigma^{(1)})^{2}\kappa^{4}}
\bm{\delta'} \varepsilon,
\end{equation}
where
\begin{equation*}
  \bm{\delta'}= \frac{1}{\rho-\rho^{(1)}}
\left(\frac{1}{\sigma-\sigma^{(1)}}\frac{1}{\sigma}+\frac{1}{\rho^{(1)}-\rho^{(2)}}
\frac{1}{\mu}\right)\frac{1}{\mu}.
\end{equation*}
By (\ref{es1})
and (\ref{es70}),
 we have
\begin{equation}\label{es103}
|||X_{g^{low}}|||^{T}_{p,D(\rho-4\tau,(\sigma-4\tau)^{2},\sigma-4\tau)\times U'}
\preceq\frac{d_{\Delta}^{d}\varepsilon}{\tau^{4}\kappa^{4}}.
\end{equation}
Following the proof of (\ref{es90}), we get from (\ref{es102}) and (\ref{es103}) that
\begin{equation}\label{es104}
\begin{aligned}
|||X_{(1-\mathcal{T}_{\Delta'})g^{low}}||&|^{T}_{p,D(\rho-5\tau,(\sigma-4\tau)^{2},\sigma-4\tau)
\times U'}\\
\preceq & ~ \frac{1}{\tau^{\# \mathcal{A}}}e^{-\frac{1}{2}\tau\Delta'}\frac{d_{\Delta}^{d}\varepsilon}{\tau^{4}\kappa^{4}}
+\frac{1}{\gamma^{d+p}}\frac{1}{\tau^{\# \mathcal{A}+1}}e^{-\frac{1}{2}\gamma\Delta'}
\frac{(\Delta\Delta')^{\exp}}{(\sigma^{(1)})^{2}\kappa^{4}}
\bm{\delta'} \varepsilon \\
\preceq & ~\frac{d_{\Delta}^{d}}{\kappa^{4}\tau^{\# \mathcal{A}+4}}e^{-\frac{1}{2}\tau\Delta'}\varepsilon
+\frac{1}{\gamma^{d+p}}\frac{1}{\tau^{\# \mathcal{A}+3}}e^{-\frac{1}{2}\gamma\Delta'}\frac{(\Delta\Delta')^{\exp}}{\sigma^{6}\kappa^{4}}\varepsilon,
\end{aligned}
\end{equation}
where $\rho^{(1)}=\rho-\tau,\rho^{(2)}=\rho-2\tau,\sigma^{(1)}=\sigma-\tau$.

\textbf{Estimate of the vector field $X_{f_{1}^{low}}$.}
Using (\ref{ho1}), (\ref{ho3}), (\ref{es76}) and (\ref{es103}), we have
\begin{equation*}\label{es105}
|||X_{\{h,s\}}|||^{T}_{p,D(\rho-5\tau,(\sigma-5\tau)^{2},\sigma-5\tau)\times U'}\preceq\frac{d_{\Delta}^{d}\varepsilon}{\tau^{4}\kappa^{4}},
\end{equation*}
which together with Proposition \ref{p} and  (\ref{es75}) implies
\begin{equation*}\label{es106}
|||X_{\{\{h,s\},s\}}|||^{T}_{p,D(\rho-6\tau,(\sigma-6\tau)^{2},\sigma-6\tau)\times U'}\preceq\frac{d_{\Delta}^{4d}\varepsilon^{2}}{\tau^{10}\kappa^{10}}.
\end{equation*}
By Proposition \ref{p}, (\ref{ho1}), (\ref{es75}), we have
\begin{equation}\label{es107}
\begin{aligned}
& |||X_{\{f^{low},s\}}|||^{T}_{p,D(\rho-6\tau,(\sigma-6\tau)^{2},\sigma-6\tau)\times U'}\preceq\frac{d_{\Delta}^{3d}\varepsilon^{2}}{\tau^{6}\kappa^{6}},\\
& |||X_{\{f^{high},s\}}|||^{T}_{p,D(\rho-6\tau,(\sigma-6\tau)^{2},\sigma-6\tau)\times U'}\preceq\frac{d_{\Delta}^{3d}\varepsilon}{\tau^{6}\kappa^{6}},
\end{aligned}
\end{equation}
which implies
\begin{equation*}\label{es109}
|||X_{\{\{f^{high},s\},s\}}|||^{T}_{p,D(\rho-7\tau,(\sigma-7\tau)^{2},\sigma-7\tau)\times U'}\preceq\frac{d_{\Delta}^{6d}\varepsilon^{2}}{\tau^{12}\kappa^{12}}.
\end{equation*}
Then applying Theorem 3.3 in \cite{CLY},
we have
\begin{equation}\label{es110}
\begin{aligned}
& |||X_{\int_{0}^{1}(1-t)\{\{h,s\},s\}\circ X^{t}_{s}dt}|||^{T}_{p,D(\rho-7\tau,(\sigma-7\tau)^{2},\sigma-7\tau)\times U'}
\preceq\frac{d_{\Delta}^{4d}\varepsilon^{2}}{\tau^{10}\kappa^{10}},\\
& |||X_{\int_{0}^{1}\{f^{low},s\}\circ X^{t}_{s}dt}|||^{T}_{p,D(\rho-7\tau,(\sigma-7\tau)^{2},\sigma-7\tau)\times U'}\preceq\frac{d_{\Delta}^{3d}\varepsilon^{2}}{\tau^{6}\kappa^{6}},\\
& |||X_{\int_{0}^{1}(1-t)\{\{f^{high},s\},s\}\circ X^{t}_{s}dt}|||^{T}_{p,D(\rho-8\tau,(\sigma-8\tau)^{2},\sigma-8\tau)\times U'}\preceq\frac{d_{\Delta}^{6d}\varepsilon^{2}}{\tau^{12}\kappa^{12}}.
\end{aligned}
\end{equation}
By (\ref{es90}), (\ref{es104}), \eqref{es110} and (\ref{es113}), we have
\begin{equation*}\label{es115}
\begin{aligned}
|||X_{f_{1}^{low}}&|||^{T}_{p,D(\rho-8\tau,(\sigma-8\tau)^{2},\sigma-8\tau)\times U'}\\
&\preceq\frac{d_{\Delta}^{d}}{\kappa^{4}\tau^{\# \mathcal{A}+4}}e^{-\frac{1}{2}\tau\Delta'}\varepsilon
+\frac{1}{\gamma^{d+p}}\frac{1}{\tau^{\# \mathcal{A}+3}}e^{-\frac{1}{2}\gamma\Delta'}\frac{(\Delta\Delta')^{\exp}}{\sigma^{6}\kappa^{4}}\varepsilon+\frac{d_{\Delta}^{6d}\varepsilon^{2}}{\tau^{12}\kappa^{12}}.
\end{aligned}
\end{equation*}

\textbf{Estimate of the vector field $X_{f_{1}^{high}}$.}
By (\ref{ho1}), \eqref{es107},
(\ref{es110})
and \eqref{es113},
we have
\begin{equation*}\label{es116}
|||X_{f_{1}^{high}}|||^{T}_{p,D(\rho-8\tau,(\sigma-8\tau)^{2},\sigma-8\tau)\times U'}\preceq 1+\frac{d_{\Delta}^{3d}\varepsilon}{\tau^{6}\kappa^{6}}+\frac{d_{\Delta}^{6d}\varepsilon^{2}}{\tau^{12}\kappa^{12}}.
\end{equation*}

\subsection{Verification of \eqref{ho16}-\eqref{ho17}.}
In this part, we shall verify the estimates \eqref{ho16}-\eqref{ho17}.
From   (\ref{es6}), (\ref{ho2}) and (\ref{es101}), we obtain
\begin{equation*}\label{es117}
[s^{0}]_{\Lambda,\gamma,\sigma^{(1)};U',\rho^{(2)},\mu^{(1)}}
\preceq\frac{1}{\kappa^{6}}(\Delta\Delta')^{\exp} \bm{\delta}
\varepsilon,
\end{equation*}
where $\bm{\delta}$ is given by \eqref{bm delta}.
Since $s^{0}$ is independent of $\zeta$, we obtain
\begin{equation}\label{es118}
[s^{0}]_{\Lambda,\gamma,\sigma;U',\rho^{(2)},\mu^{(1)}}
\preceq\frac{1}{\kappa^{6}}(\Delta\Delta')^{\exp}\frac{1}{\rho-\rho^{(1)}}\left(\frac{1}{\sigma^{2}}+\frac{1}{\rho^{(1)}-\rho^{(2)}}\frac{1}{\mu-\mu^{(1)}}\right)\frac{1}{\mu}\varepsilon.
\end{equation}

Applying Proposition 6.7 in \cite{EK} to  equation (\ref{es7}), it follows from  (\ref{ho2}) and (\ref{es101}) that
\begin{equation}\label{es119}
\begin{aligned}
& [s^{2}]_{\left\{\substack{\Lambda'+d_\Delta+2,\gamma,\sigma^{(1)};\\
U',\rho^{(2)},\mu^{(1)}}\right\}}
\preceq\frac{1}{\kappa^{7}}(\Delta\Delta')^{\exp} \bm{\delta}
\varepsilon,\\
& [h_{1}]_{\left\{\substack{\Lambda'+d_\Delta+2,\gamma,\sigma^{(1)};\\ U',\rho^{(2)},\mu^{(1)}}\right\}}
\preceq\frac{1}{\kappa^{4}}(\Delta\Delta')^{\exp} \bm{\delta}
\varepsilon.
\end{aligned}
\end{equation}
Using (\ref{es91}), (\ref{es96}), (\ref{es118}) and (\ref{es119}), we obtain
\begin{equation*}\label{es121}
[s]_{\left\{\substack{\Lambda'+d_\Delta+2,\gamma,\sigma^{(1)};\\U',\rho^{(2)},\mu^{(1)}}\right\}}
\preceq\frac{1}{\kappa^{7}}(\Delta\Delta')^{\exp}\frac{1}{\rho-\rho^{(1)}}\left(\frac{1}{\sigma-\sigma^{(1)}}\frac{1}{\sigma}+\frac{1}{\rho^{(1)}-\rho^{(2)}}\frac{1}{\mu-\mu^{(1)}}\right)\frac{1}{\mu}\varepsilon.
\end{equation*}

This completes the proof of Proposition \ref{p1}.

\section{Proof of the KAM theorem \ref{t1}}\label{sect 4}

This section is devoted to the proof of Theorem \ref{t1}.
In subsection \ref{sec 4.1}, we take the normal form computations.
In subsection \ref{sec 4.2}, we establish and prove the KAM iterative lemma,
based on which Theorem \ref{t1} is an immediate result.

\subsection{The normal form computation.}\label{sec 4.1}

For $\rho_{+}<\rho$, $\gamma_{+}<\gamma$, let
$$\Delta'=80(\log\frac{1}{\varepsilon})^{2}\frac{1}{\min(\gamma-\gamma_{+},\rho-\rho_{+})},$$
and  $n=[\log\frac{1}{\varepsilon}]$.
Assume $\rho=\sigma$, $\mu=\sigma^{2}$, $d_{\Delta}\gamma\leq1$. For $1\leq j\leq n$, let
\begin{equation*}
\varepsilon_{j}=\frac{\varepsilon}{\kappa^{20}}\varepsilon_{j-1}, \ \varepsilon_{0}=\varepsilon,
\end{equation*}
\begin{equation*}
\gamma_{j}=\gamma-j\frac{\gamma-\gamma_{+}}{n}, \ \gamma_{0}=\gamma,
\end{equation*}
\begin{equation*}
\rho_{j}=\rho-j\frac{\rho-\rho_{+}}{n}, \ \rho_{0}=\rho,
\end{equation*}
\begin{equation*}
\sigma_{j}=\sigma-j\frac{\sigma-\sigma_{+}}{n}, \ \sigma_{0}=\sigma,
\end{equation*}
\begin{equation*}
\mu_{j}=\sigma_{j}^{2}, \ \mu_{0}=\mu,
\end{equation*}
\begin{equation*}
\Lambda_{j}=\Lambda_{j-1}+d_{\Delta}+30, \ \Lambda_{0}=\mathrm{cte}.\max(\Lambda,d_{\Delta}^{2},d_{\Delta'}^{2}),
\end{equation*}
where the constant $\mathrm{cte}$. is the one in Proposition 6.7 in \cite{EK}.

We have the following lemma.
\begin{lem}\label{l1}
For $0\leq j<n$, consider the Hamiltonian $h+h_1+\cdots+h_{j}+f_{j}$, where
\begin{equation*}
h(r,\zeta;w)=\langle\omega(w),r\rangle+\frac{1}{2}\langle\zeta,(\Omega(w)+H(w))\zeta\rangle
\end{equation*}
satisfy \eqref{as1}-\eqref{as9}, $H(w),\partial_{w}H(w)$ are T\"{o}plitz at $\infty$ and $\mathcal{NF}_{\Delta}$ for all $w\in U$.
Let $U' \subset U$ satisfy \eqref{sd1}-\eqref{sd4}. For all $w\in U'$,
\begin{equation*}
h_{j}=a_{j}(w)+\langle \chi_{j}(w),r\rangle+\frac{1}{2}\langle\zeta,H_{j}(w)\zeta\rangle,
\end{equation*}
\begin{equation*}
f_{j}=f_{j}^{low}+f_{j}^{high}
\end{equation*}
satisfy
\begin{equation}\label{fi1}
|||X_{f_{j}^{low}}|||^{T}_{p,D(\rho_{j},\mu_{j},\sigma_{j})\times U'}\leq\beta^{j}\varepsilon_{j}, \  |||X_{f_{j}^{high}}|||^{T}_{p,D(\rho_{j},\mu_{j},\sigma_{j})\times U'}\leq1,
\end{equation}
\begin{equation}\label{fi2}
[f_{j}^{low}]_{\Lambda_{j},\gamma_{j},\sigma_{j};U',\rho_{j},\mu_{j}}\leq\beta^{j}\varepsilon_{j}, \ [f_{j}^{high}]_{\Lambda_{j},\gamma_{j},\sigma_{j};U',\rho_{j},\mu_{j}}\leq 1
\end{equation}
for some
\begin{equation*}
\beta\preceq \max\left(\frac{1}{\gamma-\gamma_{+}},\frac{1}{\rho-\rho_{+}},\Delta,\Lambda,\log\frac{1}{\varepsilon}\right)^{\exp_1}.
\end{equation*}

Then there exists an exponent $\exp_2$ such that if
\begin{equation*}
\varepsilon\preceq\kappa^{20}\min\left(\gamma-\gamma_{+},\rho-\rho_{+},\frac{1}{\Delta},\frac{1}{\Lambda},\frac{1}{\log\frac{1}{\varepsilon}}\right)^{\exp_2},
\end{equation*}
then for all $w\in U'$, there is a real analytic symplectic map $\Phi_{j}$ such that
\begin{equation*}
(h+h_1+\cdots+h_{j}+f_{j})\circ\Phi_{j}=h+h_1+\cdots+h_{j+1}+f_{j+1},
\end{equation*}
with the estimates
\begin{equation}\label{fi3}
|||X_{f^{low}_{j+1}}|||^{T}_{p,D(\rho_{j+1},\mu_{j+1},\sigma_{j+1})\times U'}
\preceq\beta^{j+1}\varepsilon_{j+1},
\end{equation}
\begin{equation}\label{fi4}
|||X_{f^{high}_{j+1}}|||^{T}_{p,D(\rho_{j+1},\mu_{j+1},\sigma_{j+1})\times U'}
\preceq 1+\frac{1}{\kappa^{6}}\beta^{j+1}\varepsilon_{j}+\beta^{j+1}\varepsilon_{j+1},
\end{equation}
\begin{equation}\label{fi5}
[f^{low}_{j+1}]_{\Lambda_{j+1},\gamma_{j+1},\sigma_{j+1};U',\rho_{j+1},\mu_{j+1}}
\preceq\beta^{j+1}\varepsilon_{j+1},
\end{equation}
\begin{equation}\label{fi6}
[f^{high}_{j+1}]_{\Lambda_{j+1},\gamma_{j+1},\sigma_{j+1};U',\rho_{j+1},\mu_{j+1}}
\preceq 1+\frac{1}{\kappa^{7}}\beta^{j+1}\varepsilon_{j}+\beta^{j+1}\varepsilon_{j+1},
\end{equation}
where the exponents $\exp_1$, $\exp_2$ depend on $d, \# \mathcal{A}, p$.

\end{lem}

\noindent\textbf{Proof.}~
By Proposition \ref{p1}, we can solve the homological equation
\begin{equation}\label{a1}
\{h,s_{j}\}=-\mathcal{T}_{\Delta'}f_{j}^{low}-\mathcal{T}_{\Delta'}\{f_{j}^{high},s_{j}\}^{low}+h_{j+1}
\end{equation}
with the estimates
\begin{equation}\label{a2}
\begin{aligned}
& [s_{j}]_{\left\{\substack{\Lambda_{j}+d_\Delta+2,\gamma_{j},\sigma_{j}^{(1)};\\ U',\rho_{j}^{(1)},\mu_{j}^{(1)}}\right\}}
\preceq\frac{1}{\kappa^{7}}(\Delta\Delta')^{\exp}\bm{\delta_{1}}\varepsilon_{j},\\
& [h_{j+1}]_{\left\{\substack{\Lambda_{j}+d_\Delta+2,\gamma_{j},\sigma_{j}^{(1)};\\U',\rho_{j}^{(1)},\mu_{j}^{(1)}}\right\}}
\preceq\frac{1}{\kappa^{4}}(\Delta\Delta')^{\exp}\bm{\delta_{1}}\varepsilon_{j},
\end{aligned}
\end{equation}
where
\begin{equation*}
\bm{\delta_{1}}= \frac{1}{\rho_{j}-\rho_{j}^{(1)}}
\left(\frac{1}{\sigma_{j}-\sigma_{j}^{(1)}}\frac{1}{\sigma_{j}}+\frac{1}{\rho_{j}-\rho_{j}^{(1)}}\frac{1}{\mu_{j}-\mu_{j}^{(1)}}\right)\frac{1}{\mu_{j}} \beta^{j}.
\end{equation*}

\bigskip

\textbf{Step 1: computation of $f_{j+1}$.} In this step, we compute  the new perturbation $f_{j+1}= f_{j+1}^{low}+f_{j+1}^{high}$.
Using Taylor's formula, by the homological equation (\ref{a1}), we obtain
\begin{equation}\label{a4}
\begin{aligned}
(h+h_1+\cdots&+h_{j}+f_{j})\circ X^{t}_{s_{j}}\mid_{t=1}
=h+h_{1}+\cdots+h_{j+1}\\
&+(1-\mathcal{T}_{\Delta'})f_{j}^{low}+f_{j}^{high}+(1-\mathcal{T}_{\Delta'})\{f_{j}^{high},s_{j}\}^{low}+\{f_{j}^{high},s_{j}\}^{high} \\
&+\int_{0}^{1}(1-t)\{\{h,s_{j}\},s_{j}\}\circ X^{t}_{s_{j}}dt+\int_{0}^{1}\{h_1+\cdots+h_{j},s_{j}\}\circ X^{t}_{s_{j}}dt \\
&+\int_{0}^{1}\{f_{j}^{low},s_{j}\}\circ X^{t}_{s_{j}}dt+\int_{0}^{1}(1-t)\{\{f_{j}^{high},s_{j}\},s_{j}\}\circ X^{t}_{s_{j}}dt.
\end{aligned}
\end{equation}
Let $\Phi_{j}= X^{t}_{s_{j}}\mid_{t=1}$ and
\begin{equation*}
(h+h_1+\cdots+h_{j}+f_{j})\circ X^{t}_{s_{j}}\mid_{t=1}
=h+h_{1}+\cdots+h_{j+1}+f_{j+1}.
\end{equation*}
Then we get from
(\ref{a1}) and  (\ref{a4}) that
\begin{equation*}\label{a20}
\begin{aligned}
f_{j+1}=& (1-\mathcal{T}_{\Delta'})f_{j}^{low}+f_{j}^{high}+(1-\mathcal{T}_{\Delta'})\{f_{j}^{high},s_{j}\}^{low}+\{f_{j}^{high},s_{j}\}^{high}\\
&+\int_{0}^{1}(1-t)\{-\mathcal{T}_{\Delta'}f_{j}^{low},s_{j}\}\circ X^{t}_{s_{j}}dt
+\int_{0}^{1}(1-t)\{-\mathcal{T}_{\Delta'}\{f_{j}^{high},s_{j}\}^{low},s_{j}\}\circ X^{t}_{s_{j}}dt\\
&+\int_{0}^{1}\{h_1+\cdots+h_{j}+(1-t)h_{j+1},s_{j}\}\circ X^{t}_{s_{j}}dt
+\int_{0}^{1}\{f_{j}^{low},s_{j}\}\circ X^{t}_{s_{j}}dt\\
&+\int_{0}^{1}(1-t)\{\{f_{j}^{high},s_{j}\},s_{j}\}\circ X^{t}_{s_{j}}dt.
\end{aligned}
\end{equation*}
As a result, there is
\begin{align*}
f^{low}_{j+1}=&(1-\mathcal{T}_{\Delta'})f_{j}^{low}+(1-\mathcal{T}_{\Delta'})
\{f_{j}^{high},s_{j}\}^{low}+\left(\int_{0}^{1}\{f_{j}^{low},s_{j}\}\circ X^{t}_{s_{j}}dt\right)^{low}\\
&+\left(\int_{0}^{1}(1-t)\{-\mathcal{T}_{\Delta'}f_{j}^{low},s_{j}\}\circ X^{t}_{s_{j}}dt\right)^{low}\\
&+\left(\int_{0}^{1}(1-t)\{-\mathcal{T}_{\Delta'}\{f_{j}^{high},s_{j}\}^{low},s_{j}\}\circ X^{t}_{s_{j}}dt\right)^{low}
\\
&+\left(\int_{0}^{1}(1-t)\{\{f_{j}^{high},s_{j}\},s_{j}\}\circ X^{t}_{s_{j}}dt\right)^{low}\\
&+\left(\int_{0}^{1}\{h_1+\cdots+h_{j}+(1-t)h_{j+1},s_{j}\}\circ X^{t}_{s_{j}}dt\right)^{low},
\end{align*}
and
\begin{align*}
 f^{high}_{j+1}=&f_{j}^{high}+\{f_{j}^{high},s_{j}\}^{high}
+\left(\int_{0}^{1}(1-t)\{-\mathcal{T}_{\Delta'}f_{j}^{low},s_{j}\}\circ X^{t}_{s_{j}}dt\right)^{high}\\
&+\left(\int_{0}^{1}(1-t)\{-\mathcal{T}_{\Delta'}\{f_{j}^{high},s_{j}\}^{low},s_{j}\}\circ X^{t}_{s_{j}}dt\right)^{high}\\
&+\left(\int_{0}^{1}(1-t)\{\{f_{j}^{high},s_{j}\},s_{j}\}\circ X^{t}_{s_{j}}dt\right)^{high}\\
&+\left(\int_{0}^{1}\{f_{j}^{low},s_{j}\}\circ X^{t}_{s_{j}}dt\right)^{high}\\
&+\left(\int_{0}^{1}\{h_1+\cdots+h_{j}+(1-t)h_{j+1},s_{j}\}\circ X^{t}_{s_{j}}dt\right)^{high}.
\end{align*}

\bigskip

\textbf{Step 2: estimates of $f_{j+1}^{low}$ and $f_{j+1}^{high}$.} In this part, we verify the estimates
\eqref{fi5} and \eqref{fi6}. The various estimates of the poission brackets are based on \eqref{a2}
and  \cite[ Proposition 3.3]{EK}.

To begin with, we see from (\ref{fi2}) that
\begin{equation}\label{a5}
[(1-\mathcal{T}_{\Delta'})f_{j}^{low}]_{\left\{\substack{\Lambda_{j},\gamma_{j}^{(1)},\sigma_{j};\\U',\rho_{j}^{(1)},\mu_{j}}\right\}}\preceq
\left[\left(\frac{1}{\rho_{j}-\rho_{j}^{(1)}}\right)^{\#\mathcal{A}}e^{-\frac{1}{2}(\rho_{j}-\rho_{j}^{(1)})\Delta'}+e^{-(\gamma_{j}-\gamma_{j}^{(1)})\Delta'}\right]\beta^{j}\varepsilon_{j}.
\end{equation}
By (\ref{es101}), there is
\begin{equation}\label{a6}
[\{f_{j}^{high},s_{j}\}^{low}]_{\left\{\substack{\Lambda_{j},\gamma_{j},\sigma_{j}^{(1)};\\U',\rho_{j}^{(1)},\mu_{j}^{(1)}}\right\}}
\preceq\frac{1}{\kappa^{4}}(\Delta\Delta')^{\exp}\bm{\delta_{1}}\varepsilon_{j},
\end{equation}
and thus
\begin{equation}\label{a7}
\begin{aligned}
&\quad [(1-\mathcal{T}_{\Delta'}) \{f_{j}^{high},s_{j}\}^{low}]_{\left\{\substack{\Lambda_{j},\gamma^{(1)}_{j},\sigma_{j}^{(1)};\\U',\rho_{j}^{(2)},\mu_{j}^{(1)}}\right\}}\\
\preceq & \left[\left(\frac{1}{\rho_{j}^{(1)}-\rho_{j}^{(2)}}\right)^{\#\mathcal{A}}e^{-\frac{1}{2}(\rho_{j}^{(1)}-\rho_{j}^{(2)})\Delta'}+e^{-(\gamma_{j}-\gamma_{j}^{(1)})\Delta'}\right]
\frac{1}{\kappa^{4}}(\Delta\Delta')^{\exp}\bm{\delta_{1}} \varepsilon_{j}.
\end{aligned}
\end{equation}

Let
\begin{equation*}
\bm{\delta_{2}}= (\Lambda_{j}+d_\Delta+2)^{2}\left(\frac{1}{\gamma_{j}-\gamma_{j}^{(1)}}\right)^{d+1}\frac{1}{\sigma_{j}^{(1)}-\sigma_{j}^{(2)}}\frac{1}{\sigma_{j}^{(1)}}
+\frac{1}{\rho_{j}^{(1)}-\rho_{j}^{(2)}}\frac{1}{\mu_{j}^{(1)}-\mu_{j}^{(2)}}.
\end{equation*}
By (\ref{fi2}), (\ref{a2}) and   \cite[ Proposition 3.3]{EK}, we have
\begin{equation}\label{a8}
[\{f_{j}^{high},s_{j}\}]_{\left\{\substack{\Lambda_{j}+d_\Delta+5,\gamma_{j}^{(1)},\sigma_{j}^{(2)};\\
U',\rho_{j}^{(2)},\mu_{j}^{(2)}}\right\}}\preceq
\bm{\delta_{2}} \frac{1}{\kappa^{7}}(\Delta\Delta')^{\exp} \bm{\delta_{1}} \varepsilon_{j},
\end{equation}
and
\begin{equation*}\label{a9}
[\{f_{j}^{low},s_{j}\}]_{\left\{\substack{\Lambda_{j}+d_\Delta+5,\gamma_{j}^{(1)},\sigma_{j}^{(2)};\\U',\rho_{j}^{(2)},\mu_{j}^{(2)}}\right\}}\preceq \bm{\delta_{2}}  \frac{1}{\kappa^{7}}(\Delta\Delta')^{\exp} \bm{\delta_{1}} \beta^{j}  \varepsilon_{j}^{2}.
\end{equation*}
Moreover, we obtain from
 (\ref{a2}), (\ref{a6}) and  \cite[Proposition 3.3]{EK} that
\begin{equation}\label{a10}
[\{\{f_{j}^{high},s_{j}\}^{low},s_{j}\}]_{\left\{\substack{\Lambda_{j}+d_\Delta+5,\gamma_{j}^{(1)},\sigma_{j}^{(2)}; \\U',\rho_{j}^{(2)},\mu_{j}^{(2)}
}\right\}}\preceq \bm{\delta_{2}} \frac{1}{\kappa^{11}}(\Delta\Delta')^{\exp} \bm{\delta_{1}}^{2} \varepsilon_{j}^{2}.
\end{equation}
Let
\begin{equation*}
\bm{\delta_{3}}=  (\Lambda_{j}+d_\Delta+5)^{2}\left(\frac{1}{\gamma_{j}^{(1)}-\gamma_{j}^{(2)}}\right)^{d+1}\frac{1}{\sigma_{j}^{(2)}-\sigma_{j}^{(3)}}\frac{1}{\sigma_{j}^{(2)}}
+\frac{1}{\rho_{j}^{(2)}-\rho_{j}^{(3)}}\frac{1}{\mu_{j}^{(2)}-\mu_{j}^{(3)}}.
\end{equation*}
By (\ref{a2}),  (\ref{a8}), using Proposition 3.3 in \cite{EK}, we have
\begin{equation*}\label{a11}
[\{\{f_{j}^{high},s_{j}\},s_{j}\}]_{\left\{\substack{\Lambda_{j}+d_\Delta+8,\gamma_{j}^{(2)},\sigma_{j}^{(3)};\\
U',\rho_{j}^{(3)},\mu_{j}^{(3)}}\right\}}
\preceq \bm{\delta_{3} \delta_{2}} \frac{1}{\kappa^{14}} \bm{\delta_{1}}^{2} \varepsilon_{j}^{2},
\end{equation*}
and
\begin{equation*}
\begin{aligned}
& [ \{h_{i+1},s_{j}\}]_{\left\{\substack{\Lambda_{j}+d_\Delta+5,\gamma_{j}^{(1)},\sigma_{j}^{(2)};\\U',\rho_{j}^{(2)},\mu_{j}^{(2)}}\right\}}\\
\preceq&\left[(\Lambda_{j}+d_\Delta+2)^{2}\left(\frac{1}{\gamma_{j}-\gamma_{j}^{(1)}}\right)^{d+p}\frac{\sigma_{j}^{(1)}}{\sigma_{j}^{(1)}-\sigma_{j}^{(2)}}\frac{1}{(\sigma_{i}^{(1)})^{2}}
+\frac{1}{\rho_{j}^{(1)}-\rho_{j}^{(2)}}\frac{1}{\mu_{i}^{(1)}}\right]\\
&\times\frac{1}{\kappa^{11}}(\Delta\Delta')^{\exp}\bm{\delta_{1}}
\cdot \frac{1}{\rho_{i}-\rho_{i}^{(1)}}\left(\frac{1}{\sigma_{i}-\sigma_{i}^{(1)}}\frac{1}{\sigma_{i}}+\frac{1}{\rho_{i}-\rho_{i}^{(1)}}\frac{1}{\mu_{i}-\mu_{i}^{(1)}}\right)\frac{1}{\mu_{i}}
\beta^{i}\varepsilon_{i}\varepsilon_{j}.
\end{aligned}
\end{equation*}

Take
\begin{equation*}
\begin{aligned}
&\rho_{j}^{(l)}=\rho_{j}-\frac{l}{4}(\rho_{j}-\rho_{j+1}), \quad \gamma_{j}^{(l)}=\gamma_{j}-\frac{l}{4}(\gamma_{j}-\gamma_{j+1}), \\
&\sigma_{j}^{(l)}=\sigma_{j}-\frac{l}{4}(\sigma_{j}-\sigma_{j+1}),\quad \mu_{j}^{(l)}=(\sigma_{j}^{(l)})^{2},
\end{aligned}
\end{equation*}
where  $l=1,2,3,4$.
By (\ref{a5}), we have
\begin{equation}\label{a13}
[(1-\mathcal{T}_{\Delta'})f_{j}^{low}]_{\Lambda_{j},\gamma_{j}^{(1)},\sigma_{j};U',\rho_{j}^{(1)},\mu_{j}}\preceq
\beta\varepsilon\beta^{j}\varepsilon_{j}\preceq\beta^{j+1}\varepsilon_{j+1}.
\end{equation}
By (\ref{a7}),
we have
\begin{equation*}\label{a14}
\begin{aligned}
& [(1-\mathcal{T}_{\Delta'})\{f_{j}^{high},s_{j}\}^{low}]_{\Lambda_{j},\gamma^{(1)}_{j},\sigma_{j}^{(1)};U',\rho_{j}^{(2)},\mu_{j}^{(1)}}
\preceq\frac{1}{\kappa^{4}}\beta\varepsilon\beta^{j}\varepsilon_{j}\preceq\beta^{j+1}\varepsilon_{j+1},\\
& [\{f_{j}^{high},s_{j}\}]_{\Lambda_{j}+d_\Delta+5,\gamma_{j}^{(1)},\sigma_{j}^{(2)};U',\rho_{j}^{(2)},\mu_{j}^{(2)}}
\preceq\frac{1}{\kappa^{7}}\beta\beta^{j}\varepsilon_{j}\preceq\frac{1}{\kappa^{7}}\beta^{j+1}\varepsilon_{j},\\
& [\{f_{j}^{low},s_{j}\}]_{\Lambda_{j}+d_\Delta+5,\gamma_{j}^{(1)},\sigma_{j}^{(2)};U',\rho_{j}^{(2)},\mu_{j}^{(2)}}
\preceq\frac{1}{(\Lambda_{j}+d_\Delta+5)^{14}}\beta^{j+1}\varepsilon_{j+1},\\
& [\{\{f_{j}^{high},s_{j}\}^{low},s_{j}\}]_{\Lambda_{j}+d_\Delta+5,\gamma_{j}^{(1)},\sigma_{j}^{(2)};U',\rho_{j}^{(2)},\mu_{j}^{(2)}}
\preceq\frac{1}{(\Lambda_{j}+d_\Delta+5)^{14}}\beta^{j+1}\varepsilon_{j+1},\\
& [\{\{f_{j}^{high},s_{j}\},s_{j}\}]_{\Lambda_{j}+d_\Delta+8,\gamma_{j}^{(2)},\sigma_{j}^{(3)};U',\rho_{j}^{(3)},\mu_{j}^{(3)}}
\preceq\frac{1}{(\Lambda_{j}+d_\Delta+8)^{14}}\beta^{j+1}\varepsilon_{j+1},\\
& [\{h_1+\cdots+h_{j}+(1-t)h_{j+1},s_{j}\}]_{\Lambda_{j}+d_\Delta+5,\gamma_{j}^{(1)},\sigma_{j}^{(2)};U',\rho_{j}^{(2)},\mu_{j}^{(2)}}
\preceq\frac{1}{(\Lambda_{j}+d_\Delta+5)^{14}}\beta^{j+1}\varepsilon_{j+1},
\end{aligned}
\end{equation*}
which imply the following properties
\begin{equation}\label{a23}
\begin{aligned}
&\left[\int_{0}^{1}(1-t)\{-\mathcal{T}_{\Delta'}f_{j}^{low},s_{j}\}\circ X^{t}_{s_{j}}dt\right]_{\Lambda_{j+1},\gamma_{j+1},\sigma_{j+1};U',\rho_{j+1},\mu_{j+1}}
\preceq\beta^{j+1}\varepsilon_{j+1},\\
&\left[\int_{0}^{1}\{f_{j}^{low},s_{j}\}\circ X^{t}_{s_{j}}dt\right]_{\Lambda_{j+1},\gamma_{j+1},\sigma_{j+1};U',\rho_{j+1},\mu_{j+1}}
\preceq\beta^{j+1}\varepsilon_{j+1},\\
&\left[\int_{0}^{1}(1-t)\{-\mathcal{T}_{\Delta'}\{f_{j}^{high},s_{j}\}^{low},s_{j}\}\circ X^{t}_{s_{j}}dt\right]_{\Lambda_{j+1},\gamma_{j+1},\sigma_{j+1};U',\rho_{j+1},\mu_{j+1}}
\preceq\beta^{j+1}\varepsilon_{j+1},\\
&\left[\int_{0}^{1}\{h_1+\cdots+h_{j}+(1-t)h_{j+1},s_{j}\}\circ X^{t}_{s_{j}}dt\right]_{\Lambda_{j+1},\gamma_{j+1},\sigma_{j+1};U',\rho_{j+1},\mu_{j+1}}
\preceq\beta^{j+1}\varepsilon_{j+1},\\
&\left[\int_{0}^{1}(1-t)\{\{f_{j}^{high},s_{j}\},s_{j}\}\circ X^{t}_{s_{j}}dt\right]_{\Lambda_{j+1},\gamma_{j+1},\sigma_{j+1};U',\rho_{j+1},\mu_{j+1}}
\preceq\beta^{j+1}\varepsilon_{j+1}.
\end{aligned}
\end{equation}

Then from the definition of $f_{j+1}^{low}$ and $f_{j+1}^{high}$, we get the desired estimates
\begin{equation}\label{a28}
[f^{low}_{j+1}]_{\Lambda_{j+1},\gamma_{j+1},\sigma_{j+1};U',\rho_{j+1},\mu_{j+1}}
\preceq\beta^{j+1}\varepsilon_{j+1},
\end{equation}
\begin{equation*}\label{a29}
[f^{high}_{j+1}]_{\Lambda_{j+1},\gamma_{j+1},\sigma_{j+1};U',\rho_{j+1},\mu_{j+1}}
\preceq 1+\frac{1}{\kappa^{7}}\beta^{j+1}\varepsilon_{j}+\beta^{j+1}\varepsilon_{j+1}.
\end{equation*}

\bigskip

\textbf{Step 3: estimates of the vector fields $X_{f_{j+1}^{low}}$ and $X_{f_{j+1}^{high}}$.}
By (\ref{es75}) and (\ref{es76}), we have
\begin{equation}\label{a30}
\begin{aligned}
& |||X_{s_{j}}|||^{T}_{p,D(\rho_{j}^{(1)},\mu_{j}^{(1)},\sigma_{j}^{(1)})\times U'}\preceq\frac{1}{\kappa^{6}}\Delta^{\exp}\left(\frac{1}{\rho_{j}-\rho_{j}^{(1)}}\right)^{5}\beta^{j}\varepsilon_{j},\\
& |||X_{h_{j+1}}|||^{T}_{p,D(\rho_{j}^{(1)},\mu_{j}^{(1)},\sigma_{j}^{(1)})\times U'}\preceq\frac{1}{\kappa^{4}}\Delta^{\exp}\left(\frac{1}{\rho_{j}-\rho_{j}^{(1)}}\right)^{4}\beta^{j}\varepsilon_{j}.
\end{aligned}
\end{equation}
Using (\ref{es90}), we have
\begin{equation}\label{a32}
\begin{aligned}
|||X_{(1-\mathcal{T}_{\Delta'})f_{j}^{low}}&|||^{T}_{p,D(\rho_{j}^{(1)},\mu_{j}^{(1)},\sigma_{j}^{(1)})\times U'}
\preceq \left(\frac{1}{\rho_{j}-\rho_{j}^{(1)}}\right)^{\# \mathcal{A}}e^{-\frac{1}{2}(\rho_{j}-\rho_{j}^{(1)})\Delta'}
\beta^{j}\varepsilon_{j}\\
&\qquad+\frac{1}{\gamma_{j}^{d+p}}\left(\frac{1}{\rho_{j}-\rho_{j}^{(1)}}\right)^{\# \mathcal{A}+1}\frac{1}{\sigma_{j}^{2}}e^{-\frac{1}{2}\gamma_{j}\Delta'}\beta^{j}\varepsilon_{j}\preceq\beta^{j+1}\varepsilon_{j+1}.
\end{aligned}
\end{equation}
By (\ref{es103}) and (\ref{es104}), we have
\begin{equation}\label{a33}
|||X_{\{f_{j}^{high},s_{j}\}^{low}}|||^{T}_{p,D(\rho_{j}^{(1)},\mu_{j}^{(1)},\sigma_{j}^{(1)})\times U'}
\preceq\frac{1}{\kappa^{4}}\Delta^{\exp}\left(\frac{1}{\rho_{j}-\rho_{j}^{(1)}}\right)^{4}\beta^{j}\varepsilon_{j},
\end{equation}
\begin{equation}\label{a34}
\textrm{and}\qquad |||X_{(1-\mathcal{T}_{\Delta'})\{f_{j}^{high},s_{j}\}^{low}}|||^{T}_{p,D(\rho_{j}^{(1)},\mu_{j}^{(1)},\sigma_{j}^{(1)})\times U'}
\end{equation}
\begin{equation*}
\preceq\frac{1}{\kappa^{4}}(\Delta\Delta')^{\exp}\left(\frac{1}{\rho_{j}-\rho_{j}^{(1)}}\right)^{\# \mathcal{A}+4}\left[e^{-\frac{1}{2}(\rho_{j}-\rho_{j}^{(1)})\Delta'}
+\frac{1}{\gamma_{j}^{d+p}}\frac{1}{\sigma_{j}^{6}}e^{-\frac{1}{2}\gamma_{j}\Delta'}\right]\beta^{j}\varepsilon_{j}\preceq\beta^{j+1}\varepsilon_{j+1}.
\end{equation*}

By (\ref{fi1}), (\ref{a30}) and Proposition \ref{p}, we have
\begin{equation}\label{a35}
|||X_{\{f_{j}^{low},s_{j}\}}|||^{T}_{p,D(\rho_{j}^{(2)},\mu_{j}^{(2)},\sigma_{j}^{(2)})\times U'}
\preceq\frac{1}{\rho_{j}^{(1)}-\rho_{j}^{(2)}}\frac{1}{\kappa^{6}}\Delta^{\exp}\left(\frac{1}{\rho_{j}-\rho_{j}^{(1)}}\right)^{5}\beta^{2j}\varepsilon_{j}^{2},
\end{equation}
\begin{equation}\label{a36}
|||X_{\{f_{j}^{high},s_{j}\}}|||^{T}_{p,D(\rho_{j}^{(2)},\mu_{j}^{(2)},\sigma_{j}^{(2)})\times U'}
\preceq\frac{1}{\rho_{j}^{(1)}-\rho_{j}^{(2)}}\frac{1}{\kappa^{6}}\Delta^{\exp}\left(\frac{1}{\rho_{j}-\rho_{j}^{(1)}}\right)^{5}\beta^{j}\varepsilon_{j}.
\end{equation}

Using (\ref{es105}), we have
\begin{equation}\label{a37}
|||X_{\{h,s_{j}\}}|||^{T}_{p,D(\rho_{j}^{(1)},\mu_{j}^{(1)},\sigma_{j}^{(1)})\times U'}
\preceq\frac{1}{\kappa^{4}}\Delta^{\exp}\left(\frac{1}{\rho_{j}-\rho_{j}^{(1)}}\right)^{4}\beta^{j}\varepsilon_{j}.
\end{equation}
By (\ref{a30}),
(\ref{a36}), (\ref{a37}), using Proposition \ref{p}, we have
\begin{equation}\label{a38}
\begin{aligned}
& |||X_{\{\{h,s_{j}\},s_j\}}|||^{T}_{p,D(\rho_{j}^{(2)},\mu_{j}^{(2)},\sigma_{j}^{(2)})\times U'}
\preceq\frac{\Delta^{\exp}
\beta^{2j}\varepsilon_{j}^{2}}{(\rho_{j}^{(1)}-\rho_{j}^{(2)})\cdot \kappa^{10}}
\left(\frac{1}{\rho_{j}-\rho_{j}^{(1)}}\right)^{9},\\
& |||X_{\{\{f_{j}^{high},s_{j}\},s_j\}}|||^{T}_{p,D(\rho_{j}^{(3)},\mu_{j}^{(3)},\sigma_{j}^{(3)})\times U'}
\preceq
\frac{\Delta^{\exp} \beta^{2j} \varepsilon_{j}^{2}}{(\rho_{j}^{(2)}-\rho_{j}^{(3)})(\rho_{j}^{(1)}-\rho_{j}^{(2)})
\kappa^{12}}
\left(\frac{1}{\rho_{j}-\rho_{j}^{(1)}}\right)^{10},\\
&|||X_{\{h_{i+1},s_j\}}|||^{T}_{p,D(\rho_{j}^{(2)},\mu_{j}^{(2)},\sigma_{j}^{(2)})\times U'}
\preceq\frac{\Delta^{\exp} \beta^{i+j}\varepsilon_{i}\varepsilon_{j}}{(\rho_{j}^{(1)}-\rho_{j}^{(2)})\kappa^{10}}
\left(\frac{1}{\rho_{j}-\rho_{j}^{(1)}}\right)^{5}
\left(\frac{1}{\rho_{i}-\rho_{i}^{(1)}}\right)^{4}.
\end{aligned}
\end{equation}

By (\ref{a35}) and \eqref{a38}, we see from
 \cite[Theorem 3.3]{CLY} that
 \begin{equation}\label{a44}
 \begin{aligned}
& |||X_{\int_{0}^{1}(1-t)\{\{h,s_{j}\},s_{j}\}\circ X^{t}_{s_{j}}dt}|||^{T}_{p,D(\rho_{j+1},\mu_{j+1},\sigma_{j+1})\times U'}
\preceq\beta^{j+1}\varepsilon_{j+1},\\
& |||X_{\int_{0}^{1}\{h_1+\cdots+h_{j},s_{j}\}\circ X^{t}_{s_{j}}dt}|||^{T}_{p,D(\rho_{j+1},\mu_{j+1},\sigma_{j+1})\times U'}
\preceq\beta^{j+1}\varepsilon_{j+1},\\
& |||X_{\int_{0}^{1}\{f_{j}^{low},s_{j}\}\circ X^{t}_{s_{j}}dt}|||^{T}_{p,D(\rho_{j+1},\mu_{j+1},\sigma_{j+1})\times U'}
\preceq\beta^{j+1}\varepsilon_{j+1},\\
& |||X_{\int_{0}^{1}(1-t)\{\{f_{j}^{high},s_{j}\},s_{j}\}\circ X^{t}_{s_{j}}dt}|||^{T}_{p,D(\rho_{j+1},\mu_{j+1},\sigma_{j+1})\times U'}
\preceq\beta^{j+1}\varepsilon_{j+1}.
\end{aligned}
\end{equation}

Finally, we obtain from the definitions of $f^{low}_{j+1}$ and $f^{high}_{j+1}$ that
\begin{equation*}\label{a48}
|||X_{f^{low}_{j+1}}|||^{T}_{p,D(\rho_{j+1},\mu_{j+1},\sigma_{j+1})\times U'}
\preceq\beta^{j+1}\varepsilon_{j+1},
\end{equation*}
\begin{equation*}\label{a49}
|||X_{f^{high}_{j+1}}|||^{T}_{p,D(\rho_{j+1},\mu_{j+1},\sigma_{j+1})\times U'}
\preceq 1+\frac{1}{\kappa^{6}}\beta^{j+1}\varepsilon_{j}+\beta^{j+1}\varepsilon_{j+1}.
\end{equation*}
This completes the proof of Lemma \ref{l1}. \qed

\subsection{The KAM iterative lemma}\label{sec 4.2}

Assume $\rho=\sigma$, $\mu=\sigma^{2}$, $d_{\Delta}\gamma\leq1$. For $m\geq 0$, let
\begin{equation*}
\varepsilon_{m}=e^{-\frac{1}{20}(\log\frac{1}{\varepsilon_{m-1}})^{2}}, \ \varepsilon_{0}=\varepsilon,
\end{equation*}
\begin{equation*}
\vartheta_{m}=\frac{\sum_{j=1}^{m}\frac{1}{j^{2}}}{2\sum_{j=1}^{\infty}\frac{1}{j^{2}}}, \ \vartheta_{0}=0,
\end{equation*}
\begin{equation*}
\rho_{m}=(1-\vartheta_{m})\rho, \ \rho_{0}=\rho,
\end{equation*}
\begin{equation*}
\sigma_{m}=(1-\vartheta_{m})\sigma, \ \sigma_{0}=\sigma,
\end{equation*}
\begin{equation*}
\mu_{m}=\sigma_{m}^{2}, \ \mu_{0}=\mu,
\end{equation*}
\begin{equation*}
\gamma_{m}=d^{-1}_{\Delta_{m}}, \ \gamma_{0}=\min(\gamma,d^{-1}_{\Delta}),
\end{equation*}
\begin{equation*}
\Delta_{m}=80(\log\frac{1}{\varepsilon_{m-1}})^{2}\frac{1}{\min(\gamma_{m-1},\rho_{m-1}-\rho_{m})}, \ \Delta_{0}=\Delta,
\end{equation*}
\begin{equation*}
\Lambda_{m}=\mathrm{cte}.d^{2}_{\Delta_{m}},
\end{equation*}
where the constant $\mathrm{cte}$. is the one in Proposition 6.7 in \cite{EK}.

We have the following KAM iterative lemma.

\begin{lem}\label{l2}
For $m\geq0$, consider the Hamiltonian $h_{m}+f_{m}$, where
\begin{equation*}
h_{m}=\langle \omega_{m}(w),r\rangle+\frac{1}{2}\langle\zeta,(\Omega(w)+H_{m}(w))\zeta\rangle,
\end{equation*}
$H_{m}(w),\partial_{w}H_{m}(w)$ are T\"{o}plitz at $\infty$ and $\mathcal{NF}_{\Delta_{m}}$ for all $w\in U_{m}$,
\begin{equation*}
f_{m}=f_{m}^{low}+f_{m}^{high}
\end{equation*}
satisfy
\begin{equation*}\label{ii1}
|||X_{f_{m}^{low}}|||^{T}_{p,D(\rho_{m},\mu_{m},\sigma_{m})\times U_{m}}\leq\varepsilon_{m},
 \  |||X_{f_{m}^{high}}|||^{T}_{p,D(\rho_{m},\mu_{m},\sigma_{m})\times U_{m}}\leq \varepsilon+\sum_{j=1}^{m}\varepsilon_{j-1}^{\frac{2}{3}},
\end{equation*}
\begin{equation*}\label{ii2}
[f_{m}^{low}]_{\Lambda_{m},\gamma_{m},\sigma_{m};U_{m},\rho_{m},\mu_{m}}\leq\varepsilon_{m},
\ [f_{m}^{high}]_{\Lambda_{m},\gamma_{m},\sigma_{m};U_{m},\rho_{m},\mu_{m}}\leq \varepsilon+\sum_{j=1}^{m}\varepsilon_{j-1}^{\frac{2}{3}}.
\end{equation*}
Assume for all $w\in U_{m}$,
\begin{equation*}\label{ii3}
|\omega_{m}(w)-\omega(w)|+|\partial_{w}(\omega_{m}(w)-\omega(w))|\leq\sum_{j=1}^{m}\varepsilon_{j-1}^{\frac{2}{3}},
\end{equation*}
\begin{equation*}\label{ii4}
\|H_{m}-H\|_{U_{m}}+\langle H_{m}-H\rangle_{\Lambda_{m};U_{m}}\leq\sum_{j=1}^{m}\varepsilon_{j-1}^{\frac{2}{3}}.
\end{equation*}
Then there is a subset $U_{m+1}\subset U_{m}$ such that if
\begin{equation*}
\varepsilon\preceq\min\left(\gamma,\rho,\frac{1}{\Delta},\frac{1}{\Lambda}\right)^{\exp},
\end{equation*}
then for all $w\in U_{m+1}$, there is a real analytic symplectic map $\Phi_{m}$ such that
\begin{equation*}
(h_{m}+f_{m})\circ\Phi_{m}=h_{m+1}+f_{m+1}
\end{equation*}
with the estimates
\begin{equation*}\label{ii5}
|||X_{f^{low}_{m+1}}|||^{T}_{p,D(\rho_{m+1},\mu_{m+1},\sigma_{m+1})\times U_{m+1}}
\leq\varepsilon_{m+1},
\end{equation*}
\begin{equation*}\label{ii6}
|||X_{f^{high}_{m+1}}|||^{T}_{p,D(\rho_{m+1},\mu_{m+1},\sigma_{m+1})\times U_{m+1}}
\leq \varepsilon+\sum_{j=1}^{m+1}\varepsilon_{j-1}^{\frac{2}{3}},
\end{equation*}
\begin{equation*}\label{ii7}
[f^{low}_{m+1}]_{\Lambda_{m+1},\gamma_{m+1},\sigma_{m+1};U_{m+1},\rho_{m+1},\mu_{m+1}}
\leq\varepsilon_{m+1},
\end{equation*}
\begin{equation*}\label{ii8}
[f^{high}_{m+1}]_{\Lambda_{m+1},\gamma_{m+1},\sigma_{m+1};U_{m+1},\rho_{m+1},\mu_{m+1}}
\leq \varepsilon+\sum_{j=1}^{m+1}\varepsilon_{j-1}^{\frac{2}{3}},
\end{equation*}
\begin{equation*}\label{ii9}
|\omega_{m+1}(w)-\omega(w)|+|\partial_{w}(\omega_{m+1}(w)-\omega(w))|\leq\sum_{j=1}^{m+1}\varepsilon_{j-1}^{\frac{2}{3}},
\end{equation*}
\begin{equation*}\label{ii10}
\|H_{m+1}-H\|_{U_{m+1}}+\langle H_{m+1}-H\rangle_{\Lambda_{m+1};U_{m+1}}\leq\sum_{j=1}^{m+1}\varepsilon_{j-1}^{\frac{2}{3}},
\end{equation*}
\begin{equation*}\label{ii11}
\mathrm{meas}(U_{m}\setminus U_{m+1})\preceq \varepsilon_{m}^{\exp'},
\end{equation*}
where the exponents $\exp$, $\exp'$ depend on $d, \# \mathcal{A}, p$.

\end{lem}

\noindent\textbf{Proof.}~
Take $\kappa^{20}=\varepsilon^{\frac{1}{20}}$ in Lemma \ref{l1},
there is a real analytic symplectic map $\Phi$ such that
\begin{equation*}
(h+f)\circ\Phi=h+h_1+\cdots+h_{n}+f_{n}.
\end{equation*}
Let $h_{+}=h+h_1+\cdots+h_{n}$, $f_{+}=f_{n}$.
The KAM iterative lemma follows immediately from Lemma \ref{l1}.
\qed

\section{Long time stability of the KAM tori}\label{sect 5}

In this section, we prove
Theorem \ref{t2} on the long time stability of the KAM tori.
By momentum conservation, the frequency shift is diagonal.
We will construct a partial normal form of order $M+2$ based on
$h_{\tilde{m}}+f^{high}_{\tilde{m}}$, where
\begin{equation*}
h_{\tilde{m}}=\langle \omega_{\tilde{m}}(w),r\rangle+\frac{1}{2}\langle\zeta,(\Omega(w)+H_{\tilde{m}}(w))\zeta\rangle,
\end{equation*}
$H_{\tilde{m}}(w),\partial_{w}H_{\tilde{m}}(w)$ are T\"{o}plitz at $\infty$ and $\mathcal{NF}_{\Delta_{\tilde{m}}}$ for all $w\in U_{\infty}$.

We change to complex coordinates
\begin{equation*}
z=\left(
    \begin{array}{c}
      u \\
      v \\
    \end{array}
  \right)
  =C^{-1}\left(
           \begin{array}{c}
             \xi \\
             \eta \\
           \end{array}
         \right),
         \
         C=\frac{1}{\sqrt{2}}\left(
         \begin{array}{cc}
           1 & 1 \\
           -\mathrm{i} & \mathrm{i} \\
         \end{array}
       \right),
\end{equation*}
then
\begin{equation*}
h_{\tilde{m}}=\langle \omega_{\tilde{m}}(w),r\rangle+\langle u,(\Omega(w)+H_{\tilde{m}}(w))v\rangle.
\end{equation*}

Let $\mathcal{B}=\{a\in \mathcal{L}:|a|\leq N\}$, $\tilde{u}=(u_{a})_{a\in\mathcal{B}}$, $\tilde{v}=(v_{a})_{a\in\mathcal{B}}$,
$\check{u}=(u_{a})_{a\in\mathcal{L}\backslash\mathcal{B}}$,$\check{v}=(v_{a})_{a\in\mathcal{L}\backslash\mathcal{B}}$.
Write
\begin{equation*}
h_{\tilde{m}}=\langle \omega_{\tilde{m}}(w),r\rangle+\sum_{a\in\mathcal{B}}\tilde{\lambda}_{a}(w)u_{a}v_{a}+\langle \check{u},(\Omega(w)+H_{\tilde{m}}(w))\check{v}\rangle,
\end{equation*}
where $\tilde{\lambda}_{a}(w)\in\mathrm{spec}(\Omega(w)+H_{\tilde{m}}(w))_{\mathcal{B}}$.

Let $\tilde{\Delta}>1$ and $0<\tilde{\kappa}<1$. Assume there exists $\tilde{U} \subset U_{\infty}$ such that
for all $w\in \tilde{U}$, $|k|\leq\tilde{\Delta}$, $|\tilde{l}|\leq M+2$, $|k|+|\tilde{l}|\neq 0$, the following conditions hold:
\begin{itemize}
  \item Diophantine condition:
  \begin{equation}\label{sda1}
|\langle k,\omega_{\tilde{m}}(w)\rangle+\langle \tilde{l},\tilde{\lambda}(w)\rangle|\geq \frac{\tilde{\kappa}}{4^{M}N^{(4d)^{4d}(|\tilde{l}|+4)^{2}}};
\end{equation}

\item The first Melnikov condition:
\begin{equation}\label{sda2}
|\langle k,\omega_{\tilde{m}}(w)\rangle+\langle \tilde{l},\tilde{\lambda}(w)\rangle+\alpha(w)|\geq \frac{\tilde{\kappa}}{4^{M}N^{(4d)^{4d}(|\tilde{l}|+4)^{2}}}
\end{equation}
for any $\alpha(w)\in\mathrm{spec}(((\Omega+H_{\tilde{m}})(w))_{[a]_{\Delta_{\tilde{m}}}})$
 and any $[a]_{\Delta_{\tilde{m}}}$;

 \item The second Melnikov condition with the same sign:
 \begin{equation}\label{sda3}
|\langle k,\omega_{\tilde{m}}(w)\rangle+\langle \tilde{l},\tilde{\lambda}(w)\rangle+\alpha(w)+\beta(w)|\geq \frac{\tilde{\kappa}}{4^{M}N^{(4d)^{4d}(|\tilde{l}|+4)^{2}}}
\end{equation}
for any $\alpha(w)\in\mathrm{spec}(((\Omega+H_{\tilde{m}})(w))_{[a]_{\Delta_{\tilde{m}}}}),
\beta(w)\in\mathrm{spec}(((\Omega+H_{\tilde{m}})(w))_{[b]_{\Delta_{\tilde{m}}}})$
and any $[a]_{\Delta_{\tilde{m}}}, [b]_{\Delta_{\tilde{m}}}$;

 \item The second Melnikov condition with the opposite signs:
 \begin{equation}\label{sda4}
|\langle k,\omega_{\tilde{m}}(w)\rangle+\langle \tilde{l},\tilde{\lambda}(w)\rangle+\alpha(w)-\beta(w)|\geq \frac{\tilde{\kappa}}{4^{M}N^{(4d)^{4d}(|\tilde{l}|+4)^{2}}}
\end{equation}
for any $ \alpha(w)\in\mathrm{spec}(((\Omega+H_{\tilde{m}})(w))_{[a]_{\Delta_{\tilde{m}}}}),
\beta(w)\in\mathrm{spec}(((\Omega+H_{\tilde{m}})(w))_{[b]_{\Delta_{\tilde{m}}}}),
$ and any $\mathrm{dist}([a]_{\Delta_{\tilde{m}}}, [b]_{\Delta_{\tilde{m}}})\leq \tilde{\Delta}+2d_{\Delta_{\tilde{m}}}$.
\end{itemize}

We have the following lemma.
For the sake of notations, we denote
\begin{equation*}
  \bm{\sum_{(1)}} = \sum_{2|\alpha|+|\beta|+|\upsilon|=j} ,\quad
    \bm{\sum_{(2)}}= \sum_{2|\alpha|+|\beta|+|\upsilon|=j-1},\quad
     \bm{\sum_{(3)}}= \sum_{2|\alpha|+|\beta|+|\upsilon|=j-2}.
\end{equation*}

\begin{lem}\label{l3}

For $2\leq j_0\leq M+1$, consider the partial normal form of order $j_0$
\begin{equation*}
T_{j_{0}}=h_{\tilde{m}}+Z_{j_{0}}+P_{j_{0}}+R_{j_{0}}+Q_{j_{0}},
\end{equation*}
with
\begin{equation*}
Z_{j_{0}}=\sum_{3\leq j\leq j_{0}}Z_{j_{0}j},~ P_{j_{0}}=\sum_{j\geq j_{0}+1}P_{j_{0}j},~
R_{j_{0}}=\sum_{3\leq j\leq j_{0}}R_{j_{0}j},~ Q_{j_{0}}=\sum_{j\geq 3}Q_{j_{0}j},
\end{equation*}
where
\begin{itemize}
  \item $Z_{j_{0}j}=Z_{j_{0}j}(r,u,v;w)$ equals to
  \begin{equation*}
\sum_{2|\alpha|+2|\beta|=j}Z_{j_{0}}^{\alpha\beta\beta0}(w)r^{\alpha}\tilde{u}^{\beta}\tilde{v}^{\beta}+
\sum_{2|\alpha|+2|\beta|=j-2,|a|=|b|}Z_{j_{0}}^{\alpha\beta\beta ab}(w)r^{\alpha}\tilde{u}^{\beta}\tilde{v}^{\beta}\check{u}_{a}\check{v}_{b},
\end{equation*}

\item $P_{j_{0}j}=P_{j_{0}j}(\varphi,r,u,v;w)$ equals to
\begin{align*}
&\bm{\sum_{(1)}} P_{j_{0}}^{\alpha\beta\upsilon0}(\varphi;w)r^{\alpha}\tilde{u}^{\beta}\tilde{v}^{\upsilon} +
\bm{\sum_{(2)}} \langle \check{u},P_{j_{0}}^{\alpha\beta\upsilon\check{u}}(\varphi;w)\rangle r^{\alpha}\tilde{u}^{\beta}\tilde{v}^{\upsilon}\\
+&\bm{\sum_{(2)}} \langle \check{v},P_{j_{0}}^{\alpha\beta\upsilon\check{v}}(\varphi;w)\rangle r^{\alpha}\tilde{u}^{\beta}\tilde{v}^{\upsilon}
+\bm{\sum_{(3)}} \frac{1}{2}\langle \check{u},P_{j_{0}}^{\alpha\beta\upsilon\check{u}\check{u}}(\varphi;w)\check{u}\rangle r^{\alpha}\tilde{u}^{\beta}\tilde{v}^{\upsilon}\\
+&\bm{\sum_{(3)}} \langle \check{u},P_{j_{0}}^{\alpha\beta\upsilon\check{u}\check{v}}(\varphi;w)\check{v}\rangle r^{\alpha}\tilde{u}^{\beta}\tilde{v}^{\upsilon}
+\bm{\sum_{(3)}} \frac{1}{2}\langle \check{v},P_{j_{0}}^{\alpha\beta\upsilon\check{v}\check{v}}(\varphi;w)\check{v}\rangle r^{\alpha}\tilde{u}^{\beta}\tilde{v}^{\upsilon},
\end{align*}

\item $R_{j_{0}j}=R_{j_{0}j}(\varphi,r,u,v;w)$ equals to
\begin{align*}
&\bm{\sum_{(1)}}R_{j_{0}}^{\alpha\beta\upsilon0}(\varphi;w)r^{\alpha}\tilde{u}^{\beta}\tilde{v}^{\upsilon}+
\bm{\sum_{(2)}}\langle \check{u},R_{j_{0}}^{\alpha\beta\upsilon\check{u}}(\varphi;w)\rangle r^{\alpha}\tilde{u}^{\beta}\tilde{v}^{\upsilon}
\\
+&\bm{\sum_{(2)}}\langle \check{v},R_{j_{0}}^{\alpha\beta\upsilon\check{v}}(\varphi;w)\rangle r^{\alpha}\tilde{u}^{\beta}\tilde{v}^{\upsilon}
+\bm{\sum_{(3)}}\frac{1}{2}\langle \check{u},R_{j_{0}}^{\alpha\beta\upsilon\check{u}\check{u}}(\varphi;w)\check{u}\rangle r^{\alpha}\tilde{u}^{\beta}\tilde{v}^{\upsilon}
\\
+&\bm{\sum_{(3)}}\langle \check{u},R_{j_{0}}^{\alpha\beta\upsilon\check{u}\check{v}}(\varphi;w)\check{v}\rangle r^{\alpha}\tilde{u}^{\beta}\tilde{v}^{\upsilon}
+\bm{\sum_{(3)}}\frac{1}{2}\langle \check{v},R_{j_{0}}^{\alpha\beta\upsilon\check{v}\check{v}}(\varphi;w)\check{v}\rangle r^{\alpha}\tilde{u}^{\beta}\tilde{v}^{\upsilon},
\end{align*}

\item $Q_{j_{0}j}=Q_{j_{0}j}(\varphi,r,u,v;w)$ equals to
  \begin{equation*}
\sum_{2|\alpha|+|\beta|+|\upsilon|+|\mu|+|\nu|=j,|\mu|+|\nu|\geq3}Q_{j_{0}}^{\alpha\beta\upsilon\mu\nu}(\varphi;w)r^{\alpha}\tilde{u}^{\beta}\tilde{v}^{\upsilon}\check{u}^{\mu}\check{v}^{\nu}.
\end{equation*}
\end{itemize}

Let $\rho'=\frac{\rho}{12M}, \delta'=\frac{\delta}{2M}$ and
\begin{equation*}
D_{j_{0}}=D\left(\frac{\rho}{2}-3(j_{0}-2)\rho',
(5\delta-2(j_{0}-2)\delta')^{2},5\delta-2(j_{0}-2)\delta'\right),
\end{equation*}
\begin{equation*}
D'_{j_{0}}=D\left(\frac{\rho}{2}-(3(j_{0}-2)+1)\rho',
(5\delta-2(j_{0}-2)\delta')^{2},5\delta-2(j_{0}-2)\delta'\right).
\end{equation*}

Assume
\begin{equation}\label{nf1}
|||X_{Z_{j_{0}j}}|||^{T}_{p,D_{j_{0}}\times\tilde{U}}\preceq\delta\left(\frac{1}{\tilde{\kappa}^{2}}\Delta_{\tilde{m}}^{\exp}N^{(4d)^{4d}(j_{0}+4)^{2}}\delta\right)^{j-3},
\end{equation}
\begin{equation}\label{nf2}
|||X_{R_{j_{0}j}}|||^{T}_{p,D_{j_{0}}\times\tilde{U}}\preceq\delta\left(\frac{1}{\tilde{\kappa}^{2}}\Delta_{\tilde{m}}^{\exp}N^{(4d)^{4d}(j_{0}+4)^{2}}\delta\right)^{j-3},
\end{equation}
\begin{equation}\label{nf3}
|||X_{P_{j_{0}j}}|||^{T}_{p,D_{j_{0}}\times\tilde{U}}\preceq\delta\left(\frac{1}{\tilde{\kappa}^{2}}\Delta_{\tilde{m}}^{\exp}N^{(4d)^{4d}(j_{0}+5)^{2}}\delta\right)^{j-3},
\end{equation}
\begin{equation}\label{nf4}
|||X_{Q_{j_{0}j}}|||^{T}_{p,D_{j_{0}}\times\tilde{U}}\preceq\delta\left(\frac{1}{\tilde{\kappa}^{2}}\Delta_{\tilde{m}}^{\exp}N^{(4d)^{4d}(j_{0}+5)^{2}}\delta\right)^{j-3},
\end{equation}
where $a\preceq b$ means there is a constant $c>0$ depending on $\rho,M,d,\# \mathcal{A}, p, c_1,c_2,c_3,c_4,c_5$ such that $a\leq c b$,
the exponent $\exp$ depends on $d, \# \mathcal{A}, p$.

Then there is a symplectic map $\Psi_{j_{0}}$ such that
\begin{equation*}
T_{j_{0}+1}=T_{j_{0}}\circ\Psi_{j_{0}}=h_{\tilde{m}}+Z_{j_{0}+1}+P_{j_{0}+1}+R_{j_{0}+1}+Q_{j_{0}+1},
\end{equation*}
which
is given exactly by the formula of $T_{j_{0}}$ but with $j_{0}+1$ in place of $j_{0}$.
Moreover,
the estimates \eqref{nf1}-\eqref{nf4} also hold  with $j_{0}+1$ in place of $j_{0}$,

\end{lem}

\noindent\textbf{Proof.}~
Consider the homological equation
\begin{equation}\label{n1}
\{h_{\tilde{m}},F_{j_{0}}\}=-\mathcal{T}_{\tilde{\Delta}}P_{j_{0}(j_{0}+1)}+\hat{Z}_{j_{0}},
\end{equation}
where
\begin{align*}
\mathcal{T}_{\tilde{\Delta}}P_{j_{0}(j_{0}+1)}
= & \sum_{|k|\leq\tilde{\Delta}}\Big[\sum_{2|\alpha|+|\beta|+|\upsilon|=j_{0}+1}
\hat{P}_{j_{0}}^{\alpha\beta\upsilon0}(k;w)r^{\alpha}\tilde{u}^{\beta}\tilde{v}^{\upsilon}\\
+& \sum_{2|\alpha|+|\beta|+|\upsilon|=j_{0}} \left( \langle \check{u},\hat{P}_{j_{0}}^{\alpha\beta\upsilon\check{u}}(k;w)\rangle r^{\alpha}\tilde{u}^{\beta}\tilde{v}^{\upsilon}
+\langle \check{v},\hat{P}_{j_{0}}^{\alpha\beta\upsilon\check{v}}(k;w)\rangle r^{\alpha}\tilde{u}^{\beta}\tilde{v}^{\upsilon}\right)\\
+&\sum_{2|\alpha|+|\beta|+|\upsilon|=j_{0}-1}\left(\frac{1}{2}\langle \check{u},\hat{P}_{j_{0}}^{\alpha\beta\upsilon\check{u}\check{u}}(k;w)\check{u}\rangle r^{\alpha}\tilde{u}^{\beta}\tilde{v}^{\upsilon}
+\frac{1}{2}\langle \check{v},\hat{P}_{j_{0}}^{\alpha\beta\upsilon\check{v}\check{v}}(k;w)\check{v}\rangle r^{\alpha}\tilde{u}^{\beta}\tilde{v}^{\upsilon}\right.\\
&\left.+\langle \check{u},\mathcal{T}_{\tilde{\Delta}}\hat{P}_{j_{0}}^{\alpha\beta\upsilon\check{u}\check{v}}
(k;w)\check{v}\rangle\right) r^{\alpha}\tilde{u}^{\beta}\tilde{v}^{\upsilon}\Big]
e^{\mathrm{i}\langle k,\varphi\rangle}.
\end{align*}
Let
\begin{align*}
F_{j_{0}}(\varphi,r,&u,v;w)=\sum_{2|\alpha|+|\beta|+|\upsilon|=j_{0}+1}
F_{j_{0}}^{\alpha\beta\upsilon0}(\varphi;w)r^{\alpha}\tilde{u}^{\beta}\tilde{v}^{\upsilon}\\
& +\sum_{2|\alpha|+|\beta|+|\upsilon|=j_{0}}\langle \check{u},F_{j_{0}}^{\alpha\beta\upsilon\check{u}}(\varphi;w)\rangle r^{\alpha}\tilde{u}^{\beta}\tilde{v}^{\upsilon}
+\langle \check{v},F_{j_{0}}^{\alpha\beta\upsilon\check{v}}(\varphi;w)\rangle r^{\alpha}\tilde{u}^{\beta}\tilde{v}^{\upsilon}\\
&+\sum_{2|\alpha|+|\beta|+|\upsilon|=j_{0}-1}\frac{1}{2}\langle \check{u},F_{j_{0}}^{\alpha\beta\upsilon\check{u}\check{u}}(\varphi;w)\check{u}\rangle r^{\alpha}\tilde{u}^{\beta}\tilde{v}^{\upsilon}
+\langle \check{u},F_{j_{0}}^{\alpha\beta\upsilon\check{u}\check{v}}(\varphi;w)\check{v}\rangle r^{\alpha}\tilde{u}^{\beta}\tilde{v}^{\upsilon}\\
&+\sum_{2|\alpha|+|\beta|+|\upsilon|=j_{0}-1}\frac{1}{2}\langle \check{v},F_{j_{0}}^{\alpha\beta\upsilon\check{v}\check{v}}(\varphi;w)\check{v}\rangle r^{\alpha}\tilde{u}^{\beta}\tilde{v}^{\upsilon}.
\end{align*}
In Fourier modes, we have
\begin{equation}\label{n2}
-\mathrm{i}(\langle k,\omega_{\tilde{m}}\rangle+\langle \beta-\upsilon,\tilde{\lambda}\rangle)\hat{F}_{j_{0}}^{\alpha\beta\upsilon0}(k)
=-\hat{P}_{j_{0}}^{\alpha\beta\upsilon0}(k)+\delta_{0}^{k}\delta_{\beta}^{\upsilon}
\hat{Z}_{j_{0}}^{\alpha\beta\upsilon0},\quad
\end{equation}
\begin{equation}\label{n3}
-\mathrm{i}(\langle k,\omega_{\tilde{m}}\rangle+\langle \beta-\upsilon,\tilde{\lambda}\rangle)\hat{F}_{j_{0}}^{\alpha\beta\upsilon\check{u}}(k)
-\mathrm{i}(\Omega+H_{\tilde{m}})\hat{F}_{j_{0}}^{\alpha\beta\upsilon\check{u}}(k)
=-\hat{P}_{j_{0}}^{\alpha\beta\upsilon\check{u}}(k),
\end{equation}
\begin{equation}\label{n4}
-\mathrm{i}(\langle k,\omega_{\tilde{m}}\rangle+\langle \beta-\upsilon,\tilde{\lambda}\rangle)\hat{F}_{j_{0}}^{\alpha\beta\upsilon\check{v}}(k)
+\mathrm{i}(\Omega+H_{\tilde{m}}^{T})\hat{F}_{j_{0}}^{\alpha\beta\upsilon\check{v}}(k)
=-\hat{P}_{j_{0}}^{\alpha\beta\upsilon\check{v}}(k),
\end{equation}
\begin{equation}\label{n5}
\begin{aligned}
-\mathrm{i}(\langle k,\omega_{\tilde{m}}\rangle+\langle \beta-\upsilon,\tilde{\lambda}\rangle)
&\hat{F}_{j_{0}}^{\alpha\beta\upsilon\check{u}\check{u}}(k)
-\mathrm{i}(\Omega+H_{\tilde{m}})
\hat{F}_{j_{0}}^{\alpha\beta\upsilon\check{u}\check{u}}(k)\\
-\mathrm{i}&\hat{F}_{j_{0}}^{\alpha\beta\upsilon\check{u}\check{u}}(k)
(\Omega+H_{\tilde{m}}^{T})
=-\hat{P}_{j_{0}}^{\alpha\beta\upsilon\check{u}\check{u}}(k),
\end{aligned}
\end{equation}
\begin{equation}\label{n6}
\begin{aligned}
-\mathrm{i}(\langle k,\omega_{\tilde{m}}\rangle+\langle \beta-\upsilon,\tilde{\lambda}\rangle)
&\hat{F}_{j_{0}}^{\alpha\beta\upsilon\check{v}\check{v}}(k)
+\mathrm{i}(\Omega+H_{\tilde{m}}^{T})
\hat{F}_{j_{0}}^{\alpha\beta\upsilon\check{v}\check{v}}(k)\\
+\mathrm{i}&\hat{F}_{j_{0}}^{\alpha\beta\upsilon\check{v}\check{v}}(k)(\Omega+H_{\tilde{m}})
=-\hat{P}_{j_{0}}^{\alpha\beta\upsilon\check{v}\check{v}}(k),
\end{aligned}
\end{equation}
\begin{equation}\label{n7}
\begin{aligned}
-\mathrm{i}(\langle k,\omega_{\tilde{m}}\rangle+\langle \beta-\upsilon,\tilde{\lambda}\rangle)
&\hat{F}_{j_{0}}^{\alpha\beta\upsilon\check{u}\check{v}}(k)
+\mathrm{i}\hat{F}_{j_{0}}^{\alpha\beta\upsilon\check{u}\check{v}}(k)(\Omega
+H_{\tilde{m}})\\
-\mathrm{i}(\Omega+H_{\tilde{m}})&\hat{F}_{j_{0}}^{\alpha\beta\upsilon\check{u}\check{v}}(k)
=-\hat{P}_{j_{0}}^{\alpha\beta\upsilon\check{u}\check{v}}(k)
+\delta_{0}^{k}\delta_{\beta}^{\upsilon}
\hat{Z}_{j_{0}}^{\alpha\beta\upsilon\check{u}\check{v}}.
\end{aligned}
\end{equation}
We solve equations (\ref{n2})-(\ref{n7}) as in Proposition \ref{p1}.
We have
\begin{equation}\label{n8}
\begin{aligned}
\hat{Z}_{j_{0}}(r,u,v;w)=&\sum_{2|\alpha|+2|\beta|=j_{0}+1}
\hat{P}_{j_{0}}^{\alpha\beta\beta0}(0;w)r^{\alpha}\tilde{u}^{\beta}\tilde{v}^{\beta}\\
+&
\sum_{2|\alpha|+2|\beta|=j_{0}-1,|a|=|b|}\hat{P}_{j_{0}}^{\alpha\beta\beta ab}(0;w)r^{\alpha}\tilde{u}^{\beta}\tilde{v}^{\beta}\check{u}_{a}\check{v}_{b}.
\end{aligned}
\end{equation}
Let $\Psi_{j_{0}}=X_{F_{j_{0}}}^{t}|_{t=1}$. Using Taylor's formula, there is
\begin{equation}\label{n9}
\begin{aligned}
T_{j_{0}+1}=&T_{j_{0}}\circ\Psi_{j_{0}}=(h_{\tilde{m}}+Z_{j_{0}}+P_{j_{0}}+R_{j_{0}}+Q_{j_{0}})
\circ X_{F_{j_{0}}}^{t}|_{t=1}\\
=&h_{\tilde{m}}+\{h_{\tilde{m}},F_{j_{0}}\}+\int_{0}^{1}(1-t)
\{\{h_{\tilde{m}},F_{j_{0}}\},F_{j_{0}}\}\circ X_{F_{j_{0}}}^{t}dt\\
&+Z_{j_{0}}+\int_{0}^{1}\{Z_{j_{0}},F_{j_{0}}\}\circ X_{F_{j_{0}}}^{t}dt
+P_{j_{0}(j_{0}+1)}+\int_{0}^{1}\{P_{j_{0}(j_{0}+1)},F_{j_{0}}\}\circ X_{F_{j_{0}}}^{t}dt\\
&+P_{j_{0}}-P_{j_{0}(j_{0}+1)}+\int_{0}^{1}\{P_{j_{0}}-P_{j_{0}(j_{0}+1)},F_{j_{0}}\}\circ X_{F_{j_{0}}}^{t}dt\\
&+R_{j_{0}}+\int_{0}^{1}\{R_{j_{0}},F_{j_{0}}\}\circ X_{F_{j_{0}}}^{t}dt+Q_{j_{0}}+\int_{0}^{1}\{Q_{j_{0}},F_{j_{0}}\}\circ X_{F_{j_{0}}}^{t}dt.
\end{aligned}
\end{equation}
By (\ref{n1}) and (\ref{n9}), we have
\begin{equation}\label{n10}
\begin{aligned}
T_{j_{0}+1}=&h_{\tilde{m}}+Z_{j_{0}}+\hat{Z}_{j_{0}}
+\int_{0}^{1}(1-t)\{\{h_{\tilde{m}},F_{j_{0}}\},F_{j_{0}}\}\circ X_{F_{j_{0}}}^{t}dt\\
&+\int_{0}^{1}\{Z_{j_{0}},F_{j_{0}}\}\circ X_{F_{j_{0}}}^{t}dt
+\int_{0}^{1}\{P_{j_{0}(j_{0}+1)},F_{j_{0}}\}\circ X_{F_{j_{0}}}^{t}dt\\
&+P_{j_{0}}-P_{j_{0}(j_{0}+1)}+\int_{0}^{1}\{P_{j_{0}}-P_{j_{0}(j_{0}+1)},F_{j_{0}}\}\circ X_{F_{j_{0}}}^{t}dt\\
&+R_{j_{0}}+(1-\mathcal{T}_{\tilde{\Delta}})P_{j_{0}(j_{0}+1)}
+\int_{0}^{1}\{R_{j_{0}},F_{j_{0}}\}\circ X_{F_{j_{0}}}^{t}dt+Q_{j_{0}}\\
&+\int_{0}^{1}\{Q_{j_{0}},F_{j_{0}}\}\circ X_{F_{j_{0}}}^{t}dt.
\end{aligned}
\end{equation}
Hence
\begin{equation}\label{n11}
Z_{j_{0}+1}=Z_{j_{0}}+\hat{Z}_{j_{0}},
\end{equation}
\begin{equation}\label{n12}
R_{j_{0}+1}=R_{j_{0}}+(1-\mathcal{T}_{\tilde{\Delta}})P_{j_{0}(j_{0}+1)},
\end{equation}
and
\begin{equation}\label{n13}
\begin{aligned}
& P_{j_{0}+1}+Q_{j_{0}+1}=\int_{0}^{1}(1-t)\{\{h_{\tilde{m}},F_{j_{0}}\},F_{j_{0}}\}\circ X_{F_{j_{0}}}^{t}dt
+\int_{0}^{1}\{Z_{j_{0}}+P_{j_{0}},F_{j_{0}}\}\circ X_{F_{j_{0}}}^{t}dt\\
&+\int_{0}^{1}\{R_{j_{0}},F_{j_{0}}\}\circ X_{F_{j_{0}}}^{t}dt+P_{j_{0}}-P_{j_{0}(j_{0}+1)}
+Q_{j_{0}}
+\int_{0}^{1}\{Q_{j_{0}},F_{j_{0}}\}\circ X_{F_{j_{0}}}^{t}dt.
\end{aligned}
\end{equation}

By (\ref{sda1})-(\ref{sda4}) and (\ref{nf3}), following the proof of Proposition \ref{p1}, we have
\begin{equation}\label{n14}
\begin{aligned}
|||X_{F_{j_{0}}}|||^{T}_{p,D'_{j_{0}}\times\tilde{U}}\preceq & \frac{1}{\tilde{\kappa}^{2}}
\Delta_{\tilde{m}}^{\exp}N^{(4d)^{4d}(j_{0}+5)^{2}}
|||X_{P_{j_{0}(j_{0}+1)}}|||^{T}_{p,D_{j_{0}}\times\tilde{U}}\\
\preceq & \left(\frac{1}{\tilde{\kappa}^{2}}\Delta_{\tilde{m}}^{\exp}N^{(4d)^{4d}(j_{0}+5)^{2}}\delta\right)^{j_{0}-1}.
\end{aligned}
\end{equation}

Write
\begin{equation*}
Z_{j_{0}+1}=\sum_{3\leq j\leq j_{0}+1}Z_{(j_{0}+1)j},
\end{equation*}
where
\begin{equation*}
Z_{(j_{0}+1)j}=Z_{j_{0}j}, \ 3\leq j\leq j_{0},
\end{equation*}
\begin{equation*}
Z_{(j_{0}+1)(j_{0}+1)}=\hat{Z}_{j_{0}}.
\end{equation*}
For $3\leq j\leq j_{0}$, using (\ref{nf1}), we have
\begin{equation}\label{n15}
|||X_{Z_{(j_{0}+1)j}}|||^{T}_{p,D_{j_{0}+1}\times\tilde{U}}\preceq\delta\left(\frac{1}{\tilde{\kappa}^{2}}\Delta_{\tilde{m}}^{\exp}N^{(4d)^{4d}(j_{0}+5)^{2}}\delta\right)^{j-3}.
\end{equation}
By (\ref{nf3}) and (\ref{n8}), we have
\begin{equation}\label{n16}
|||X_{\hat{Z}_{j_{0}}}|||^{T}_{p,D_{j_{0}+1}\times\tilde{U}}\preceq\delta\left(\frac{1}{\tilde{\kappa}^{2}}\Delta_{\tilde{m}}^{\exp}N^{(4d)^{4d}(j_{0}+5)^{2}}\delta\right)^{j_{0}-2}.
\end{equation}
 Then the estimate of (\ref{n11}) follows from (\ref{n15}) and (\ref{n16}).

Write
\begin{equation*}
R_{j_{0}+1}=\sum_{3\leq j\leq j_{0}+1}R_{(j_{0}+1)j},
\end{equation*}
where
\begin{equation*}
R_{(j_{0}+1)j}=R_{j_{0}j}, \ 3\leq j\leq j_{0},
\end{equation*}
\begin{equation*}
R_{(j_{0}+1)(j_{0}+1)}=(1-\mathcal{T}_{\tilde{\Delta}})P_{j_{0}(j_{0}+1)}.
\end{equation*}
For $3\leq j\leq j_{0}$, using (\ref{nf2}), we have
\begin{equation}\label{n17}
|||X_{R_{(j_{0}+1)j}}|||^{T}_{p,D_{j_{0}+1}\times\tilde{U}}\preceq\delta\left(\frac{1}{\tilde{\kappa}^{2}}\Delta_{\tilde{m}}^{\exp}N^{(4d)^{4d}(j_{0}+5)^{2}}\delta\right)^{j-3}.
\end{equation}
By (\ref{nf3}), we have
\begin{equation}\label{n18}
|||X_{R_{(j_{0}+1)(j_{0}+1)}}|||^{T}_{p,D_{j_{0}+1}\times\tilde{U}}\preceq\delta\left(\frac{1}{\tilde{\kappa}^{2}}\Delta_{\tilde{m}}^{\exp}N^{(4d)^{4d}(j_{0}+5)^{2}}\delta\right)^{j_{0}-2}.
\end{equation}
By (\ref{n17}) and (\ref{n18}), we obtain the estimate of (\ref{n12}).

Let
\begin{equation*}
h_{\tilde{m}}^{(0)}=h_{\tilde{m}}, \ h_{\tilde{m}}^{(j)}=\{h_{\tilde{m}}^{(j-1)},F_{j_{0}}\}, ~ j\geq1.
\end{equation*}
We have
\begin{equation*}
h_{\tilde{m}}^{(j)}(\varphi,r,u,v;w)=\sum_{\substack{2|\alpha|+|\beta|+|\upsilon|+|\mu|+|\nu|\\
=j(j_{0}-1)+2}}
h_{\tilde{m}}^{(j)\alpha\beta\upsilon\mu\nu}(\varphi;w)r^{\alpha}\tilde{u}^{\beta}\tilde{v}^{\upsilon}\check{u}^{\mu}\check{v}^{\nu},
\end{equation*}
and
\begin{equation}\label{n19}
\int_{0}^{1}(1-t)\{\{h_{\tilde{m}},F_{j_{0}}\},F_{j_{0}}\}\circ X_{F_{j_{0}}}^{t}dt=\sum_{j\geq2}\frac{1}{j!}h_{\tilde{m}}^{(j)}.
\end{equation}

By (\ref{nf3}), (\ref{n1}) and (\ref{n8}), we have
\begin{equation}\label{n20}
|||X_{h_{\tilde{m}}^{(1)}}|||^{T}_{p,D_{j_{0}}\times\tilde{U}}\preceq\delta\left(\frac{1}{\tilde{\kappa}^{2}}\Delta_{\tilde{m}}^{\exp}N^{(4d)^{4d}(j_{0}+5)^{2}}\delta\right)^{j_{0}-2}.
\end{equation}
For $j\geq2$, let $\rho_{j}=\frac{\rho}{6jM}$, $\delta_{j}=\frac{\delta}{2jM}$. By (\ref{n14}), (\ref{n20}) and Proposition \ref{p}, we have
\begin{equation}\label{n21}
\begin{aligned}
&~ \frac{1}{j!}|||X_{h_{\tilde{m}}^{(j)}}|||^{T}_{p,D_{j_{0}+1}\times\tilde{U}} \\
\preceq &  \frac{1}{j!}\left(C\max\left(\frac{1}{\rho_{j}},\frac{\delta}{\delta_{j}}\right)\right)^{j-1}
\left(|||X_{F_{j_{0}}}|||^{T}_{p,D'_{j_{0}}\times\tilde{U}}\right)^{j-1}|||X_{h_{\tilde{m}}^{(1)}}|||^{T}_{p,D_{j_{0}}\times\tilde{U}}\\
\preceq & \frac{j^{j-1}}{j!}\delta\left(\frac{1}{\tilde{\kappa}^{2}}\Delta_{\tilde{m}}^{\exp}N^{(4d)^{4d}(j_{0}+5)^{2}}\delta\right)^{j_{0}-2}
\left(C\left(\frac{1}{\tilde{\kappa}^{2}}\Delta_{\tilde{m}}^{\exp}N^{(4d)^{4d}(j_{0}+5)^{2}}\delta\right)^{j_{0}-1}\right)^{j-1}\\
\preceq & \delta\left(\frac{1}{\tilde{\kappa}^{2}}\Delta_{\tilde{m}}^{\exp}
N^{(4d)^{4d}(j_{0}+6)^{2}}\delta\right)^{j(j_{0}-1)-1},
\end{aligned}
\end{equation}
since $j, j_{0}\geq2$, $j(j_{0}-1)+2\geq2j_{0}\geq j_{0}+2$. By (\ref{n21}), we obtain the estimate of (\ref{n19}).

Let
\begin{equation*}
W_{i}=Z_{j_{0}i}, \ 3\leq i\leq j_{0}, \ W_{i}=P_{j_{0}i}    \  i\geq j_{0}+1.
\end{equation*}
We have
\begin{equation*}
Z_{j_{0}}+P_{j_{0}}=\sum_{i\geq3}W_{i}.
\end{equation*}
Let
\begin{equation*}
W_{i}^{(0)}=W_{i},~ \ W_{i}^{(j)}=\{W_{i}^{(j-1)},F_{j_{0}}\},~ j\geq1.
\end{equation*}
We have
\begin{equation*}
W_{i}^{(j)}(\varphi,r,u,v;w)=\sum_{2|\alpha|+|\beta|+|\upsilon|+|\mu|+|\nu|=j(j_{0}-1)+i}
W_{i}^{(j)\alpha\beta\upsilon\mu\nu}(\varphi;w)r^{\alpha}\tilde{u}^{\beta}\tilde{v}^{\upsilon}\check{u}^{\mu}\check{v}^{\nu},
\end{equation*}
and
\begin{equation}\label{n22}
\int_{0}^{1}\{Z_{j_{0}}+P_{j_{0}},F_{j_{0}}\}\circ X_{F_{j_{0}}}^{t}dt=\sum_{i\geq3}\sum_{j\geq1}\frac{1}{j!}W_{i}^{(j)}.
\end{equation}
Following the proof of (\ref{n21}),  we have
\begin{equation}\label{n23}
\frac{1}{j!}|||X_{W_{i}^{(j)}}|||^{T}_{p,D_{j_{0}+1}\times\tilde{U}}\preceq
\delta\left(\frac{1}{\tilde{\kappa}^{2}}
\Delta_{\tilde{m}}^{\exp}N^{(4d)^{4d}(j_{0}+6)^{2}}\delta\right)^{j(j_{0}-1)+i-3},
\end{equation}
since $j\geq1, j_{0}\geq2, i\geq3$, $j(j_{0}-1)+i\geq j_{0}+2$. By (\ref{n23}), we obtain the estimate of (\ref{n22}).
Using the same method, we can estimate
\begin{equation*}
\int_{0}^{1}\{R_{j_{0}},F_{j_{0}}\}\circ X_{F_{j_{0}}}^{t}dt
+\int_{0}^{1}\{Q_{j_{0}},F_{j_{0}}\}\circ X_{F_{j_{0}}}^{t}dt.
\end{equation*}
Hence, we obtain the estimate of (\ref{n13}).

\qed

Now we  prove  Theorem \ref{t2} concerning the long
time stability of the KAM tori of  the infinitely dimensional
Hamiltonian system.

\noindent\textbf{Proof of Theorem \ref{t2}.}~
Since the difference between $h_{\infty}+f_{\infty}$ and $h_{\tilde{m}}+f^{high}_{\tilde{m}}$ is $\varepsilon^{\frac{2}{3}}_{\tilde{m}}$,
if we choose $\tilde{m}$ such that $\varepsilon^{\frac{2}{3}}_{\tilde{m}}\sim \delta^{M+1}$, then
we can construct a partial normal form of order $M+2$ based on $h_{\tilde{m}}+f^{high}_{\tilde{m}}$.

By Lemma \ref{l3},
there is a symplectic map $\Psi$ such that
\begin{equation*}
(h_{\tilde{m}}+f^{high}_{\tilde{m}})\circ\Psi=h_{\tilde{m}}+Z_{M+2}+P_{M+2}+R_{M+2}+Q_{M+2}.
\end{equation*}

Taking $N=\delta^{-\frac{M+1}{p-1}}$, $\tilde{\kappa}=\delta^{\frac{1}{900M}}$, we have
\begin{equation*}\label{n24}
|||X_{P_{M+2}}|||^{T}_{p,D(\frac{\rho}{4},(4\delta)^{2},4\delta)\times\tilde{U}}\preceq\delta\left(\frac{1}{\tilde{\kappa}^{2}}\Delta_{\tilde{m}}^{\exp}N^{(4d)^{4d}(M+7)^{2}}\delta\right)^{M}\leq\delta^{M+\frac{1}{2}}.
\end{equation*}
Taking $\tilde{\Delta}=800M(\log\frac{1}{\delta})^{2}\frac{1}{\min(\gamma_{\tilde{m}},\rho)}$, we have
\begin{equation*}\label{n25}
|||X_{R_{M+2}}|||^{T}_{p,D(\frac{\rho}{4},(4\delta)^{2},4\delta)\times\tilde{U}}\leq\delta^{M+\frac{1}{2}}.
\end{equation*}

Since
\begin{equation*}
\|\check{z}\|_{1}=\left(\sum_{|a|>N}|z_{a}|^{2}|a|^{2}\right)^{\frac{1}{2}}=\left(\sum_{|a|>N}|z_{a}|^{2}\frac{|a|^{2p}}{|a|^{2p-2}}\right)^{\frac{1}{2}}\leq\frac{\|\check{z}\|_{p}}{N^{p-1}}\leq
\delta^{M+1}\|\check{z}\|_{p},
\end{equation*}
we have
\begin{equation*}\label{n26}
|||X_{Q_{M+2}}|||_{\mathcal{P}^{p},D(\frac{\rho}{4},(4\delta)^{2},4\delta)\times\tilde{U}}\leq\delta^{M+\frac{1}{2}}.
\end{equation*}

As done in Section 5.3 of \cite{CLY}, we can prove the long time stability for the KAM tori obtained in Theorem \ref{t1}. The measure estimate can be done as in Theorem \ref{t1}.
\qed

%
%


\newcommand{\etalchar}[1]{$^{#1}$}

\end{document}